\documentclass{preprint}
\textwidth=\paperwidth
\evensidemargin=\oddsidemargin
\topmargin=-\headheight
\makeatletter
\newcommand\footnotetext@\relax
\let\footnotetext@\@footnotetext
\renewcommand{\thanks}[1]{%
  \unskip\protected@xdef\@thefnmark{}%
  \protected@xdef\@thanks{\@thanks\protect\footnotetext@{#1}}}
\newcommand{\MSC}[2][2010]{%
  \unskip\protected@xdef\@thefnmark{}%
  \protect\footnotetext@{\kern-1.8em{\slshape #1 MSC:}\/ #2.}}
\newcommand{\keywords}[1]{%
  \unskip\protected@xdef\@thefnmark{}%
  \protect\footnotetext@{\kern-1.8em{\slshape Keywords:}\/ #1.}}
\newcommand{\address}[1]{\\[1.2ex]%
  \parbox{.96\textwidth}{\centering\footnotesize\slshape #1}}
\gdef\@date{}
\makeatother
\usepackage{preprint}
\entrymodifiers={!!<.em,.6ex>+}
\renewcommand{\bydef}[1][=]{\overset{\mathit{def}}#1}
\renewcommand{\GL}{\operatorname{GL}}
\renewcommand{\GO}{\operatorname{O}}
\renewcommand{\SO}{\operatorname{SO}}
\makeatletter%
\newcommand\m@themphpunct[1]{\mathchoice{%
  \mathord{#1}\mkern5mu}%
 {\mathord{#1}\allowbreak\mskip\thickmuskip}{#1}{#1}}
\renewcommand{\emphpunct}[2][3000]{%
  \ifmmode
    \m@themphpunct{#2}%
  \else
    #2\spacefactor#1{}%
  \fi}
\newcommand\sidetext@{%
  \@ifnextchar(\sidetext@p
 {\@ifnextchar[\sidetext@b
  \relax}}
\newcommand\@sidetext{%
  \mskip.mu\@plus5mu\sidetext@}
\renewcommand{\sidetext}{%
  \@ifstar\sidetext@\@sidetext}
\newcommand{\longapproxto}{\overset{\textstyle\approx}{\sm@shrelto.5ex\longto}}
\newcommand{\longsimto}{\overset{\textstyle\sim}{\sm@shrelto.3ex\longto}}
\makeatother%
\theoremstyle{definition}
  \newtheorem{exmp}[stmt]{Example}
\newtheoremstyle{problem}{\topsep}{\topsep}{\itshape}%
  {}{\bfseries}{:}{.5em}{}
\theoremstyle{problem}
  \newtheorem*{PROBLEM}{Extension Problem}
\DeclareMathOperator{\Fr}{Fr}
\newcommand{\tto}{\rightrightarrows}
\DeclareMathOperator{\Mcon}{\mathit{Mcon}}
\DeclareMathOperator{\Rep}{\mathit{Rep}}
\DeclareMathOperator{\Conn}{\mathit{Conn}}
\DeclareMathOperator{\Psr}{\mathit{Psr}}
\usepackage{accents}
\newcommand{\nef}[1]{\accentset{\ast}{#1}}
\newcommand{\Mfd}{\mathit{Mfd}}
\newcommand{\Set}{\mathit{Set}}
\newcommand{\op}{\mathrm{op}}
\DeclareMathOperator{\Hilb}{Hilb}
\DeclareMathOperator{\Pos}{Pos}
\DeclareMathOperator{\Sym}{Sym}
\newcommand{\idc}[1]{\accentset{\circ}{#1}}
\newcommand{\tor}[1]{#1^\mathrm{tor}}
\newcommand{\free}[1]{#1^\mathrm{free}}
\newcommand{\inv}{\mathit{inv}}
\newcommand{\hol}{\mathit{hol}}
\DeclareMathOperator{\Skew}{Skew}
\title{Regular Cartan groupoids %
    and longitudinal representations%
  \MSC{Primary 22A22; %
    Secondary 53C05, 53C10, 55S40}
  \keywords{Lie groupoid, proper groupoid, regular groupoid, %
    multiplicative connection, longitudinal representation, %
    obstruction}
}%
\author{Giorgio Trentinaglia%
  \address{Center for Mathematical Analysis, Geometry and Dynamical Systems, %
    Instituto Superior T\'ecnico, Universidade de Lisboa, %
    Av.~Rovisco Pais, 1049-001 Lisboa, Portugal}
  \thanks{The author acknowledges support %
    from the Portuguese Foundation for Science and Technology %
    (Fun\-da\-\c c\~ao pa\-ra a Ci\-\^en\-cia e a Tec\-no\-lo\-gia) %
    through the grants %
    SFRH/\allowbreak BPD/\allowbreak 81810/\allowbreak 2011 and %
    UID/\allowbreak MAT/\allowbreak 04459/\allowbreak 2013.}
}%
\begin{document}
\maketitle

\begin{abstract} With the intent of laying the groundwork for a program that aims at explicitly describing the space of Cartan (i.e.~multiplicative) connections on a general proper Lie groupoid, we begin to investigate the space of such connections in the regular case. We point out that there is a close relationship between Cartan connections on a proper regular groupoid and representations of the groupoid on its own longitudinal bundle (i.e.~on the vector distribution tangent to its orbits). This observation enables us to reduce the original problem to a simpler one. We carry out a prospective study of the latter problem, and apply the resulting analysis to produce a number of examples in rank two which serve to illustrate the diversity of the possible obstructions to the existence of multiplicative connections. \end{abstract}

\tableofcontents

\section*{Introduction}

In order to give the reader an idea of the subject of this paper probably the best way is to start with an example. Suppose $G\to\GL(n,\R)$ is a Lie group mapping to the general linear group of degree $n$ over the reals. Recall that a \emph{$G$\mdash structure} on an $n$\mdash dimensional manifold $M$ is a reduction of the general linear frame bundle $\Fr(M)$ from $\GL(n,\R)$ to $G$ i.e.~an equivariant fi\-ber-bun\-dle map $P\to\Fr(M)$ from some right principal $G$\mdash bundle $P\to M$. This classical notion~\cite{Chern,Koba,Stern} admits a purely groupoid theoretic counterpart, namely, the notion of a \emph{tangent representation} $\rho:\varGamma\to\GL(TM)$ of a {\em transitive}\/ Lie groupoid $\varGamma\tto M$. In one direction, given a $G$\mdash structure $\sigma:P\to\Fr(M)$, you obtain a tangent representation of $\varGamma=(P\times P)/G$\textemdash the \emph{gauge groupoid} associated with the principal $G$\mdash bundle $P\to M$~\cite{MM}\textemdash upon setting \[%
	\rho(q\otimes p)=\sigma(q)\circ\sigma(p)^{-1}:T_xM\simto\R^n\simto T_yM
\] for all $p\in P_x\emphpunct, q\in P_y$. Conversely, up to isomorphism, every tangent representation of a transitive Lie groupoid over $M$ arises in this way from some $G$\mdash structure on $M$, for a suitable (essentially unique) choice of $G$. Observe that the case $G$~compact corresponds to the case $\varGamma$~proper. Next, in the study of $G$\mdash structures, it is customary to assume that the underlying principal $G$\mdash bundles say $\pi:P\to M$ be endowed with suitable principal (i.e.~equivariant) Ehresmann connections $h:\pi^*TM\to TP$. Now, given any tangent representation $\rho:\varGamma\to\GL(TM)$ of the gauge groupoid $\varGamma=(P\times P)/G$, you obtain a \emph{multiplicative connection} $\eta:s^*TM\to T\varGamma$ on $\varGamma$ whose \emph{effect} $t_*\circ\eta:s^*TM\simto t^*TM$ coincides with $\rho$~\cite{2016a}, when you set \[%
	\eta_{q\otimes p}=\pr_*\circ(h_q\circ\rho(q\otimes p),h_p):T_xM\to T_qP\oplus T_pP=T_{(q,p)}(P\times P)\onto T_{q\otimes p}((P\times P)/G).
\] Conversely any multiplicative connection $\eta$ on $\varGamma$ with effect $\rho$ will give rise to an essentially unique Ehresmann connection $h$ on $\pi:P\to M$ from which $\eta$ can be recovered by means of the preceding formula. In principle, every question one might ask about $G$\mdash structures could be translated into a question about tangent representations of, or multiplicative connections on, transitive Lie groupoids.

The above picture generalizes naturally to the intransitive case. Namely, let our transitive Lie groupoid $\varGamma\tto M$ be now an arbitrary regular groupoid. Any multiplicative connection on $\varGamma$ will determine a representation of $\varGamma$ on $TM$, the effect of the connection, obtained by composing the connection horizontal lift $s^*TM\to T\varGamma$ with the tangent groupoid target map $t_*:T\varGamma\to t^*TM$. This representation will necessarily carry the \emph{longitudinal bundle} $\varLambda$ of $\varGamma$\textemdash the subbundle \[%
	\varLambda=\bigcup_{x\in M}T_x\incl_x\left(T_xO_x\right)\subset TM
\] consisting of all those vectors that are tangent to the (immersed) $\varGamma$~orbits \(\incl_x:O_x\hookto M\)\textemdash into itself, thus giving rise to the \emph{longitudinal effect} of the connection, an example of what we call a \emph{longitudinal representation} $\varGamma\to\GL(\varLambda)$. As a matter of terminology, we shall refer to the rank of the smooth vector bundle $\varLambda$ also as the rank of the regular Lie groupoid $\varGamma$. In the transitive case, `longitudinal' is the same as `tangent', and, as we have seen in the preceding paragraph, every tangent representation is the effect of some multiplicative connection. In the more general intransitive (regular) case, there may be longitudinal representations which do not come from any multiplicative connections. However, if one restricts attention to the proper (regular) case, one can by an easy averaging argument still show that every longitudinal representation is the longitudinal effect of some (in general, non-unique) multiplicative connection; this relationship between multiplicative connections and longitudinal representations in the proper case is reviewed in \S\ref{sec:1} below. Multiplicative connections on Lie groupoids are also known as \emph{Cartan connections}~\cite{Bla16,CSS}; accordingly, by a \emph{Cartan groupoid} we mean a Lie groupoid endowed with a multiplicative connection.

Having introduced our main objects of interest\textemdash those referred to in the title\textemdash we turn to the motivations. The present study is centered around the following%

\begin{PROBLEM} Let\/ $\varGamma\tto M$ be a Lie groupoid which is regular and proper. Let\/ $U$ and\/ $V$ be invariant open subsets of\/ $M$ such that\/ $\overline U\subset V$, and let\/ $\varPhi$ be a multiplicative connection on\/ $\varGamma\mathbin|V\tto V$. What are the obstructions to extending\/ $\varPhi\mathbin|U$ to a (globally defined) multiplicative connection on\/ $\varGamma\tto M$? \end{PROBLEM}

\noindent In the light of our previous considerations the reader might suspect the above problem to be equivalent to an analogous problem only involving longitudinal representations; this is in fact the case, as we are going to explain at the beginning of \S\ref{sec:1}. Notice that already in the transitive case, the problem (with $U=V=\emptyset$) is quite interesting from a geometric viewpoint, and also quite non-triv\-i\-al. Indeed, translating back into the language of $G$\mdash structures, it corresponds to the problem of finding the obstructions to the existence of arbitrary $G$\mdash structures on manifolds~\cite{Chern}; for instance, Cartan connections on a pair groupoid $M\times M\tto M$ (which is proper, and source proper for compact $M$) are in one-to-one correspondence with \emph{absolute parallelisms} on $M$ i.e.~vector bundle trivializations of $TM$. For more information about the significance of the general extension problem\textemdash which in fact subsumes a number of at first sight unrelated obstruction problems in geometry\textemdash we refer the reader to the introduction of \cite{2016a}. (The material presented here was originally intended to appear in \cite{2016a}. Reasons of size, plus the fact that this material was essentially self-con\-tained, prompted us to write a separate paper.) We mention that the limitation coming from the regularity hypothesis is only apparent, as demonstrated in \cite{2016a}.

The contribution of this article to the study of the above-men\-tioned problem is two-fold. First, we carry out a conceptual analysis of the problem, showing that it can be rephrased as a standard problem in equivariant (orbifold) obstruction theory\textemdash a circumstance which opens the possibility of adapting existing techniques~\cite{Bred,tomD} to its study. We then proceed to describe some general technical tools that to a certain extent enable one to simplify the formulation of the problem. All this is done in Sections~\ref{sec:1}~through \ref{sec:4}. Second, in \S\ref{sec:5} and \S\ref{sec:6}, we put our analysis at work and produce examples (although probably we should rather call them results) featuring complete descriptions of the obstructions to the existence of extensions for specific classes of (proper regular) Lie groupoids. Altogether these examples, which we expect to serve in the future as aids to intuition in the elaboration of a comprehensive obstruction theory for Cartan connections, provide a fine although incomplete picture of the rank-two case.

In what follows, we shall stick to the overall conventions of \cite{2016a}. We have tried to make our exposition as self-con\-tained as possible; hopefully the reader will never need to refer back to \cite{2016a}. For any Lie groupoid $\varGamma\tto M$, we let $\Mcon(\varGamma)$ denote the space of multiplicative connections on $\varGamma$ (which could be empty). Moreover, for any smooth vector bundle $E$ over $M$, we let $\Rep(\varGamma;E)$ denote the space of representations of $\varGamma$ on $E$.


\section{The orbibundle of longitudinal representations}\label{sec:1}

Throughout this section $\varGamma\tto M$ will be a Lie groupoid that is regular and proper. As before we shall refer to the common dimension of its orbits as its \emph{rank}, and let $\varLambda\subset TM$ denote its \emph{longitudinal bundle}.

Let us begin with an observation\emphpunct: \em The mapping
\begin{equation}
\label{eqn:14A.6.1}
	\lambda_\varLambda:\Conn(\varGamma)\longto\Psr(\varGamma;\varLambda),
\quad	H\mapsto\lambda^H_\varLambda
\end{equation}
which to any groupoid connection\/ $H$ on\/ $\varGamma$ associates the pseudo-rep\-re\-sen\-ta\-tion of\/ $\varGamma$ on\/ $\varLambda$ induced upon restriction by\/ $\lambda^H\in\Psr(\varGamma;TM)$ (the effect of\/ $H$) is a trivializable affine fibration. \em Indeed, if one fixes a Riemannian metric $\mu$ on $M$, an arbitrary groupoid connection $H_0\in\Conn(\varGamma)$ and a splitting $\xi:t^*\varLambda\hookto\ker Ts\subset T\varGamma$ to the tangent target mapping $t_*:\ker Ts\onto t^*\varLambda$, the correspondence that to any pseudo-rep\-re\-sen\-ta\-tion $\alpha\in\Psr(\varGamma;\varLambda)$ associates the connection $\varPsi(\alpha)=\varPsi^{\mu,H_0,\xi}(\alpha)\in\Conn(\varGamma)$ characterized by the identity
\begin{equation}
\label{eqn:14A.6.2}
	\eta^{\varPsi(\alpha)}=\eta^{H_0}-\xi\circ(\lambda^{H_0}-%
\begin{pmatrix}
	\alpha & 0
\\	0      & \nu
\end{pmatrix}),
\end{equation}
where $\nu:\varGamma\to\GL(\varLambda^\bot)$ denotes the unique representation that corresponds to the intrinsic \emph{infinitesimal effect} of $\varGamma$ on the transversal bundle $TM/\varLambda\cong\varLambda^\bot$ (see \cite[\S 1]{Tre6}), gives a global section to $\lambda_\varLambda$,
\begin{equation}
\label{eqn:14A.6.3}
	\varPsi=\varPsi^{\mu,H_0,\xi}:\Psr(\varGamma;\varLambda)\longto\Conn(\varGamma)
\end{equation}
which carries representations to effective connections. This global section determines a choice of affine origin in each fiber, and thus enables one to define a global trivializing chart for $\lambda_\varLambda$;
\begin{equation}
\label{eqn:14A.6.4}
	\tau^\varPsi:\Conn(\varGamma)\approxto\Psr(\varGamma;\varLambda)\times\Gamma^\infty\bigl(\varGamma;L_{s^*\varLambda,\ker t_*}(s^*TM,\ker Ts)\bigr)
\end{equation}
the notation $L_{E',F'}(E,F)$ here refers to the smooth vector subbundle of $L(E,F)$ consisting of all those linear maps that carry a given subbundle $E'\subset E$ into another subbundle $F'\subset F$. Every change of charts will be a fiberwise affine transformation of vector spaces\textemdash a translation in each fiber, in fact.

\begin{stmt}\label{stmt:10/17} Let\/ $U$ and\/ $V$ be invariant open subsets of\/ $M$ such that\/ $\overline U\subset V$. Let\/ $\varTheta\in\Mcon(\varGamma\mathbin|V)$ be a multiplicative connection on\/ $\varGamma\mathbin|V\tto V$, and let\/ $\lambda^\varTheta_\varLambda\in\Rep(\varGamma\mathbin|V;\varLambda\mathbin|V)$ denote the longitudinal effect of\/ $\varTheta$. Then\/ $\varTheta\mathbin|U\in\Mcon(\varGamma\mathbin|U)$ can be extended to a multiplicative connection defined on all of\/ $\varGamma\tto M$ if, and only if, $\lambda^\varTheta_\varLambda\mathbin|U\in\Rep(\varGamma\mathbin|U;\varLambda\mathbin|U)$ can be extended to some global longitudinal representation\/ $\varGamma\to\GL(\varLambda)$. \end{stmt}

\begin{proof} The argument makes use of the averaging operator introduced with \cite{2016a}; in fact, the truth of our statement depends in an essential way on the hypothesis that our groupoid $\varGamma$ be proper, compare \cite[Counterexample~2.9]{2016a}. Suppose that it is possible to prolong $\lambda^\varTheta_\varLambda\mathbin|U$ to some global longitudinal representation say $\alpha\in\Rep(\varGamma;\varLambda)$. Then by using a global trivializing chart of the form \eqref{eqn:14A.6.4} for the trivializable affine fibration \eqref{eqn:14A.6.1} we may lift $\alpha$ to a groupoid connection $\varPhi\in\Conn(\varGamma)$ so that $\varPhi\mathbin|U=\varTheta\mathbin|U$. The \emph{multiplicative average} $\hat\varPhi$ of $\varPhi$\textemdash whose horizontal lift at any $g\in\varGamma$ recall is given by the integral expression
\begin{equation*}
	\eta^{\hat\varPhi}_g=\integral_{k\in\varGamma(-,sg)}[\eta^\varPhi_{gk}\cdot(\eta^\varPhi_k)^{-1}]\circ(\lambda^\varPhi_k)^{-1}dk
\end{equation*}
(depending on the preliminary choice of some left invariant normalized Haar system on $\varGamma$)\textemdash will be a multiplicative connection on $\varGamma$ satisfying $\hat\varPhi\mathbin|U=\varTheta\mathbin|U$, as desired. \end{proof}

\em Motivated by the preceding statement, we shall henceforth focus on the analysis of the longitudinal representations of\/ $\varGamma$.\em

\begin{npar}\label{npar:14A.6.1} For each base point $x\in M$, let $\nef\varGamma^x_x$ denote the \emph{ineffective part} of the isotropy group $\varGamma^x_x$ (cf.~\cite[Section~1]{Tre6}), consisting by definition of all those arrows whose intrinsic infinitesimal effects on the transversal tangent space $T_xM/\varLambda_x$ are trivial. As $x$ ranges over $M$, the closed normal subgroups $\nef\varGamma^x_x\subset\varGamma^x_x$ trace out a \emph{kernel} within $\varGamma$, that is, a totally isotropic subgroupoid of $\varGamma$ that contains all the units and is normal. Because of the regularity and properness of $\varGamma$, this kernel $\nef\varGamma$ must in fact be a closed, smooth subgroupoid of $\varGamma$ of constant dimension (hence in particular a Lie groupoid $\nef\varGamma\tto M$ in its own right). The quotient
\begin{equation}
\label{eqn:14A.6.5}
	P:=\varGamma/\nef\varGamma\tto M
\end{equation}
makes sense therefore as a {\em Lie}\/ groupoid (cf.~\cite[Appendix~A]{Tre6}), and as such, it is a proper foliation groupoid.

There is a well-de\-fined, canonical (left) Lie-group\-oid action
\begin{equation}
\label{eqn:14A.6.6}
\xymatrix@C=.em@R=3ex{%
 \varGamma
 \ar@<+.4ex>[dr]\ar@<-.4ex>[dr]
 &	\circlearrowright
	&	P\makebox[.em][l]{,\qquad $\smash{\varGamma\ftimes stP}\to P$}
		\ar[dl]^t
\\ &	M}
\end{equation}
given by $g\cdot[h]=[gh]$ (arrow composition), which determines a corresponding (smooth) translation groupoid, hereafter denoted
\begin{equation}
\label{eqn:14A.6.7}
	\varPi:=\varGamma\ltimes P\tto P.
\end{equation}
The orbits of this action (that is to say, the $\varPi$~orbits) coincide with the fibers $P^x$ of the source map $s:P\to M$. The submersions $t^x:P^x\twoheadto O_x$ are covering projections with finite fibers $P^x_x$, so one has $\dim P^x=\dim O_x$ for all $x$. Accordingly, the Lie groupoid $\varPi\tto P$ is regular, of the same rank $q$ as $\varGamma\tto M$ (even though in general $\dim P\geqq\dim M$, with equality only for $q=0$). It is moreover proper, as it is easy to check. For each $h\in\varGamma^x$, one has a canonical identification \(T_{[h]}P^x=T_{th}O_x\) of vector spaces, and another one \(\varPi^{[h]}_{[h]}=\nef\varGamma^{th}_{th}\) of Lie groups. The isotropy kernel of $\varPi$ is consequently a locally trivial bundle of compact Lie groups (a \emph{proper Lie bundle}), equivalently, $\varPi$ is \emph{principal} (every isotropic arrow has trivial intrinsic infinitesimal effect). Intuitively, $\varPi$ should be thought of as a suitable ``resolution'' or ``desingularization'' of $\varGamma$:
\begin{equation}
\label{eqn:14A.6.8}
 \begin{split}
\xymatrix@R=3.3ex{%
 \varPi
 \ar@<+.4ex>[d]\ar@<-.4ex>[d]
 \ar@{->>}[r]^\pr
 &	\varGamma
	\ar@<+.4ex>[d]\ar@<-.4ex>[d]
\\ P
 \ar@{->>}[r]^t
 &	M
}\end{split}
\end{equation} \end{npar}

\begin{npar}\label{npar:14A.6.2} We shall let \(\varPi^x\tto P^x\) denote the (proper, transitive) Lie groupoid $\varPi\mathbin|P^x\tto P^x$, and
\begin{equation}
\label{eqn:14A.6.9}
	\mathscr R^\varGamma_x\bydef\Rep(\varPi^x;TP^x)
\end{equation}
the space of its smooth representations on $TP^x$. We shall regard the union of all these spaces as a (set-the\-o\-ret\-ic) fiber bundle $\mathscr R^\varGamma=\bigcup_{x\in M}\mathscr R^\varGamma_x$ over $M$.

We have a canonical groupoid action
\begin{equation}
\label{eqn:14A.6.10}
\xymatrix@C=.em@R=3ex{%
 P
 \ar@<+.4ex>[dr]\ar@<-.4ex>[dr]
 &	\circlearrowright
	&	\mathscr R^\varGamma\makebox[.em][l]{,\qquad $\smash{P\ftimes s\pr\mathscr R^\varGamma}\to\mathscr R^\varGamma$}
		\ar[dl]^\pr
\\ &	M}
\end{equation}
arising as follows. Any arrow $k\in\varGamma$ determines a Lie-group\-oid isomorphism
\begin{equation}
\label{eqn:14A.6.11}
 \begin{split}
\xymatrix@R=3.3ex{%
 \varPi^{tk}
 \ar@<+.4ex>[d]\ar@<-.4ex>[d]
 \ar[r]^-k_-*{\approx}
 &	\varPi^{sk}
	\ar@<+.4ex>[d]\ar@<-.4ex>[d]
\\ P^{tk}
 \ar[r]^-k_-*{\approx}
 &	P^{sk}
}\end{split}
\end{equation}
acting on base points resp.~arrows by the rule $k([h])=[hk]$ resp.~$k(g,[h])=(g,[hk])$. This isomorphism is the same for all representatives $k$ of a given class $[k]\in P$. To each arrow $[k]\in P$ and representation $\rho\in\mathscr R^\varGamma_{sk}=\Rep(\varPi^{sk};TP^{sk})$ the action \eqref{eqn:14A.6.10} then assigns the representation%
\begin{subequations}
\label{eqn:10/19}
\begin{equation}
\label{eqn:14A.6.12}
	[k]\cdot\rho\bydef(t^*k_*)^{-1}(k^*\rho)(s^*k_*)\in\mathscr R^\varGamma_{tk}=\Rep(\varPi^{tk};TP^{tk}),
\end{equation}
where $k_*:TP^{tk}\simto k^*TP^{sk}$ denotes the vec\-tor-bun\-dle isomorphism that corresponds to $Tk:TP^{tk}\simto TP^{sk}$. Explicitly,
\begin{equation}
\label{eqn:14A.6.13}
	([k]\cdot\rho)(g,[h])\bydef(T_{[gh]}k)^{-1}\circ\rho(g,[hk])\circ T_{[h]}k:T_{[h]}P^{tk}\to T_{[gh]}P^{tk}.
\end{equation}
\end{subequations} \end{npar}

\begin{npar}\label{npar:14A.6.3} Next, we want to endow the total space $\mathscr R^\varGamma$ of our abstract fiber bundle $\mathscr R^\varGamma\to M$ with a suitable ``smooth'' structure (in the ``diffeological'' sense). This structure is going to make $\mathscr R^\varGamma\to M$ into a ``smooth'' fibration and \eqref{eqn:14A.6.10} into a ``smooth'' groupoid action (the fiber product $P\ftimes s\pr\mathscr R^\varGamma$ being given the ``smooth'' structure relative to which a map $S\to P\ftimes s\pr\mathscr R^\varGamma$ is ``smooth'' whenever so are its components $S\to P$, $S\to\mathscr R^\varGamma$).

Let $\Mfd$ stand for the category of smooth manifolds. For any given ``space'' $\mathscr X$, let $\Set(-,\mathscr X):\Mfd^\op\to\Set$ denote the presheaf that to every smooth manifold $S$ assigns the set of all maps from $S$ to $\mathscr X$. By a \emph{\(C^\infty\)\mdash structure} on $\mathscr X$, we shall mean a subsheaf $\mathcal S\subset\Set(-,\mathscr X)$ which contains all the constant maps\textemdash the sheaf property being with respect to arbitrary open covers. We shall refer to the pair $\mathscr X=(\mathscr X,\mathcal S)$ as a \emph{\(C^\infty\)\mdash space}. By a \emph{\(C^\infty\) mapping} $f:\mathscr X\to\mathscr Y$ between two \(C^\infty\)\mdash spaces we shall mean one which by forward composition gives rise to a natural transformation between the corresponding \(C^\infty\)\mdash structures.

In order that a map $S\xto\rho\mathscr R^\varGamma$ may be declared to be ``smooth'', it is first of all required that the composition $S\xto\rho\mathscr R^\varGamma\to M$ be smooth. Whenever that is the case, the pullback $P\times_MS\xto{\pr_P}P$ will exist as a smooth manifold, and the groupoid $\varPi\times_MS\overset{s\times\id}{\underset{t\times\id}\tto}P\times_MS$ [with composition law $(g',[gh];s)(g,[h];s)=(g'g,[h];s)$] will make sense as a Lie groupoid. We then declare $S\xto\rho\mathscr R^\varGamma$ to be ``smooth'' if, and only if, so is the representation
\begin{equation*}
	\varPi\times_MS\longto\GL(\pr_P^*L),\quad(g,[h];s)\mapsto\rho(s)(g,[h]),
\end{equation*}
where $L\to P$ denotes the longitudinal bundle of the regular groupoid $\varPi\tto P$.

It follows at once from this definition that the fi\-ber-bun\-dle projection $\mathscr R^\varGamma\to M$ is a \(C^\infty\) mapping. Also, it is not hard to see that the operation $P\times_M\mathscr R^\varGamma\to\mathscr R^\varGamma$ underlying the groupoid action \eqref{eqn:14A.6.10} is a \(C^\infty\) mapping. Our definition also provides $\mathscr R^\varGamma\to M$ with a notion of local \(C^\infty\) section; it is immediate to recognize that a local \(C^\infty\) section $U\xto\rho\mathscr R^\varGamma$ is the same thing as a smooth representation $\rho:\varPi\mathbin|P^U\to\GL(L\mathbin|P^U)$, where $P^U\subset P$ is the open set $\{[h]:sh\in U\}$. \end{npar}

\begin{stmt}\label{stmt:14A.6.4} The correspondence\/ $\Rep(\varGamma;\varLambda)\to\Gamma(M;\mathscr R^\varGamma)$ that to every longitudinal representation\/ $\alpha:\varGamma\to\GL(\varLambda)$ associates the global section\/ $M\ni x\mapsto\rho^\alpha(x)\in\mathscr R^\varGamma_x$, where
\begin{equation}
\label{eqn:14A.6.15}
	\rho^\alpha(x)(g,[h])\bydef(T_{[gh]}t^x)^{-1}\circ\alpha(g)\circ T_{[h]}t^x:T_{[h]}P^x\to T_{[gh]}P^x
\end{equation}
($t^x$ being the covering projection\/ $P^x\twoheadto O_x$), is one-to-one onto the equivariant\/ \(C^\infty\) sections of\/ $\mathscr R^\varGamma$.\qed \end{stmt}

\begin{npar}\label{npar:14A.6.5} The preceding statement may be paraphrased informally by saying that the longitudinal representations of $\varGamma$ are in one-to-one correspondence with the ``smooth'' sections of the ``orbibundle'' \([\mathscr R^\varGamma/P]\to[M/P]\). In particular, let $U,V\subset M$ be as in \ref{stmt:10/17}. Suppose that $\alpha\in\Rep(\varGamma\mathbin|V;\varLambda\mathbin|V)$ is a longitudinal representation of $\varGamma\mathbin|V\tto V$, and let $\rho^\alpha\in\Gamma^\infty(V;\mathscr R^\varGamma)^P$ denote the equivariant local $C^\infty$ section of $\mathscr R^\varGamma$ defined over $V$ which, by \ref{stmt:14A.6.4}, corresponds to $\alpha$. Then, it will be possible to prolong $\alpha\mathbin|U\in\Rep(\varGamma\mathbin|U;\varLambda\mathbin|U)$ to a global longitudinal representation of $\varGamma$ if, and only if, there exists some extension of $\rho^\alpha\mathbin|U\in\Gamma^\infty(U;\mathscr R^\varGamma)^P$ to a global equivariant $C^\infty$ section of $\mathscr R^\varGamma$. Thus, within the geometric setup that we are developing for the analysis of longitudinal representations, the extension problem that we want to solve can be reformulated as a special case of a standard problem in equivariant obstruction theory. We shall see in \S\ref{sec:3} that under relatively mild assumptions on $\varGamma$ such as, for example, source properness, the ``orbibundle'' \([\mathscr R^\varGamma/P]\to[M/P]\) is ``locally trivial'' in a sense to be clarified. \end{npar}

\subsection*{The orthogonalization trick}

We proceed to expunge a few irrelevant ``degrees of freedom'' from our analysis of the longitudinal representations of $\varGamma$.

Let $\phi$ be an arbitrary metric on the longitudinal bundle $L\subset TP$ of the regular groupoid $\varPi\tto P$, and let $\phi^x$ stand for the Riemannian metric induced on $P^x$ along the canonical identification of vector bundles $TP^x=L\mathbin|P^x$. Then, let $\mathscr R^{\varGamma,\phi}\subset\mathscr R^\varGamma$ denote the subbundle whose fiber at $x$ is the space of $\phi^x$~orthogonal tangent representations of $\varPi^x\tto P^x$,
\begin{equation}
\label{eqn:14A.6.16}
	\mathscr R^{\varGamma,\phi}_x\bydef\Rep(\varPi^x;TP^x,\phi^x).
\end{equation}
We are going to construct a $P$~equivariant, strong deformation retraction of fiber bundles
\begin{equation}
\label{eqn:10/21}
	\{-\}_t:\mathscr R^\varGamma\rightsquigarrow\mathscr R^{\varGamma,\phi}
\end{equation}
which is also $C^\infty$ as a map $\R\times\mathscr R^\varGamma\to\mathscr R^\varGamma$, in the sense that its composition with every $C^\infty$ map $S\to\R\times\mathscr R^\varGamma$ is a $C^\infty$ map. Of course, in order for this to be possible we need $\mathscr R^{\varGamma,\phi}$ to be a $P$~invariant subbundle of $\mathscr R^\varGamma$ i.e.~to be carried into itself by the $P$~action \eqref{eqn:14A.6.10}. In order for that to be the case, we need the metric $\phi$ to be invariant under the $P$~action, in other words, we need \[%
	(P^{tk},\phi^{tk})\xto[\textstyle\approx]k(P^{sk},\phi^{sk})
\] to be a Riemannian isometry for each arrow $k\in\varGamma$. Now the inverse image of any metric on $\varLambda$ along the vec\-tor-bun\-dle projection $L=t^*\varLambda\to\varLambda$ will be one such $P$~invariant metric on $L$. So $L$ does admit $P$~invariant metrics trivially. Let us once and for all fix one such metric $\phi$.

Since by hypothesis $\varGamma\tto M$ is proper, so will be $\varPi\tto P$ and hence a~fortiori each one of the transitive Lie groupoids $\varPi^x\tto P^x$. Let us fix an arbitrary {\em right invariant}, normalized, Haar system $\nu$ on $\varGamma\tto M$. This will give rise to similar Haar systems on $\varPi\tto P$ and $\varPi^x\tto P^x$, which we may denote $\nu$ as well. Our first move is to use the Haar system on $\varPi^x\tto P^x$ to average the metric $\phi^x$ on $TP^x$ relative to any given representation $\rho\in\Rep(\varPi^x;TP^x)$. We do this by writing down the following Haar integral expression, which depends parametrically on $[h]\in P^x$, where $v,w$ stand for arbitrary tangent vectors in $T_{[h]}P^x$:
\begin{equation*}
	\hat\phi^\rho_{[h]}(v,w):=\integral_{sg=th}\phi^x_{[gh]}\bigl(\rho(g,[h])v,\rho(g,[h])w\bigr)d\nu^{th}(g).
\end{equation*}
By construction, $\rho\in\Rep(\varPi^x;TP^x,\hat\phi^\rho)$. Moreover, because of the right invariance of $\nu$,
\begin{equation}
\label{eqn:14A.6.19}
	\hat\phi^{k\cdot\rho}=k^*\hat\phi^\rho.
\end{equation}

Next, let \(H^\rho_{[h]}P^\rho_{[h]}\) denote the \emph{polar decomposition} (see e.g.~\cite[\S VII.2]{Lang}) of the identity map
\begin{equation*}
	\id:(T_{[h]}P^x,\phi^x_{[h]})\simto(T_{[h]}P^x,\hat\phi^\rho_{[h]});
\end{equation*}
by definition, $H^\rho_{[h]}\in\Hilb(T_{[h]}P^x,\phi^x_{[h]};T_{[h]}P^x,\hat\phi^\rho_{[h]})$ and $P^\rho_{[h]}\in\Pos(T_{[h]}P^x,\phi^x_{[h]})$ are respectively an inner product preserving isomorphism and a symmetric positive definite automorphism. Since
\begin{equation*}
	\exp:\Sym(T_{[h]}P^x,\phi^x_{[h]})\approxto\Pos(T_{[h]}P^x,\phi^x_{[h]})
\end{equation*}
is a $C^\infty$\mdash isomorphism with domain a vector space (namely the space of all symmetric linear operators), there will be a canonical path $P^\rho_{t,[h]}:=\exp(t\log P^\rho_{[h]})$ deforming $P^\rho_{[h]}=P^\rho_{1,[h]}$ into $\id=P^\rho_{0,[h]}$ within $\Pos(T_{[h]}P^x,\phi^x_{[h]})$, and therefore also a canonical path
\begin{equation*}
	A^\rho_{t,[h]}:=H^\rho_{[h]}P^\rho_{t,[h]}
\end{equation*}
deforming $\id=A^\rho_{1,[h]}$ into $H^\rho_{[h]}=A^\rho_{0,[h]}$ through invertible linear operators. Then, setting
\begin{equation*}
	\rho_t(g,[h]):=(A^\rho_{t,[gh]})^{-1}\circ\rho(g,[h])\circ A^\rho_{t,[h]}
\end{equation*}
will yield a path $t\mapsto\rho_t$ within $\Rep(\varPi^x;TP^x)$ deforming $\rho=\rho_1$ into a representation $\rho_0$ which is orthogonal with respect to the original metric $\phi^x$. You can check that $(t,\rho)\mapsto\rho_t$ provides the desired deformation retraction, the equivariance of every fi\-ber-bun\-dle map $\rho\mapsto\rho_t$ being a consequence of the equation \eqref{eqn:14A.6.19}.

\begin{stmt}\label{stmt:10/21} Let\/ $U,V\subset M$ be as usual. Let\/ $\rho\in\Gamma^\infty(V;\mathscr R^\varGamma)^P$ be an equivariant local\/ $C^\infty$ section of\/ $\mathscr R^\varGamma$ with domain\/ $V$. Then, there exist vector bundle metrics\/ $\phi$ on\/ $\varLambda$ for which\/ $\rho\mathbin|U\in\Gamma^\infty(U;\mathscr R^{\varGamma,\phi})^P$. For any such metric\/ $\phi$, global equivariant\/ $C^\infty$ sections of\/ $\mathscr R^\varGamma$ extending\/ $\rho\mathbin|U$ exist if and only if there exist global equivariant\/ $C^\infty$ sections of\/ $\mathscr R^{\varGamma,\phi}$ with the same property. \end{stmt}

\begin{proof} Let $\alpha$ denote the unique element of $\Rep(\varGamma\mathbin|V;\varLambda\mathbin|V)$ for which, in the notations of \ref{stmt:14A.6.4}, $\rho=\rho^\alpha$. By properness, we can find some $\alpha$~invariant vector bundle metric on $\varLambda\mathbin|V$. Then, by a standard partition of unity argument, we can construct some vector bundle metric $\phi$ on $\varLambda$ with respect to which $\alpha\mathbin|U$ is an orthogonal representation. Evidently, for any such $\phi$ we must have $\rho\mathbin|U=\rho^{\alpha\mathbin|U}\in\Gamma^\infty(U;\mathscr R^{\varGamma,\phi})^P$. The orthogonalization trick can then be used to deform any global equivariant $C^\infty$ section of $\mathscr R^\varGamma$ which extends $\rho\mathbin|U$ into one of $\mathscr R^{\varGamma,\phi}$ with the same property. \end{proof}

\begin{exmp}\label{exmp:10/21*} Suppose that $\varGamma$ has rank zero i.e.~that every $\varGamma$~orbit is ze\-ro-di\-men\-sion\-al. Let $U\subset M$ be an arbitrary open set and let $\rho\in\Gamma^\infty(U;\mathscr R^\varGamma)$ be an arbitrary local $C^\infty$ section to $\mathscr R^\varGamma\to M$ defined over $U$. Then, $\rho$ admits a unique extension to an element of $\Gamma^\infty(M;\mathscr R^\varGamma)^P$; indeed, for every $x\in M$ we have $\dim P^x=\dim O_x=0$ and therefore $\mathscr R^\varGamma_x=\Rep(\varPi^x;TP^x)=\{\ast\}$. \end{exmp}

\begin{exmp}\label{exmp:10/21**} Suppose that $\varGamma$ has rank one i.e.~that its orbits are $1$\mdash dimensional. Suppose in addition that it is source connected. Let $U,V\subset M$ be invariant open sets such that $\overline U\subset V$. Let $\rho\in\Gamma^\infty(V;\mathscr R^\varGamma)^P$ be an arbitrary equivariant local $C^\infty$ section to $\mathscr R^\varGamma\to M$ with domain $V$. Then, it will always be possible to prolongate $\rho\mathbin|U$ to some element of $\Gamma^\infty(M;\mathscr R^\varGamma)^P$. Indeed, because of our assumptions, for any vector bundle metric $\phi$ on $\varLambda$ the $C^\infty$ fibration $\mathscr R^{\varGamma,\phi}\to M$ will be a 1:1 covering. Hence, on the account of \ref{stmt:10/21}, an extension of $\rho\mathbin|U$ can be obtained by choosing $\phi$ so that $\rho\mathbin|U\in\Gamma^\infty(U;\mathscr R^{\varGamma,\phi})^P$. \end{exmp}

\section{Equivariant contractions}\label{sec:2}

We now take a little time to look back into the constructs of the previous section from a slightly more abstract point of view, so as to highlight certain ideas implicit in the theory that lend themselves to application in situations more general than, for instance, the source-con\-nect\-ed, rank-one case considered in Example~\ref{exmp:10/21**}. After a brief elaboration on these ideas, we shall present a generalization of Example~\ref{exmp:10/21**}, namely Proposition~\ref{prop:14A.7.9} below, which eventually will enable us to say something about the rank-two case. Throughout this section, we shall be dealing with an arbitrary {\em proper}\/ Lie groupoid \(\varOmega\tto X\) along with an arbitrary \emph{\(\varOmega\)~equivariant\/ \(C^\infty\)\mdash fibration} \(\mathscr R\to X\), by which the following set of data should be understood: a \(C^\infty\)\mdash space \(\mathscr R\); a \(C^\infty\) mapping \(\mathscr R\to X\); a \(C^\infty\) (left) Lie-group\-oid action \(\varOmega\times_X\mathscr R\to\mathscr R\).

Let \(V\subset X\) be an invariant open set. By an \emph{equivariant contraction} of \(\mathscr R\to X\) over \(V\) we shall mean a one-pa\-ram\-e\-ter family of fiber preserving \(\varOmega\mathbin|V\tto V\)~equivariant maps \(\{-\}_t:\mathscr R\mathbin|V\to\mathscr R\mathbin|V\sidetext(t\in\R)\) which give rise to a \(C^\infty\) total mapping \(\R\times\mathscr R\mathbin|V\to\mathscr R\mathbin|V\) and which for every \(x\in V\) restrict for \(t=1\) to the identity on \(\mathscr R_x\) and for \(t=0\) to a retraction of \(\mathscr R_x\) onto a single point. We shall also require, as part of the definition, that the correspondence \(V\ni x\mapsto\bar\rho(x)\in\mathscr R_x\) which to each \(x\) assigns the constant value \(\bar\rho(x)\) of the \(t=0\)~contraction of \(\mathscr R_x\) provides an (equivariant) local \(C^\infty\) section to \(\mathscr R\to X\) over \(V\). (We find it convenient, and also correct from a formal point of view, to define {\em contractible}\/ spaces to be {\em non-emp\-ty}.)

\begin{lem}\label{lem:14A.7.1} Let\/ \(\{-\}^V_t:\R\times\mathscr R\mathbin|V\to\mathscr R\mathbin|V\) and\/ \(\{-\}^{V'}_t:\R\times\mathscr R\mathbin|V'\to\mathscr R\mathbin|V'\) be two equivariant contractions of\/ \(\mathscr R\to X\) over two invariant open subsets\/ \(V\) and\/ \(V'\) of\/ \(X\). Then, for any invariant open subsets\/ \(U\subset\overline U\subset V\) and\/ \(U'\subset\overline{U'}\subset V'\) it is possible to construct an equivariant contraction\/ \(\{-\}_t:\R\times\mathscr R\mathbin|U\cup U'\to\mathscr R\mathbin|U\cup U'\) that coincides with\/ \(\{-\}^V_t\) over\/ \(U\). \end{lem}

\begin{proof} For any invariant $C^\infty$ function $\varphi:X\to[0,1]$ with $\supp\varphi\subset V$ let $\{-\}^\varphi_t:\R\times\mathscr R\mathbin|U\cup U'\to\mathscr R\mathbin|U\cup U'$ denote the fiber preserving, equivariant, $C^\infty$ mapping given by \[%
	\{-\}^\varphi_{t,x}=%
\begin{cases}
	\{-\}^V_{\varphi(x)(t-1)+1,x} &\text{for $x\in U\cup(U'\cap V)$}
\\	\id                           &\text{for $x\in U'\smallsetminus\supp\varphi$}
\end{cases}
:\mathscr R_x\to\mathscr R_x
\] (the two alternative expressions agree on the intersection $U'\cap V\smallsetminus\supp\varphi$ so this is a good definition). Notice that the map $\{-\}^\varphi_{t,x}$ equals $\id$ for $t=1$ or when $x\notin\supp\varphi$ and that it equals $\{-\}^V_{t,x}$ when $\varphi(x)=1$. Symmetrically, for any invariant $C^\infty$ function $\varphi':X\to[0,1]$ with $\supp\varphi'\subset V'$ there will be a fi\-ber-pre\-serv\-ing and equivariant $C^\infty$ mapping $\{-\}^{\varphi'}_t:\R\times\mathscr R\mathbin|U\cup U'\to\mathscr R\mathbin|U\cup U'$ whose restriction to $\mathscr R_x$ is the identity for $t=1$ or when $x\notin\supp\varphi'$ and equals $\{-\}^{V'}_{t,x}$ when $\varphi'(x)=1$.

Now, pick $\varphi$ so that $\varphi=1$ in a neighborhood of $\overline U$, and then take $\varphi'$ with $\supp\varphi'\subset V'\smallsetminus\overline U$ so that $\varphi'=1$ in a neighborhood of $\overline{U'}\smallsetminus\varphi^{-1}(1)$; properness of $\varOmega\tto X$ guarantees that such a choice of invariant functions $\varphi$ and $\varphi'$ can be made. Then setting \[%
	\{\rho\}_t=\{\{\rho\}^\varphi_t\}^{\varphi'}_t
\] defines the desired equivariant contraction $\{-\}_t:\R\times\mathscr R\mathbin|U\cup U'\to\mathscr R\mathbin|U\cup U'$; indeed for every $x\in U\cup U'$ either $\varphi=1$ near $x$, in which case $\{-\}_{0,x}=\{-\}^{\varphi'}_{0,x}\circ\{-\}^V_{0,x}$, or $\varphi'=1$ near $x$, in which case $\{-\}_{0,x}=\{-\}^{V'}_{0,x}\circ\{-\}^\varphi_{0,x}$. \end{proof}

\begin{lem}\label{lem:14A.7.2} Let\/ \(\mathcal V\) be an open cover of\/ \(X\) by invariant open sets. Let\/ \(\bigl\{\{-\}^V_t:\R\times\mathscr R\mathbin|V\to\mathscr R\mathbin|V\bigr\}\) be any family of equivariant contractions of\/ \(\mathscr R\to X\) over the open sets\/ \(V\in\mathcal V\) of this cover. Then, it is possible to construct a global equivariant contraction of\/ \(\mathscr R\to X\). \end{lem}

\begin{proof} Choose any countable open cover of $X$ by invariant open sets $\{U_n\}$ subordinated to $\mathcal V$ in the sense that the closure $\overline{U_n\mkern-6mu}\mkern+6mu$ of each $U_n$ is contained in some $V\in\mathcal V$. Select $V_0\supset\overline{U_0\mkern-6mu}\mkern+6mu$ and set $\{-\}^{(0)}_t=\{-\}^{V_0}_t$. By induction on $n$ suppose that an equivariant contraction $\{-\}^{(n)}_t$ has been constructed over a suitable invariant open neighborhood $V_n$ of the union $\overline{U_0\mkern-6mu}\mkern+6mu\cup\overline{U_1\mkern-6mu}\mkern+6mu\cup\dotsb\cup\overline{U_n\mkern-6mu}\mkern+6mu$. Then, by the preceding lemma, for any choice of invariant open neighborhoods $U\subset\overline U\subset V_n$ of $\overline{U_0\mkern-6mu}\mkern+6mu\cup\overline{U_1\mkern-6mu}\mkern+6mu\cup\dotsb\cup\overline{U_n\mkern-6mu}\mkern+6mu$ and $U'\subset\overline{U'}\subset V'\in\mathcal V$ of $\overline{U_{n+1}\mkern-19mu}\mkern+19mu$ there will be some equivariant contraction $\{-\}^{(n+1)}_t$ over $V_{n+1}=U\cup U'$ that coincides with $\{-\}^{(n)}_t$ over $U$. \end{proof}

\begin{lem}\label{lem:14A.7.3} Let\/ \(i:S\to X\) be any smooth mapping completely transversal to the orbits of\/ \(\varOmega\tto X\). Suppose that the pullback\/ \(C^\infty\)\mdash fibration\/ \(i^*\mathscr R=S\times_X\mathscr R\to S\) admits a global\/ \(i^*\varOmega\tto S\)~equivariant contraction\/ \(\{-\}^S_t:\R\times i^*\mathscr R\to i^*\mathscr R\). Then, it is possible to construct a global\/ \(\varOmega\)~equivariant contraction of\/ \(\mathscr R\to X\). \end{lem}

\begin{proof} Given $\rho\in\mathscr R_x$, pick $h\in\varOmega$ with $sh=x$ and $th=i(s)\in i(S)$, and set \[%
	\{\rho\}_t=h^{-1}\cdot\pr\{(s,h\cdot\rho)\}^S_t.
\] On the account of the $i^*\varOmega$~equivariance of $\{-\}^S_t$, this definition will not depend on the choice of $h,s$. \end{proof}

We shall call \(\mathscr R\) \emph{locally transversely trivializable} if through each base-point \(x\) there exist slices \(i:S\hookto X\) over which the pullback \(C^\infty\)\mdash fibration \(i^*\mathscr R\to S\) admits equivariant trivializations of the form%
\begin{subequations}
\label{eqn:10/25}
\begin{equation}
\label{eqn:14A.7.1}
 \begin{split}
\xymatrix@C=5em@R=9ex{%
 S\times\mathscr R_x
 \ar[d]|(.42)\hole
 \ar[r]^(.5)*-<.3ex>{\approx}_(.5){C^\infty}
 \save[]+<-2.33em,-5ex>
	*+{\varOmega^x_x\ltimes S}
	\ar@{}[]|*{\circlearrowright}
	\ar@<+.4ex>[d]\ar@<-.4ex>[d]
	\ar[r]!/l1.1em/+<-2.33em,-5ex>^-*-<.6ex>{\sim}
 \restore
 &	i^*\mathscr R
	\ar[d]
	\save[]+<-2.33em,-5ex>
		*+{i^*\varOmega}
		\ar@{}[]|*{\circlearrowright}
		\ar@<+.4ex>[d]\ar@<-.4ex>[d]
	\restore
\\ S
 \ar@{=}[r]
 &	S
}\end{split}
\end{equation}
with \(\varOmega^x_x\ltimes S\tto S\) acting on \(S\times\mathscr R_x\to S\) by the rule
\begin{equation}
\label{eqn:14A.7.2}
	(k,s)\cdot(s,\rho)=(ks,k\cdot\rho).
\end{equation}
\end{subequations}

\begin{lem}\label{lem:14A.7.4} Suppose that\/ \(\mathscr R\) is locally transversely trivializable. Let\/ \(x\) be a base-point for which there exists an\/ \(\varOmega^x_x\)~equivariant\/ \(C^\infty\) contraction\/ \(\{-\}^x_t:\R\times\mathscr R_x\to\mathscr R_x\) of the fiber\/ \(\mathscr R_x\). Then, for some slice\/ $i:S\hookto X$ through\/ $x$, it will be possible to construct an\/ \(i^*\varOmega\tto S\)~equivariant contraction of\/ \(i^*\mathscr R\to S\). \end{lem}

\begin{proof} Setting $\{(s,\rho)\}^S_t=(s,\{\rho\}^x_t)$ gives an $\varOmega^x_x\ltimes S\tto S$~equivariant contraction of $S\times\mathscr R_x\to S$ onto the constant section $s\mapsto(s,\bar\rho)$, where $\bar\rho\in\mathscr R_x$ denotes the value of the constant map $\{-\}^x_0$. \end{proof}

\begin{prop}\label{prop:14A.7.5} Let an\/ \(\varOmega\)~equivariant\/ \(C^\infty\)\mdash fibration\/ \(\mathscr R\to X\) be given which is locally transversely trivializable and whose fiber\/ \(\mathscr R_x\) over each base-point\/ \(x\) is\/ \(\varOmega^x_x\)~equivariantly\/ \(C^\infty\)\mdash contractible. Let\/ \(U\subset\overline U\subset V\) be any invariant open subsets of\/ \(X\). Then, the equivariant extension problem that for any equivariant local\/ \(C^\infty\) section\/ $\rho$ of\/ \(\mathscr R\) defined over\/ \(V\) requires prolonging\/ \(\rho\mathbin|U\) to some element of\/ \(\Gamma^\infty(X;\mathscr R)^\varOmega\) always admits solutions. \end{prop}

\begin{proof} By the preceding lemmas there must be some global equivariant contraction $\{-\}_t$ of $\mathscr R\to X$. Then, essentially the same argument as in \ref{exmp:10/21**} applies with $\{-\}_t$ in place of the orthogonalization trick retraction. \end{proof}

\begin{npar}\label{npar:14A.7.6} By a \emph{discrete obstruction} for \(\mathscr R\to X\) we shall mean the following set of data: an \'etale mapping (local diffeomorphism) \(R\to X\); a $C^\infty$ (left) Lie-group\-oid action \(\varOmega\times_XR\to R\); a fi\-ber-pre\-serv\-ing \(\varOmega\)~equivariant \(C^\infty\) mapping \(\mathscr R\to R\). For any element \(r\in R_x\) we shall let \(\mathscr R_x[r]\subset\mathscr R_x\) denote the preimage of \(r\) under the mapping \(\mathscr R\to R\), and for any equivariant global \(C^\infty\) section \(r\in\Gamma^\infty(X;R)^\varOmega\) we shall let \(\mathscr R[r]\subset\mathscr R\) denote the invariant subbundle whose fiber at \(x\) is \(\mathscr R_x[r(x)]\).

Given an arbitrary equivariant trivialization of the form \eqref{eqn:14A.7.1}, the pullback fibration \(i^*R\to S\) will be a discrete obstruction for the \(i^*\varOmega\tto S\)~equivariant \(C^\infty\)\mdash fibration \(i^*\mathscr R\to S\), and
\begin{equation}
\label{eqn:14A.7.3}
	i^*(\mathscr R[r])=(i^*\mathscr R)[i^*r].
\end{equation}
Assuming without loss of generality that \(S\) is connected, a $C^\infty$ section of \(i^*\mathscr R\) will lie entirely within \((i^*\mathscr R)[i^*r]\) if, and only if, so does its value at any given point of \(S\). Since the constant sections of the trivial fibration $S\times\mathscr R_x\to S$ are $C^\infty$, it follows that an element $(s,\rho)$ of $S\times\mathscr R_x$ corresponds to one of $(i^*\mathscr R)[i^*r]$ if, and only if, $\rho\in\mathscr R_x[r(x)]=(\mathscr R[r])_x$. \em Thus, whenever\/ \(\mathscr R\) is locally transversely trivializable, so must for any\/ \(r\in\Gamma^\infty(X;R)^\varOmega\) be\/ \(\mathscr R[r]\).\em \end{npar}

\begin{npar}\label{npar:14A.7.7} Let us for an arbitrary \(\varOmega\)~equivariant \(C^\infty\)\mdash fibration \(\mathscr R\to X\) consider the equivariant extension problem formulated with Proposition~\ref{prop:14A.7.5}. Given any discrete obstruction \(R\to X\) for \(\mathscr R\to X\) the composition of \(\rho\) with the equivariant fiber bundle map \(\mathscr R\to R\) will be an equivariant local $C^\infty$ section $r:V\to R$ to the \'etale map $R\to X$. Being able to prolong \(r\mathbin|U\) to an element of \(\Gamma^\infty(X;R)^\varOmega\) lying in the image of \(\mathscr R\to R\) is evidently a necessary condition for the resolvability of the equivariant extension problem for $\rho$. We express the circumstance that $r$ can be so prolonged by saying that ``the $R$~obstruction to extending $\rho\mathbin|U$ vanishes''. \end{npar}

\section{The primary obstruction}\label{sec:3}

The present section is substantially a continuation of \S\ref{sec:1}, the notational conventions of which carry over. We intend to apply the abstract theoretical framework of the previous section to the $P\tto M$~equivariant \(C^\infty\)\mdash fibration \(\mathscr R^\varGamma\to M\). To begin with, we need a criterion for this particular equivariant \(C^\infty\)\mdash fibration to be locally transversely trivializable.

\begin{npar}\label{npar:14A.7.8} We start by considering the following special situation. Let \(\nu:\varSigma\to\GL(N)\) be a representation of a proper transitive Lie groupoid \(\varSigma\tto Z\) on a vector bundle \(p:N\to Z\). Suppose that \(\varGamma=\varSigma\ltimes^\nu M\tto M\) arises as the restriction of the action groupoid \(\varSigma\ltimes^\nu N\tto N\) over some invariant open neighborhood \(M\subset N\) of the zero section \(\{0_z\in N_z:z\in Z\}\). We want to show that through any given point \(x\in M\) there exists a slice \(i:S\hookto M\) for which the pullback fibration \(i^*\mathscr R^\varGamma\to S\) admits equivariant trivializations of the form \eqref{eqn:14A.7.1}.

Let \(K=\ker\nu\) denote the kernel of our representation. This is a closed Lie subgroupoid of \(\varSigma\tto Z\), and we have \(\nef\varGamma=K\times_ZM\). The proper foliation groupoid \(P\tto M\) may thus be identified with the action groupoid \(\bar\varSigma\ltimes^{\bar\nu}M\tto M\), where \(\bar\varSigma\) denotes the quotient Lie groupoid $\varSigma/K\tto Z$ and \(\bar\nu:\bar\varSigma\to\GL(N)\) the induced representation. The action of \(\varGamma=\varSigma\ltimes^\nu M\tto M\) along \(P=\bar\varSigma\ftimes spM\xto{\bar\nu}M\) is accordingly given by
\begin{equation*}
	(g,\bar\nu([h])m)\cdot([h],m)=([gh],m).
\end{equation*}

Given \(x\in M_z=M\cap N_z\), where \(z=p(x)\in Z\), let us fix an open neighborhood \(B\subset Z\) of \(z\) so small that there exist local smooth sections \(\kappa:B\to\varSigma^z\) to \(t^z:\varSigma^z\to Z\) through \(1_z=\kappa(z)\). Let us put \(U=M\cap p^{-1}(B)\). Any choice of one such section \(\kappa\) will determine an isomorphism of Lie groupoids [`$\kappa(u)$' hereafter being short for `$\kappa(p(u))$']
\begin{equation}
\label{eqn:14A.7.4}
	U\times\varPi^x\simto\varPi\mathbin|P^U,
\quad	\bigl(u\emphpunct;g,\bar\nu([h])x;[h],x\bigr)\mapsto\bigl(g,\bar\nu([h\kappa(u)^{-1}])u;[h\kappa(u)^{-1}],u\bigr),
\end{equation}
which in turn will give rise to a trivialization of \(C^\infty\)\mdash fibrations over \(U\)
\begin{equation}
\label{eqn:14A.7.5}
	U\times\mathscr R^\varGamma_x\mathrel{%
\xymatrix@1@C=1.33em@M=.em{%
 \ar[r]^*{\approx}_{C^\infty}
 &}}\mathscr R^\varGamma\mathbin|U
\end{equation}
obtained by sending each pair $(u,\rho)\in U\times\Rep(\varPi^x;TP^x)$ to the element of $\Rep(\varPi^u;TP^u)$ whose inverse image along the Lie groupoid isomorphism induced by \eqref{eqn:14A.7.4} \(\varPi^x=\{u\}\times\varPi^x\simto\varPi^u\) corresponds to \(\rho\).

There is a canonical action of the isotropy group \(P^x_x\) on the fiber \(M_z\)
\begin{equation*}
	([k],x)\cdot u=\nu(k)u.
\end{equation*}
For any \(P^x_x\)~invariant open subset \(S\) of \(M_z\), let us define \(i:S\hookto M\) to be the set-the\-o\-ret\-ic inclusion \(S\subset M_z\subset M\). We have a canonical Lie groupoid homomorphism
\begin{equation*}
	P^x_x\ltimes S\to i^*P,
\quad	([k],x;s)\mapsto(\nu(k)s;[k],s;s).
\end{equation*}
The trivialization of \(C^\infty\)\mdash fibrations that \eqref{eqn:14A.7.5} induces upon restriction
\begin{equation*}
	S\times\mathscr R^\varGamma_x\mathrel{%
\xymatrix@1@C=1.33em@M=.em{%
 \ar[r]^*{\approx}_{C^\infty}
 &}}\mathscr R^\varGamma\mathbin|S=i^*\mathscr R^\varGamma
\end{equation*}
will be equivariant relative to this Lie groupoid homomorphism. We conclude by observing that one can always choose \(S\) so small around \(x\) as to make the same homomorphism an isomorphism.

Similar considerations apply to the \(P\tto M\)~equivariant \(C^\infty\)\mdash fibration \(\mathscr R^{\varGamma,\phi}\to M\) for an arbitrary \(P\tto M\)~invariant metric \(\phi\) on the longitudinal bundle of \(\varPi\tto P\). \end{npar}

For \(\varGamma\tto M\) general, now, it follows from the preceding analysis that our equivariant \(C^\infty\)\mdash fibration \(\mathscr R^\varGamma\to M\) must be locally transversely trivializable whenever \(\varGamma\tto M\) is \emph{locally invariantly linearizable} in the sense that it is possible to cover \(M\) with invariant open subsets \(U\) over which the restrictions \(\varGamma\mathbin|U\tto U\) are isomorphic to Lie groupoids of the special type considered in \ref{npar:14A.7.8}. Ditto for $\mathscr R^{\varGamma,\phi}\to M$.

\begin{prop}\label{prop:14A.7.9} Let\/ \(\varGamma\tto M\) be any Lie groupoid which is regular, proper, and locally invariantly linearizable. Let\/ \(R\to M\) be a discrete obstruction for the\/ \(P\tto M\)~equivariant\/ \(C^\infty\)\mdash fibration\/ \(\mathscr R^\varGamma\to M\). Suppose that the preimage\/ \(\mathscr R^\varGamma_x[r]\subset\mathscr R^\varGamma_x\) of every\/ \(r\in R_x\) under the equivariant fi\-ber-bun\-dle projection\/ \(\mathscr R^\varGamma\to R\) is (if non-emp\-ty)\/ \(P^x_x\)~equivariantly\/ \(C^\infty\)\mdash contractible. Then, given\/ $U,V\subset M$ invariant, open, and with\/ $\overline U\subset V$, for any equivariant local\/ \(C^\infty\) section\/ $\rho\in\Gamma^\infty(V;\mathscr R^\varGamma)^P$ it will be possible to prolong\/ $\rho\mathbin|U$ to some element of\/ $\Gamma^\infty(M;\mathscr R^\varGamma)^P$ if, and only if, the\/ $R$~obstruction to extending\/ $\rho\mathbin|U$ vanishes. \end{prop}

\begin{proof} Proposition~\ref{prop:14A.7.5} applies to any equivariant $C^\infty$\mdash fibration of the form $\mathscr R^\varGamma[r]$, which must itself be locally transversely trivializable in virtue of the remarks near the end of \ref{npar:14A.7.6}. \end{proof}

By Weinstein's linearization theorem for proper regular groupoids \cite{Wein}, the hypotheses on $\varGamma\tto M$ will be satisfied when $\varGamma\tto M$ is source proper.

\begin{npar}\label{npar:14A.7.10} We next want to describe our main example of a discrete obstruction for $\mathscr R^\varGamma\to M$, which we call the \emph{primary obstruction} and denote $R^\varGamma\to M$. For every $x\in M$,%
\begin{subequations}
\label{eqn:14A.7.8}
\begin{equation}
\label{eqn:14A.7.8a}
	R^\varGamma_x\bydef\{(z,\sigma)\mathrel|z\in P^x\text{, }\sigma:\varPi_z^z\to\GL(T_zP^x)\text{ Lie group homomorphism}\}\big/\mathord\sim,
\end{equation}
where we write $(z_0,\sigma_0)\sim(z_1,\sigma_1)$ to express the circumstance that we can find arrows $h\in\varPi(z_0,z_1)$ and invertible linear maps $A\in\Lis(T_{z_0}P^x,T_{z_1}P^x)$ for which $\sigma_1(hkh^{-1})=A\sigma_0(k)A^{-1}~\forall k\in\varPi^{z_0}_{z_0}$. If we fix some point $z\in P^x$, every class in $R^\varGamma_x$ will be representable in the form $(z,\sigma)$. For any two homomorphisms $\sigma_0,\sigma_1:\varPi^z_z\to\GL(T_zP^x)$, the relation $(z,\sigma_0)\sim(z,\sigma_1)$ will hold if, and only if, $\sigma_0\sim_{\GL}\sigma_1$ i.e.~the two are intertwined by some element of $\GL(T_zP^x)$. Therefore, letting $1:M\to P$ denote the unit map of $P\tto M$, and recalling the natural identifications $\varPi^{1x}_{1x}=\nef\varGamma^x_x$ and $T_{1x}P^x=\varLambda_x$, we could also write
\begin{equation}
\label{eqn:14A.7.8b}
	R^\varGamma_x\bydef\{\sigma:\nef\varGamma_x^x\to\GL(\varLambda_x)\text{ Lie group homomorphism}\}\big/\mathord\sim_{\GL}.
\end{equation}
\end{subequations}

We proceed to introduce a collection of local trivializations for the fiber bundle \(R^\varGamma:=\bigcup_{x\in M}R^\varGamma_x\to M\) which are going to make \(R^\varGamma\) a $C^\infty$ covering space of \(M\). As $\nef\varGamma$ is a proper Lie bundle, it must be locally trivial. It must therefore be possible to cover $M$ with open subsets $U$ over which there exist local trivializations $\varphi:\nef\varGamma\mathbin|U\simto G\times U$ of Lie groupoids ($G$ some compact Lie group) and $\xi:\varLambda\mathbin|U\simto U\times\R^q$ of vector bundles. For any such $\varphi$ and $\xi$, we declare the following bijection to be a diffeomorphism
\begin{equation*}
	U\times\Rep(G;\R^q)\mathopen/\GL(q,\R)\longapproxto R^\varGamma\mathbin|U,
\quad	(u,[A])\mapsto[\bigl(1_u,{\xi_u}^{-1}(A\circ\varphi^u_u)\xi_u\bigr)].
\end{equation*}
Note that the quotient space $\Rep(G;\R^q)\mathopen/\GL(q,\R)$ is {\em discrete}, in the sense that every $C^\infty$ mapping of a smooth manifold into this space is locally constant (because any two representations of a compact Lie group $G$ which are close enough to one another must be equivalent), and may therefore be regarded as a discrete differentiable manifold. It is easy to check that the local trivializations obtained in this way are pairwise \(C^\infty\)\mdash compatible.

The action of \(P\tto M\) on \(R^\varGamma\to M\) makes arrows \([k]\in P\) operate by the rule
\begin{equation}
\label{eqn:14A.7.10}
	R^\varGamma_{sk}\ni[\bigl(k(z),\sigma\bigr)]\bydef[\mapsto][\bigl(z,{T_zk}^{-1}(\sigma\circ k^z_z)T_zk\bigr)]\in R^\varGamma_{tk}.
\end{equation}
Finally, the fi\-ber-bun\-dle projection \(\mathscr R^\varGamma\to R^\varGamma\) is given by
\begin{equation}
\label{eqn:14A.7.11}
	\mathscr R^\varGamma_x\ni\rho\bydef[\mapsto][(z,\rho^z_z)]\in R^\varGamma_x.
\end{equation} \end{npar}

\noindent{\it Comment.} We mention in passing that the ``monodromy map'' of \cite{JM} is somehow related to the problem of detecting possible non-van\-ish\-ing primary obstructions to the extension of locally given equivariant $C^\infty$ sections of $\mathscr R^\varGamma\to M$.

\begin{npar}\label{npar:14A.7.11} We end the current section with a few observations concerning the space \(\mathscr R^\varGamma_x[r]\) for an arbitrary primary obstruction class \(r\in R^\varGamma_x\). A more systematic analysis will be carried out in the next section.

a\spacefactor3000) Given any \(P\tto M\)~invariant metric \(\phi\) on the longitudinal bundle of \(\varPi\tto P\), the orthogonalization trick of \S\ref{sec:1} provides an equivariant strong deformation retraction of the space \(\mathscr R^\varGamma_x[r]\) onto its subspace \(\mathscr R^{\varGamma,\phi}_x[r]:=\mathscr R^\varGamma_x[r]\cap\mathscr R^{\varGamma,\phi}\). The two \(P^x_x\)~orbispaces \(\mathscr R^\varGamma_x[r]\) and \(\mathscr R^{\varGamma,\phi}_x[r]\) are therefore equivalent from the viewpoint of the study of their homotopical properties such as, for instance, equivariant contractibility.

b\spacefactor3000) One can always represent $r$ at any given point $z\in P^x$ by a {\em $\phi^x_z$~orthogonal}\/ tangent isotropy representation $\sigma:\varPi_z^z\to\GO(T_zP^x,\phi^x_z)$. [{\it Proof\emphpunct:} To obtain $\sigma$, apply the orthogonalization trick to any tangent isotropy representation $\alpha:\varPi^z_z\to\GL(T_zP^x)$ such that $(z,\alpha)\in r$; explicitly, write down the polar decomposition $HP$ of the identity map on $T_zP^x$ relative to $\phi^x_z$ and to some $\alpha$\mdash invariant metric, and then put $\sigma=H^{-1}\alpha H$.]

c\spacefactor3000) By definition, $\mathscr R^{\varGamma,\phi}_x=\Rep(\varPi^x;TP^x,\phi^x)$ and, for every tangent isotropy representation $\sigma:\varPi_z^z\to\GL(T_zP^x)$ such that $(z,\sigma)\in r$, $\mathscr R^{\varGamma,\phi}_x[r]=\{\rho\in\Rep(\varPi^x;TP^x,\phi^x)\mathrel|\rho^z_z\sim_{\GL}\sigma\}$. In fact, provided $\sigma$ is chosen $\phi^x_z$~orthogonal,
\begin{equation}
\label{eqn:14A.7.12}
	\mathscr R^{\varGamma,\phi}_x[r]=\{\rho\in\Rep(\varPi^x;TP^x,\phi^x)\mathrel|\rho^z_z=H\sigma H^{-1}~\exists H\in\GO(T_zP^x,\phi^x_z)\}.
\end{equation}
[{\it Proof\emphpunct:} Our claim is that whenever $\sigma_0,\sigma_1:\varPi^z_z\to\GO(T_zP^x,\phi^x_z)$ are intertwined by, say, $A\in\GL(T_zP^x)$, they are also intertwined by some element of $\GO(T_zP^x,\phi^x_z)$. Write down the polar decomposition $HP$ of $A$ relative to the metric $\phi^x_z$. At the expense of replacing $\sigma_1$ with an orthogonally equivalent representation, it will not be restrictive to assume that $A=P$ is symmetric positive definite. Then for $P_t:=\exp(t\log P)$ consider the homotopy $\alpha_t:=P_t\sigma_0{P_t}^{-1}:\varPi^z_z\to\GL(T_zP^x)$ linking $\alpha_0=\sigma_0$ to $\alpha_1=\sigma_1$. By the usual averaging argument, there will be a $1$\mdash parameter family of metrics $\hat\phi^x_{z,t}$ on $T_zP^x$ with $\hat\phi^x_{z,0}=\hat\phi^x_{z,1}=\phi^x_z$ such that each $\alpha_t$ is a $\hat\phi^x_{z,t}$\mdash orthogonal representation. Now for each $t$ write down the polar decomposition $\hat H_t\hat P_t$ of the identity map on $T_zP^x$ relative to $\phi^x_z$ and $\hat\phi^x_{z,t}$, and set
\begin{equation*}
	\sigma_t:={\hat H_t}^{-1}\alpha_t\hat H_t:\varPi^z_z\to\GO(T_zP^x,\phi^x_z).
\end{equation*}
This is a $1$\mdash parameter family of $\phi^x_z$\mdash orthogonal tangent isotropy representations linking $\sigma_0$ to $\sigma_1$. In view of the compactness of $\varPi^z_z$, any two representations belonging to such a family must be orthogonally equivalent.] \end{npar}

\section{The transitive case}\label{sec:4}

Throughout this section we shall be dealing with an arbitrary proper transitive Lie groupoid \(\varSigma\tto Z\) over a Riemannian manifold $Z=(Z,\gamma)$.

Because of transitivity, we have a canonical trivialization of \(C^\infty\)\mdash fibrations over \(Z\)
\begin{equation}
\label{eqn:14A.8.1}
	\mathscr R^\varSigma\approxto Z\times\Rep(\varSigma;TZ)
\end{equation}
arising as follows. First, the quotient \(\varSigma/\nef\varSigma\tto Z\) may be identified with the pair groupoid \(Z\times Z\tto Z\), so for any base point \(z\) the action \eqref{eqn:14A.6.10} yields a \(C^\infty\) trivialization \(Z\times\mathscr R^\varSigma_z\approxto\mathscr R^\varSigma\); second, the projection \(\varSigma\ltimes\varSigma^z/\varSigma^z_z\to\varSigma\emphpunct, (g,[h])\mapsto g\) gives an isomorphism of Lie groupoids which permits to identify \(\mathscr R^\varSigma_z\) with \(\Rep(\varSigma;TZ)\). Making no distinction between the given Riemannian metric \(\gamma\) on \(Z\) and the associated \(Z\times Z\tto Z\)~invariant metric on \(\ker T\pr_2\subset T(Z\times Z)\), we have that the two subbundles \(\mathscr R^{\varSigma,\gamma}\subset\mathscr R^\varSigma\) and \[Z\times\Rep(\varSigma;TZ,\gamma)\subset Z\times\Rep(\varSigma;TZ)\] correspond to each other under this trivialization.

The preceding considerations serve to legitimate the following notational abuses,%
\begin{subequations}
	\label{eqn:14A.8.2}
\begin{gather}
	\label{eqn:14A.8.2a}	\mathscr R^\varSigma=\Rep(\varSigma;TZ)
\\	\label{eqn:14A.8.2b}	\mathscr R^{\varSigma,\gamma}=\Rep(\varSigma;TZ,\gamma)
\end{gather}
which we take for granted in the sequel. For similar reasons we are entitled to write
\begin{equation}
	\label{eqn:14A.8.2c}	R^\varSigma=\{(z,\sigma)\mathrel|z\in Z\emphpunct, \sigma:\varSigma_z^z\to\GL(T_zZ)\text{ Lie group homomorphism}\}\big/\mathord\sim,
\end{equation}
where as in \ref{npar:14A.7.10} the relation $(z,\sigma)\sim(z',\sigma')$ signifies the existence of arrows $h\in\varSigma(z,z')$ and invertible linear maps $A\in\Lis(T_zZ,T_{z'}Z)$ such that $\sigma'(hkh^{-1})=A\sigma(k)A^{-1}\emphpunct{ }\forall k\in\varSigma^z_z$. On the account of \eqref{eqn:14A.7.12}, for every orthogonal tangent isotropy representation \(\sigma:\varSigma_z^z\to\GO(T_zZ)\) we write
\begin{equation}
	\label{eqn:14A.8.2d}	\mathscr R^{\varSigma,\gamma}[\sigma]=\{\alpha\in\Rep(\varSigma;TZ,\gamma)\mathrel|\alpha^z_z=H\sigma H^{-1}\emphpunct{ }\exists H\in\GO(T_zZ)\}.
\end{equation}
\end{subequations}
As the notation itself suggests, the space thus defined only depends on the class $[(z,\sigma)]\in R^\varSigma$.

\begin{npar}\label{npar:14A.8.1} We say that a tangent isotropy representation $\sigma:\varSigma_z^z\to\GL(T_zZ)$ is \emph{realizable} if an $\alpha\in\Rep(\varSigma;TZ)$ can be found so that $\alpha^z_z=\sigma$. Not every tangent isotropy representation is realizable; for example, the trivial (i.e.~constant) representation is realizable if and only if the manifold $Z$ is parallelizable. Any representation \(A\sigma A^{-1}\) intertwined to a realizable tangent isotropy representation \(\sigma\) by some orientation preserving invertible linear map \(A\in\GL_+(T_zZ)\) must itself be realizable. Furthermore, if an orthogonal tangent isotropy representation \(\sigma:\varSigma^z_z\to\GO(T_zZ)\) coincides with $\alpha^z_z$ for some tangent representation $\alpha:\varSigma\to\GL(TZ)$ then it is always possible to find \(\alpha:\varSigma\to\GO(TZ)\) orthogonal too. \end{npar}

\begin{npar}\label{npar:14A.8.2} Let \(\mathscr G^Z:=\Gamma^\infty\bigl(Z;\Hilb(TZ)\bigr)\) denote the \(C^\infty\)\mdash group of smooth sections of the Lie-group bundle
\begin{equation*}
	\Hilb(TZ):=\bigcup_{z\in Z}\GO(T_zZ)\to Z,
\end{equation*}
and \(\mathscr G^Z_+:=\Gamma^\infty\bigl(Z;\Hilb_+(TZ)\bigr)\) its normal subgroup consisting of those sections that take values in the open subbundle of orientation preserving tangent space isometries
\begin{equation*}
	\Hilb_+(TZ):=\bigcup_{z\in Z}\GO(T_zZ)\cap\GL_+(T_zZ)\to Z.
\end{equation*}

There is a canonical (left) \(C^\infty\) action of our \(C^\infty\)\mdash group \(\mathscr G^Z\) on the \(C^\infty\)\mdash space \(\mathscr R^{\varSigma,\gamma}\)~\eqref{eqn:14A.8.2b} of orthogonal tangent representations of \(\varSigma\tto Z\) which to \(f\in\mathscr G^Z\emphpunct, \rho\in\mathscr R^{\varSigma,\gamma}\) associates \(\rho^f\in\mathscr R^{\varSigma,\gamma}\) given by
\begin{equation}
\label{eqn:14A.8.4}
	\rho^f(h)\bydef f(th)\rho(h)f(sh)^{-1}.
\end{equation}
Note that for any orthogonal tangent isotropy representation \(\sigma:\varSigma^z_z\to\GO(T_zZ)\) the subspace \(\mathscr R^{\varSigma,\gamma}[\sigma]\subset\mathscr R^{\varSigma,\gamma}\)~\eqref{eqn:14A.8.2d} is \(\mathscr G^Z\)~invariant. \end{npar}

\begin{stmt}\label{stmt:14A.8.3} Suppose that\/ $Z$ is connected. Then for any given orthogonal tangent isotropy representation\/ \(\sigma:\varSigma^z_z\to\GO(T_zZ)\) the\/ $\mathscr G^Z$ invariant subspace\/ \(\mathscr R^{\varSigma,\gamma}[\sigma]\subset\mathscr R^{\varSigma,\gamma}\) consists of at most two\/ $\mathscr G^Z_+$~orbits, the exact number of which depends on how many\/ $\SO(T_zZ)$~conjugacy classes of realizable orthogonal tangent isotropy representations\/ \(\varSigma^z_z\to\GO(T_zZ)\) compose the\/ \(\GO(T_zZ)\)~conjugacy class of\/ \(\sigma\). \end{stmt}

\begin{proof} Given $\alpha,\rho\in\mathscr R^{\varSigma,\gamma}[\sigma]$ for which we can find $H\in\SO(T_zZ)$ such that $\rho^z_z=H\alpha^z_zH^{-1}$, let us pick an arbitrary $g\in\Gamma^\infty\bigl(Z;\Hilb_+(TZ)\bigr)$ such that $g(z)=H$ (evidently such global sections exist) and put $\varrho=\alpha^g$. Since $\varrho^z_z=\rho^z_z$, setting $f(th)=\rho(h)\varrho(h)^{-1}$ for every $h\in\varSigma^z$ yields a well-de\-fined $f\in\Gamma^\infty\bigl(Z;\Hilb_+(TZ)\bigr)$ such that $f(z)=\id$. Now $\rho=\varrho^f=(\alpha^g)^f=\alpha^{fg}$. \end{proof}

\begin{npar}\label{npar:14A.8.4} Any orthogonal tangent representation \(\rho\in\Rep(\varSigma;TZ,\gamma)\) determines a closed subbundle \((\rho)'_{\SO}\) of the proper Lie bundle \(\Hilb_+(TZ)\to Z\), the \emph{special (orthogonal) commutant} of \(\rho\), whose fiber at any base point \(z\) is the closed subgroup of $\SO(T_zZ)$
\begin{equation}
\label{eqn:14A.8.5}
	(\rho^z_z)'_{\SO}\bydef\{P\in\SO(T_zZ)\mathrel|P\rho(k)=\rho(k)P\emphpunct{ }\forall k\in\varSigma^z_z\}.
\end{equation}
For any local target section \(U\ni u\mapsto h_u\in\varSigma(z,-)\) the correspondence%
\begin{subequations}
\label{eqn:14A.8.6}
\begin{equation}
\label{eqn:14A.8.6a}
	\SO(T_zZ)\times U\longapproxto\Hilb_+(TZ)\mathbin|U\emphpunct, (P,u)\mapsto\rho(h_u)P\rho(h_u)^{-1}
\end{equation}
gives a local trivialization for \(\Hilb_+(TZ)\to Z\) over \(U\) which makes the two subbundles \((\rho^z_z)'_{\SO}\times U\subset\SO(T_zZ)\times U\) and \((\rho)'_{\SO}\mathbin|U\subset\Hilb_+(TZ)\mathbin|U\) correspond to each other. The bijection \((\rho^z_z)'_{\SO}\times U\approxto(\rho)'_{\SO}\mathbin|U\) induced by \eqref{eqn:14A.8.6a} turns out to be actually independent of the specific local smooth section \(u\mapsto h_u\) that we use in order to define it. We therefore have a well-de\-fined {\em global}\/ trivialization
\begin{equation}
\label{eqn:14A.8.6b}
	\tau^{\rho,z}:(\rho^z_z)'_{\SO}\times Z\approxto(\rho)'_{\SO}.
\end{equation}
\end{subequations}
The trivialization \(\tau^{\rho,z'}\) resulting from a different choice of base point \(z'\) will relate to \eqref{eqn:14A.8.6b} as follows,
\begin{equation*}
	\tau^{\rho,z'}\circ(c_{z,z'}\times\id)=\tau^{\rho,z}
\end{equation*}
where \(c_{z,z'}:(\rho^z_z)'_{\SO}\simto(\rho^{z'}_{z'})'_{\SO}\) denotes the Lie-group isomorphism that for any \(h\in\varSigma(z,z')\) operates as \(P\mapsto\rho(h)P\rho(h)^{-1}\). If now by a ``parallelization'' of a general Lie-group bundle \(K\to Z\) we understand, equivalently,
\begin{itemize}
 \item an equivalence class \([\tau]\) of Lie bundle isomorphisms \(\tau:G\times Z\simto K\) where \(\tau':G'\times Z\simto K\) and \(\tau\) are in the same class when there is a Lie-group isomorphism \(\vartheta:G\simto G'\) such that \(\tau'\circ(\vartheta\times\id)=\tau\),
 \item a flat~(= integrable) multiplicative connection with trivial holonomy on \(K\to Z\),
\end{itemize}
then we may summarize the previous remarks by saying that the special commutant \((\rho)'_{\SO}\) possesses an intrinsic parallelization say \(\pi^\rho:=[\tau^{\rho,z}]\) which depends only on \(\rho\).

For any parallelization $\pi$ of a general Lie bundle $K\to Z$, we let $\Gamma(K,\pi)\subset\Gamma^\infty(Z;K)$ denote the group of \(\pi\)\mdash parallel sections of \(K\to Z\). \end{npar}

\begin{stmt}\label{stmt:14A.8.5} The subgroup\/ \(\Gamma\bigl((\rho)'_{\SO},\pi^\rho\bigr)\subset\mathscr G^Z_+\) coincides with the stabilizer of\/ $\rho$ under the action\/ $\mathscr G^Z_+\circlearrowright\mathscr R^{\varSigma,\gamma}$~\eqref{eqn:14A.8.4},
\begin{equation}
\label{eqn:14A.8.8}
	\Stab_{\mathscr G^Z_+}(\rho)=\Gamma\bigl((\rho)'_{\SO},\pi^\rho\bigr)\makebox[.em][l]{.	\mathqed}
\end{equation} \end{stmt}

\begin{stmt}\label{stmt:14A.8.6} The map below is a\/ \(\mathscr G^Z_+\)~equivariant, open, embedding of\/ \(C^\infty\)\mdash spaces (cf.~Appendix~\ref{sec:A}), for any choice of ``reference point'' \(\alpha\in\mathscr R^{\varSigma,\gamma}\).
\begin{equation}
\label{eqn:14A.8.9}
\newdir{C}{{}*!/-5pt/@^{(}}%
\xymatrix@C=1.33em@M=.33em{%
	\mathscr G^Z_+\big/\left.\Gamma\bigl((\alpha)'_{\SO},\pi^\alpha\bigr)\right.
	\ar@{C->}[r]
	&	\mathscr R^{\varSigma,\gamma}\emphpunct, [f]\mapsto\alpha^f}
\end{equation}
For every realizable orthogonal tangent isotropy representation\/ \(\sigma:\varSigma^z_z\to\GO(T_zZ)\) such that\/ \(\mathscr R^{\varSigma,\gamma}[\sigma]\ni\alpha\), the image of this map is a subset of\/ \(\mathscr R^{\varSigma,\gamma}[\sigma]\). \end{stmt}

\begin{proof} In general, for an arbitrary $C^\infty$ action $\mathscr G\times\mathscr X\to\mathscr X$ of a \(C^\infty\)\mdash group $\mathscr G$ on a \(C^\infty\)\mdash space $\mathscr X$, the choice of a point $x\in\mathscr X$ determines an injective $C^\infty$ map $[g]\mapsto gx$ from the right-co\-set space $\mathscr G/\Stab_{\mathscr G}(x)$ into $\mathscr X$. In order for this map to be a \(C^\infty\)\mdash embedding, every $C^\infty$ map $S\to\mathscr Gx\subset\mathscr X$ ought to lift to a $C^\infty$ map $S\to\mathscr G$ locally. In the specific situation at hand, upon composition with the $C^\infty$ correspondence \[%
	\mathscr R^{\varSigma,\gamma}=\Rep(\varSigma;TZ,\gamma)\xto{(-)^z_z}\Hom\bigl(\varSigma^z_z,\GO(T_zZ)\bigr)
\] out of a $C^\infty$ map $\rho:S\to\mathscr G^Z_+\alpha\subset\mathscr R^{\varSigma,\gamma}$ we get a $C^\infty$ map $\rho^z_z:S\to\Hom\bigl(\varSigma^z_z,\GO(T_zZ)\bigr)$ with values in the $\SO(T_zZ)$~conjugacy class of $\alpha^z_z$. By the lemma~\ref{lem:A.1} of the appendix, it will be possible locally to find $C^\infty$ maps $H:S\to\SO(T_zZ)$ for which $\rho(s)^z_z=H(s)\alpha^z_zH(s)^{-1}$. We may then argue as in the proof of \ref{stmt:14A.8.3}, after observing that all the constructions therein involved are $C^\infty$. \end{proof}

\begin{npar}\label{npar:14A.8.7} Suppose that, besides \(\varSigma\tto Z\), we are given another proper transitive Lie groupoid \(\varSigma'\tto Z'\) over a Riemannian base \(Z'=(Z',\gamma')\) along with a homomorphism \(\kappa:\varSigma'\to\varSigma\) of Lie groupoids which on base level covers a Riemannian local diffeomorphism \(c:Z'\to Z\) (in other words an \'etale, smooth map $c$ such that \(c^*\gamma=\gamma'\)). The isometric, \'etale, smooth map \(c\) lifts to a homomorphism of Lie bundles
\begin{equation}
\label{eqn:14A.8.11}
	\Hilb_+(TZ')\longto\Hilb_+(TZ)\emphpunct, \GO(T_{z'}Z')\ni P'\mapsto(T_{z'}c)P'(T_{z'}c)^{-1}\in\GO(T_{c(z')}Z),
\end{equation}
from which one obtains a \(C^\infty\) homomorphism between the corresponding \(C^\infty\)\mdash groups of global sections
\begin{equation}
\label{eqn:14A.8.12}
	c^*:\mathscr G^Z_+\to\mathscr G^{Z'}_+,
\quad	(c^*f)(z')\bydef(T_{z'}c)^{-1}f(c(z'))(T_{z'}c).
\end{equation}
On the other hand, from the homomorphism \(\kappa:\varSigma'\to\varSigma\) of Lie groupoids one obtains a \(C^\infty\) map between the corresponding orthogonal tangent representation spaces
\begin{equation}
\label{eqn:14A.8.13}
	\kappa^*:\mathscr R^{\varSigma,\gamma}\to\mathscr R^{\varSigma',\gamma'},
\quad	(\kappa^*\rho)(h')\bydef(T_{th'}c)^{-1}\rho(\kappa(h'))(T_{sh'}c).
\end{equation}
The latter map is evidently equivariant with respect to the group homomorphism \eqref{eqn:14A.8.12}:
\begin{equation}
\label{eqn:14A.8.14}
	\kappa^*(\rho^f)=(\kappa^*\rho)^{c^*f}.
\end{equation}

For any choice of a ``reference point'' $\alpha\in\mathscr R^{\varSigma,\gamma}$, the homomorphism \eqref{eqn:14A.8.12} carries the subgroup $\Gamma\bigl((\alpha)'_{\SO},\pi^\alpha\bigr)$ of $\mathscr G^Z_+$ into the subgroup $\Gamma\bigl((\kappa^*\alpha)'_{\SO},\pi^{\kappa^*\alpha}\bigr)$ of $\mathscr G^{Z'}_+$, and thus descends to a well-de\-fined \(C^\infty\) mapping, which we shall denote \(c^*/\alpha\), between the right coset spaces $\mathscr G^Z_+\big/\Gamma\bigl((\alpha)'_{\SO},\pi^\alpha\bigr)$ and $\mathscr G^{Z'}_+\big/\Gamma\bigl((\kappa^*\alpha)'_{\SO},\pi^{\kappa^*\alpha}\bigr)$. From the equivariance \eqref{eqn:14A.8.14} of the mapping \eqref{eqn:14A.8.13}, we deduce the commutativity of the following diagram of \(C^\infty\) mappings:
\begin{equation}
\label{eqn:14A.8.15}
 \begin{split}
 \newdir{C}{{}*!/-5pt/@^{(}}%
\xymatrix@C=3em{%
	\mathscr G^Z_+\big/\left.\Gamma\bigl((\alpha)'_{\SO},\pi^\alpha\bigr)\right.
	\ar[d]^(.55){\smash{c^*/\alpha}}
	\ar@{C->}[r]^(.66){\alpha^{(-)}}
	&	\mathscr R^{\varSigma,\gamma}
		\ar[d]^(.55){\smash{\kappa^*}}
\\	\mathscr G^{Z'}_+\big/\left.\Gamma\bigl((\kappa^*\alpha)'_{\SO},\pi^{\kappa^*\alpha}\bigr)\right.
	\ar@{C->}[r]^(.69){(\kappa^*\alpha)^{(-)}}
	&	\mathscr R^{\varSigma',\gamma'}
}\end{split}
\end{equation} \end{npar}

\begin{npar}\label{npar:14A.8.8} Let \(\Aut(\varSigma)\) denote the group of Lie groupoid automorphisms of \(\varSigma\tto Z\), the group law being conventionally the opposite of composition: \(\kappa_1\kappa_2=\kappa_2\circ\kappa_1\). Let \(\Aut(\varSigma,\gamma)\subset\Aut(\varSigma)\) denote the subgroup of automorphisms which cover Riemannian self-i\-som\-e\-tries of the base \(Z=(Z,\gamma)\). \end{npar}

\begin{prop}\label{prop:14A.8.9} Let\/ \(\varSigma\tto Z\) be any Lie groupoid which is proper and transitive over a Riemannian manifold\/ \(Z=(Z,\gamma)\). Let\/ \(K\subset\Aut(\varSigma,\gamma)\) be a group of automorphisms\/ \(\kappa\) of\/ \(\varSigma\tto Z\) that cover Riemannian isometries\/ \(c\) on base level. Let\/ \(\sigma:\varSigma^z_z\to\GO(T_zZ)\) be any realizable orthogonal tangent isotropy representation of\/ \(\varSigma\tto Z\). Finally let\/ \(\alpha\) be a\/ \(K\)~invariant element of\/ \(\mathscr R^{\varSigma,\gamma}[\sigma]\) (viz.~one satisfying\/ \(\kappa^*\alpha=\alpha\) for all\/ \(\kappa\in K\)). Then the map\/ \eqref{eqn:14A.8.9}
\begin{equation*}
\newdir{C}{{}*!/-5pt/@^{(}}%
\xymatrix@C=1.33em@M=.33em{%
	\mathscr G^Z_+\big/\left.\Gamma\bigl((\alpha)'_{\SO},\pi^\alpha\bigr)\right.
	\ar@{C->}[r]
	&	\mathscr R^{\varSigma,\gamma}[\sigma]\emphpunct, [f]\mapsto\alpha^f}
\end{equation*}
is a\/ \(K\)~equivariant, open, \(C^\infty\)\mdash embedding onto a\/ \(K\)~invariant, closed (open) neighborhood of\/ \(\alpha\) within\/ \(\mathscr R^{\varSigma,\gamma}[\sigma]\), the action of\/ \(K\) on\/ \(\mathscr G^Z_+\big/\Gamma\bigl((\alpha)'_{\SO},\pi^\alpha\bigr)\) being given by\/ \(\kappa\cdot[f]:=[c^*f]\), that on\/ \(\mathscr R^{\varSigma,\gamma}[\sigma]\) by\/ \(\kappa\cdot\rho:=\kappa^*\rho\). Under the extra assumption that\/ \(Z\) be connected and that any other realizable orthogonal tangent isotropy representation in the\/ \(\GO(T_zZ)\)~conjugacy class of\/ \(\sigma\) be\/ \(\SO(T_zZ)\)~conjugated to\/ \(\sigma\), the same map will be a\/ \(C^\infty\)\mdash isomorphism between\/ \(\mathscr G^Z_+\big/\Gamma\bigl((\alpha)'_{\SO},\pi^\alpha\bigr)\) and all of\/ \(\mathscr R^{\varSigma,\gamma}[\sigma]\). \end{prop}

\begin{proof} The equivariance of the map follows immediately from the commutativity of the square diagram \eqref{eqn:14A.8.15}. All the remaining assertions have already been verified; cf.~\ref{stmt:14A.8.6}. \end{proof}

The proposition is going to be applied to the following situation: in the notations of \S\ref{sec:1}, let \(\varSigma\tto Z=(Z,\gamma)\) be \(\varPi^x\tto P^x=(P^x,\phi^x)\) for some \(P\tto M\)~invariant metric \(\phi\) on the longitudinal bundle of \(\varPi\tto P\), and take \(K\) to be the image of the group homomorphism \(P^x_x\to\Aut(\varPi^x,\phi^x)\) that makes every arrow \([k]\in P^x_x\) act as the Lie-group\-oid isomorphism \eqref{eqn:14A.6.11}. On the account of \ref{npar:14A.7.11}, the proposition will then provide an explicit ``homogeneous model'' for (a significant portion of) the \(P^x_x\)\mdash orbispace \(\mathscr R^{\varGamma,\phi}_x[r]\) for any primary obstruction class \(r\in R^\varGamma_x\) that can be realized as the image of some \(P^x_x\)~invariant element of \(\mathscr R^{\varGamma,\phi}_x\) under the projection \(\mathscr R^{\varGamma,\phi}_x\subset\mathscr R^\varGamma_x\to R^\varGamma_x\).

\section{The rank-two case}\label{sec:5}

In this section our analysis specializes further down to the case when our proper transitive Lie groupoid \(\varSigma\tto Z\) is source connected and its base $Z=(Z,\gamma)$ is a (connected) two-di\-men\-sion\-al Riemannian manifold. In the first place, we provide explicit criteria for the realizability of orthogonal tangent isotropy representations. After that, we study the associated representation spaces \eqref{eqn:14A.8.2d} showing, among other things, that the path-connected components of these spaces are always equivariantly contractible.

\subsection*{Criteria for realizability}

Given a Lie-group homomorphism $\sigma:\varSigma^z_z\to\GO(T_zZ)$, the problem of finding orthogonal tangent representations $\alpha:\varSigma\to\GO(TZ)$ for which $\sigma=\alpha^z_z$ is equivalent to the problem of extending $\sigma$ to some equivariant map of (right) principal bundles $\alpha^z:\varSigma^z\to\GO(TZ)^z$ from the $\varSigma^z_z$~bundle $t^z:\varSigma^z\to Z$ into the $\GO(T_zZ)$~bundle $t^z:\GO(TZ)^z\to Z$. Let us set $G=\varSigma^z_z$, $H=\GO(T_zZ)$, $P=\varSigma^z$ and $Q=\GO(TZ)^z$. The $\sigma:G\to H$~equivariant fi\-ber-bun\-dle maps $P\to Q$ can be identified with the global sections to the fibration
\begin{equation*}
	P\otimes Q:=(P\times_ZQ)/G\to Z,
\end{equation*}
where $G$ operates freely and properly on $P\times_ZQ$ from the right by the rule $(p,q)\cdot g:=(pg,q\sigma(g))$. Let us write $p\otimes q$ for the $G$\mdash class of a pair $(p,q)$; upon setting $gq:=q\sigma(g)^{-1}$, we have $pg\otimes q=p\otimes gq$, which justifies our notation. For any local section $(p,q):U\to P\times_ZQ$ to the fibration $P\times_ZQ\to Z$ we have a local fi\-ber-bun\-dle trivialization
\begin{equation*}
	\tau_{p,q}:U\times H\approxto(P\otimes Q)\mathbin|U\emphpunct, (u,h)\mapsto p(u)\otimes q(u)h.
\end{equation*}
For any other local section $(p',q'):U\to P\times_ZQ$, writing $p'(u)=p(u)g(u)$ and $q'(u)=q(u)h(u)$, the transition mapping ${\tau_{p,q}}^{-1}\circ\tau_{p',q'}$ is given by \[(u,h)\mapsto(u,h(u)h\sigma\bigl(g(u)\bigr)^{-1}).\]

\begin{prop}\label{prop:14A.9.1} Let\/ $\varSigma\tto Z$ be a transitive Lie groupoid which is compact and source connected over a (closed, connected) Riemannian\/ $2$\mdash manifold\/ $Z$ with finite fundamental group. Then, an arbitrary orthogonal tangent isotropy representation\/ $\sigma:\varSigma^z_z\to\GO(T_zZ)$ is realizable if, and only if, both diagrams below commute:
\begin{flalign}
\label{eqn:14A.9.1}
&&
 \begin{split}
\xymatrix@C=1.2em{%
	\pi_2(Z)
	\ar@{=}[d]
	\ar[r]^-{\partial_2}
	&	\pi_1(\varSigma^z_z)
		\ar[d]^{\pi_1(\sigma)}
\\	\pi_2(Z)
	\ar[r]^-{\partial_2}
	&	\pi_1\bigl(\GO(T_zZ)\bigr)
}\end{split}
 \begin{split}
\xymatrix@C=1.2em{%
	\pi_1(Z)
	\ar@{=}[d]
	\ar[r]^-{\partial_1}
	&	\pi_0(\varSigma^z_z)
		\ar[d]^{\pi_0(\sigma)}
\\	\pi_1(Z)
	\ar[r]^-{\partial_1}
	&	\pi_0\bigl(\GO(T_zZ)\bigr)
}\end{split}
&&
\end{flalign} \end{prop}

\begin{proof} Whenever there exists some orthogonal tangent representation $\alpha:\varSigma\to\GO(TZ)$ whose isotropy representation at $z$ coincides with $\sigma$, this must give rise to a mapping of pointed fibrations
\begin{equation}
\label{eqn:2016b.1}
 \begin{split}
 \newdir{C}{{}*!/-5pt/@^{(}}%
\xymatrix@C=1.2em{%
	(\varSigma^z_z,1_z)
	\ar[d]^{\sigma=\alpha^z_z}
	\ar[r]^\subset
	&	(\varSigma^z,1_z)
		\ar[d]^{\alpha^z}
		\ar@{->>}[r]^-{t^z}
		&	(Z,z)
			\ar@{=}[d]
\\	\bigl(\GO(T_zZ),\id\bigr)
	\ar[r]^-\subset
	&	\bigl(\GO(TZ)^z,\id\bigr)
		\ar@{->>}[r]^-{t^z}
		&	(Z,z)
}\end{split}
\end{equation}
along with a corresponding transformation of long exact sequences of homotopy groups which makes the necessity of the two commutativity conditions clear.

For the converse, we argue on the basis of the observations preceding the proposition. Since the continuous sections of a smooth fiber bundle can be approximated to any degree of precision by its smooth sections, we are allowed to reason topologically. By the classification of closed connected $2$\mdash manifolds (see \cite[Theorem~9.3.11]{Hir} or \cite[Theorems~8.3, 8.5 and 22.9]{Moise}), $Z$ must be homeomorphic to either the two-sphere, in which case $\pi_1(Z)$ is trivial and $Z$ is orientable, or the real projective plane, in which case $\pi_1(Z)=\Z/2\Z$ and $Z$ is not orientable.

Suppose that $Z$ is homeomorphic to the two-sphere. Let $D=\{\zeta\in\C:\abs\zeta\leqq 1\}$ denote the unit disk. As a CW-com\-plex, $Z$ is obtained by attaching a single 2-cell $c:D\to Z$ to $\{z\}$ via the map $\partial c:S^1\to Z$ which collapses the boundary of the disk $S^1=\partial D$ to $z$. The pullbacks $c^*P\to D\emphpunct, c^*Q\to D$ are topological fiber bundles over a contractible base, hence they are trivializable and thus admit global sections, say, $p:D\to c^*P\emphpunct, q:D\to c^*Q$ which may be assumed to take any desired value at $\zeta=1$. Define $g:S^1\to G\emphpunct, h:S^1\to H$ by setting $p(\zeta)=p(1)g(\zeta)\emphpunct, q(\zeta)=q(1)h(\zeta)$. We contend that the loop $S^1\ni\zeta\mapsto l(\zeta):=h(\zeta)^{-1}\sigma(g(\zeta))\in H$ must be homotopically trivial (contractible). Indeed, from the explicit construction of the long exact sequence of homotopy groups, we see that the two homotopy classes $[g]\in\pi_1(G)\emphpunct, [h]\in\pi_1(H)$ coincide with the images of the same element of $\pi_2(Z)$ under the two connecting homomorphisms $\partial_2:\pi_2(Z)\to\pi_1(G)\emphpunct, \partial_2:\pi_2(Z)\to\pi_1(H)$, whence by the first of the two diagrams \eqref{eqn:14A.9.1} $[l]=[\sigma\circ g]-[h]=\pi_1(\sigma)([g])-[h]=0\in\pi_1(H)$. Now, if we pick a contraction $l_t:S^1\to H$ from $l_1=l$ to $l_0=1_H$ which is locally constant near $t=0$, the section \[%
	D\ni\zeta\mapsto p(\zeta)\otimes q(\zeta)l_{\abs\zeta}(\zeta/\abs\zeta)\in c^*P\otimes c^*Q=c^*(P\otimes Q)
\] will coincide with $p\otimes q$ in a neighborhood of zero and be constant of value $p(1)\otimes q(1)$ along $S^1$, descending therefore to a well-de\-fined global section of $P\otimes Q$.

On the other hand, suppose $Z$ is homeomorphic to the real projective plane. In this case, $Z$ admits a double cover $c:S^2\to Z$ with $c(1,0,0)=z$ and the antipodal map as the only non-triv\-i\-al deck transformation. Let $e_\pm:D\hookto S^2$ denote the two embeddings \[%
	x+iy\mapsto(x,y,\pm\sqrt{1-x^2-y^2}).
\] Put $c_\pm=c\circ e_\pm$, and let $r:D\approxto D$ stand for the map $\zeta\mapsto-\zeta$. The pullbacks $c_\pm^*P\to D\emphpunct, c_\pm^*Q\to D$ are trivializable fiber bundles and hence admit global sections say $p_\pm:D\to c_\pm^*P\emphpunct, q_\pm:D\to c_\pm^*Q$. It is not restrictive to assume that $p_-$ coincides with $r^*p_+$ under the fi\-ber-bun\-dle identification $c_-^*P=(c_+\circ r)^*P=r^*c_+^*P$ and similarly that $q_-$ coincides with $r^*q_+$ under the fi\-ber-bun\-dle identification $c_-^*Q=(c_+\circ r)^*Q=r^*c_+^*Q$; we may express this by writing $p_-(\zeta)=p_+(-\zeta)\emphpunct, q_-(\zeta)=q_+(-\zeta)$ for $\zeta\in D$. Each one of $p_+,q_+$ can be made to take any desired value at $\zeta=1$. Define $g:S^1\to G\emphpunct, h:S^1\to H$ by setting $p_+(\zeta)=p_-(\zeta)g(\zeta)\emphpunct, q_+(\zeta)=q_-(\zeta)h(\zeta)$. We point out that $h:S^1\to H=\GO(T_zZ)$ must take values outside $\idc H=\SO(T_zZ)$, for otherwise $Z$ would have to be orientable. By the commutativity of the right-hand diagram \eqref{eqn:14A.9.1}, $\sigma\circ g:S^1\to H$ must then also range outside $\idc H$. In particular, $h(1)\sigma(g(1))^{-1}\in\idc H$. Thus, at the expense of deforming $q_+$ over a small neighborhood of $\zeta=-1$, we may without loss of generality assume that $h(1)=\sigma(g(1))$. We have fiber bundle trivializations
\begin{gather*}
	\tau_+:D\times H\approxto c_+^*P\otimes c_+^*Q=c_+^*(P\otimes Q)\emphpunct, (\zeta,h)\mapsto p_+(\zeta)\otimes q_+(\zeta)h
\\	\tau_-:D\times H\approxto c_-^*P\otimes c_-^*Q=c_-^*(P\otimes Q)\emphpunct, (\zeta,h)\mapsto p_-(\zeta)\otimes q_-(\zeta)h
\end{gather*}
with boundary transition map $\tau_-^{-1}\circ\tau_+:S^1\times H\approxto c_+^*(P\otimes Q)\mathbin|S^1=c_-^*(P\otimes Q)\mathbin|S^1\approxto S^1\times H\emphpunct, (\zeta,h)\mapsto(\zeta,h(\zeta)h\sigma\bigl(g(\zeta)\bigr)^{-1})$ and antipodal symmetry transition map
\begin{equation*}
	D\times H\xto{\tau_+}c_+^*(P\otimes Q)=r^*c_-^*(P\otimes Q)\xto{r^*(\tau_-)^{-1}}r^*(D\times H)\approxto D\times H\emphpunct, (\zeta,h)\mapsto(-\zeta,h).
\end{equation*}
With respect to $\tau_+$, the global sections of $P\otimes Q$ correspond to the functions $\eta:D\to H$ with the property that $\eta(-\zeta)=h(\zeta)\eta(\zeta)\sigma(g(\zeta))^{-1}$ for all $\zeta\in S^1$. We proceed to construct one such function $\eta$ which, moreover, satisfies $\eta(1)=1_H$. Since $\dim T_zZ=2$, every element of $H\smallsetminus\idc H\cong\GO(2)\smallsetminus\SO(2)$ squares to $1_H$, so we have $h(\zeta)=h(\zeta)^{-1}=h(-\zeta)$ and $\sigma(g(\zeta))=\sigma(g(\zeta))^{-1}=\sigma(g(\zeta)^{-1})=\sigma(g(-\zeta))$. The loop
\begin{equation*}
	S^1\ni\zeta\mapsto h(\zeta)\sigma(g(\zeta))^{-1}=h(\zeta)^{-1}\sigma(g(\zeta))=\bigl(h(1)^{-1}h(\zeta)\bigr)^{-1}\sigma\bigl(g(1)^{-1}g(\zeta)\bigr)\in\idc H
\end{equation*}
must therefore be 2-pe\-ri\-od\-ic. By the commutativity of the left-hand diagram \eqref{eqn:14A.9.1}, the same loop must be homotopically trivial. The unique loop $l:S^1\to\idc H$ such that $l(\zeta^2)=h(\zeta)\sigma(g(\zeta))^{-1}$ must then also be homotopically trivial; indeed $2[l]=0\in\pi_1(H)\cong\Z$ and thus $[l]=0$. Now, pick a homotopy of loops $l_t:S^1\to\idc H$ from $l_1=l$ to $l_0=1_H$ which is locally constant near $t=0$ and let $\eta:D\to\idc H$ be characterized by the two conditions $\eta(0)=1_H\emphpunct, \eta(\zeta)^2=l_{\abs\zeta}(\zeta^2/\abs\zeta^2)\sidetext(\zeta\in D)$. Then necessarily $\eta(-\zeta)=\eta(\zeta)\emphpunct{ }\forall\zeta\in D$, and $\eta(\zeta)=1_H$ for real $\zeta$. Furthermore, for $\zeta\in S^1$,
\begin{equation*}
	\eta(-\zeta)=\eta(\zeta)^2\eta(\zeta)^{-1}=l(\zeta^2)\eta(\zeta)^{-1}=h(\zeta)\sigma(g(\zeta))^{-1}\eta(\zeta)^{-1}=h(\zeta)\eta(\zeta)\sigma(g(\zeta))^{-1},
\end{equation*}
because conjugation by an element of $H\smallsetminus\idc H$ equals inversion on $\idc H$. \end{proof}

Under the same hypotheses as in the proposition (namely, $Z$ compact with finite fundamental group) the isotropy group $\varSigma^z_z$ will have at most two connected components, this being clear from the portion
\begin{equation*}
\xymatrix@C=1.2em{%
 \dotsb \ar[r]
 &	\pi_1(\varSigma^z)
	\ar[r]
	&	\pi_1(Z)
		\ar@{->>}[r]^-{\partial_1}
		&	\pi_0(\varSigma^z_z)
			\ar[r]
			&	\pi_0(\varSigma^z)=\ast}
\end{equation*}
of the long exact sequence of homotopy groups associated with the 1st pointed fibration \eqref{eqn:2016b.1}, since by the classification of closed connected $2$\mdash manifolds either $\pi_1(Z)=1$ (two-sphere) or $\pi_1(Z)=\Z/2\Z$ (real projective plane). In particular, $\varSigma^z_z$ will be connected when $Z$ is simply connected. The hypothesis that $\pi_1(Z)$ be finite is certainly satisfied when the source fibers themselves have finite fundamental groups. Indeed in this case $\pi_1(Z)$ is an extension of a finite group by a finite group and is thus itself finite.

The fundamental group $\pi_1(G)$ of any Lie group $G$ is abelian so it makes sense to take its torsion subgroup $\tor\pi_1(G)$ (consisting of all finite order elements) and its tor\-sion-free quotient $\free\pi_1(G):=\pi_1(G)/\tor\pi_1(G)$. If $\vartheta:G\to H$ is any Lie-group homomorphism then $\pi_1(\vartheta):\pi_1(G)\to\pi_1(H)$ carries $\tor\pi_1(G)$ into $\tor\pi_1(H)$ and hence induces a homomorphism of tor\-sion-free abelian groups $\free\pi_1(\vartheta):\free\pi_1(G)\to\free\pi_1(H)$.

\begin{prop}\label{prop:14A.9.2} Let\/ $\varSigma\tto Z$ be any source-con\-nect\-ed, compact, transitive, Lie groupoid over a Riemannian two-man\-i\-fold\/ $Z$ the fundamental group of which is finite. Let\/ $z\in Z$ be an arbitrary base point, and let\/ $\idc{\varSigma^z_z}$ denote the identity connected component of\/ $\varSigma^z_z$. Then, Lie group homomorphisms\/ $\vartheta:\idc{\varSigma^z_z}\to\SO(T_zZ)$ such that the left-hand diagram \eqref{eqn:14A.9.1} commutes exist if, and only if, the two conditions below are simultaneously satisfied:
\begin{enumerate}
\itemsep=0pt%
\def\labelenumi{\rm(\alph{enumi})}%
 \item The second homotopy group\/ $\pi_2(\varSigma^z)$ of the source fiber at\/ $z$ is trivial.
 \item The implication\/ $\Z e\cap\im(\free\partial_2)\neq 0\seq 2e\in\im(\free\partial_2)$ holds for every\/ $e\in\free\pi_1(\varSigma^z_z)$.
\end{enumerate}

When the fundamental group\/ $\pi_1(\varSigma^z)$ of the source fiber at\/ $z$ is finite, the problem of finding\/ $\vartheta$ so as to make the above-men\-tioned diagram commute admits at most one solution\/ $\vartheta$. By contrast, when\/ $\pi_1(\varSigma^z)$ is infinite, the same problem admits either denumerably many solutions, or no solution at all. \end{prop}

\begin{proof} As a preliminary step, we need to compute $\partial_2$ in the long exact sequence of homotopy groups associated with the pointed fibration at the bottom of \eqref{eqn:2016b.1},
\begingroup\small%
\begin{equation*}
 \newdir{ >}{{}*!/-5pt/@{>}}%
\xymatrix@C=1.2em{%
 \dotsb \ar[r]
 &	\pi_2\bigl(\GO(T_zZ)\bigr)
	\save[]+<.em,+2.7ex>*{=0}\restore
	\ar[r]
	&	\pi_2\bigl(\GO(TZ)^z\bigr)
		\ar@{ >->}[r]
		&	\pi_2(Z)
			\save[]+<.em,+2.7ex>*{\cong\Z}\restore
			\ar[r]^-{\partial_2}
			&	\pi_1\bigl(\GO(T_zZ)\bigr)
				\save[]+<.em,+2.7ex>*{\cong\Z}\restore
				\ar[r]
				&	\pi_1\bigl(\GO(TZ)^z\bigr)
					\ar[r]
					&	\pi_1(Z)
						\save[]+<.em,+2.7ex>*{\txt{finite}}\restore
						\ar[r]
						&	\dotsb}
\end{equation*}
\endgroup more precisely, we need to determine the integer $l=l(Z,z)\geqq 0$ for which $\partial_2=\pm l\cdot\id_\Z$. For this purpose, we can without loss of generality assume that $Z$ is orientable, for when it is not, it must admit an orientable, Riemannian, connected, double cover $\tilde Z$, and it is easy to see that $l(Z,z)=l(\tilde Z,\tilde z)$. Also, it not hard to see that $l(Z,z)$ is independent of the particular metric that we have on $Z$. We may henceforth assume that $Z$ is the standard two-sphere $S^2$ endowed with the standard metric and moreover that $z=(1,0,0)$. Now the god-given tangent representation of the action groupoid $\SO(3)\ltimes S^2\tto S^2$ provides a Lie-group\-oid isomorphism between this groupoid and the source-con\-nect\-ed orthogonal linear groupoid $\SO(TS^2)\tto S^2$, so we can compute $l(S^2,z)$ simply by looking at the long exact sequence of homotopy groups associated with the pointed fibration $(\SO(2),1)\hookto(\SO(3),1)\twoheadto(S^2,z)$, namely,
\begingroup\small%
\begin{equation*}
 \newdir{ >}{{}*!/-5pt/@{>}}%
\xymatrix@C=1.2em{%
 \dotsb \ar[r]
 &	\pi_2\bigl(\SO(2)\bigr)
	\save[]+<.em,+2.7ex>*{=0}\restore
	\ar[r]
	&	\pi_2\bigl(\SO(3)\bigr)
		\save[]+<.em,+2.7ex>*{=0}\restore
		\ar[r]
		&	\pi_2(S^2)
			\save[]+<.em,+2.7ex>*{=\Z}\restore
			\ar@{ >->}[r]
			&	\pi_1\bigl(\SO(2)\bigr)
				\save[]+<.em,+2.7ex>*{=\Z}\restore
				\ar@{->>}[r]
				&	\pi_1\bigl(\SO(3)\bigr)
					\save[]+<.em,+2.53ex>*{=\Z/2\Z}\restore
					\ar[r]
					&	\pi_1(S^2)
						\save[]+<.em,+2.7ex>*{=1}\restore
						\ar[r]
						&	\dotsb}
\end{equation*}
\endgroup from which we conclude that $l(S^2,z)=2$.

Clearly, the group homomorphisms $\varphi:\free\pi_1(\varSigma^z_z)\to\pi_1\bigl(\GO(T_zZ)\bigr)$ that satisfy the condition $\varphi\circ\free\partial_2=\partial_2:\pi_2(Z)\to\pi_1\bigl(\GO(T_zZ)\bigr)$ are in bijection with the solutions \textcircled{?}\ to the following problem,
\begingroup\small%
\begin{equation*}
 \newdir{ >}{{}*!/-5pt/@{>}}%
\xymatrix@C=1.4em{%
 \dotsb \ar[r]
 &	\pi_2(\varSigma^z_z)
	\save[]+<.em,+2.7ex>*{=0}\restore
	\ar[r]
	&	\pi_2(\varSigma^z)
		\ar@{ >->}[r]
		&	\pi_2(Z)
			\save[]+<.em,+2.7ex>*{\cong\Z}\restore
			\ar@{=}[d]
			\ar[r]^-{\partial_2}
			&	\pi_1(\varSigma^z_z)
				\ar@{.>}[d]^{\text{\textcircled{?}}}
				\ar[r]
				&	\pi_1(\varSigma^z)
					\ar[r]
					&	\pi_1(Z)
						\save[]+<.em,+2.7ex>*{\txt{finite}}\restore
						\ar[r]
						&	\dotsb
\\&	&	\dotsb \ar[r]
		&	\pi_2(Z)
			\save[]+<.em,-3.0ex>*{\cong\Z}\restore
			\ar@{ >->}[r]^-{\pm 2}
			&	\pi_1\bigl(\GO(T_zZ)\bigr)
				\save[]+<.em,-3.0ex>*{\cong\Z}\restore
				\ar[r]
				&	\dotsb}
\end{equation*}
\endgroup solutions which can possibly exist only when the group homomorphism $\partial_2:\pi_2(Z)\to\pi_1(\varSigma^z_z)$ is injective, equivalently, the homotopy group $\pi_2(\varSigma^z)$ vanishes. Suppose that this is the case. Then necessarily $\free\partial_2$ itself is an injective homomorphism of groups, because $\pi_2(Z)$ is tor\-sion-free. Its image $F=\im(\free\partial_2)\subset\free\pi_1(\varSigma^z_z)$ is a rank-one submodule of a free, finitely generated, $\Z$\mdash module, hence by a basic result in commutative algebra there exists some basis $e_1,\dotsc,e_r$ of $\free\pi_1(\varSigma^z_z)$ and some positive integer $d$ for which $F=\Z de_1$. Of course, it will not be restrictive to assume that $de_1=\free\partial_2(1)$. The group homomorphisms $\varphi:\free\pi_1(\varSigma^z_z)\to\Z$ satisfying the condition $\varphi\circ\free\partial_2=\pm 2\id_\Z$ correspond one-to-one to the $r$\mdash tuples of integers $\bigl(\varphi(e_1),\dotsc,\varphi(e_r)\bigr)\in\Z^r$ such that $d\varphi(e_1)=\varphi(de_1)=\varphi\bigl(\free\partial_2(1)\bigr)=\pm 2$. Clearly, such $r$\mdash tuples will not exist unless $2\Z e_1\subset F$, and they will not be unique unless $r=1$ equivalently $\pi_1(\varSigma^z)$ is finite. We finish the proof by invoking Lemma~\ref{lem:14A.9.8} from the appendix, which guarantees that an arbitrary homomorphism of groups $\varphi:\free\pi_1(\varSigma^z_z)\to\pi_1\bigl(\GO(T_zZ)\bigr)$ can be realized as $\varphi=\free\pi_1(\vartheta)$ for a unique Lie-group homomorphism $\vartheta:\idc{\varSigma^z_z}\to\SO(T_zZ)$. \end{proof}

We turn to the problem of extending $\vartheta:\idc{\varSigma^z_z}\to\SO(T_zZ)$ to the entire isotropy group $\varSigma^z_z$ so that the realizability condition \eqref{eqn:14A.9.1} is fulfilled. Observe that $\idc{\varSigma^z_z}=\varSigma^z_z$ when $\pi_1(Z)=1$ (two-sphere), so the problem in question only possibly arises when $\pi_1(Z)=\Z/2\Z$ (real projective plane), in which case the orthogonal linear groupoid $\GO(TZ)\tto Z$ is source connected, and the connecting homomorphism $\partial_1:\pi_1(Z)\to\pi_0\bigl(\GO(T_zZ)\bigr)$ is bijective. It follows that a realizable extension $\sigma$ of $\vartheta$, if it exists, has to carry $\varSigma^z_z\smallsetminus\idc{\varSigma^z_z}$ into $\GO(T_zZ)\smallsetminus\SO(T_zZ)$. Note that the realizability condition can never be fulfilled when $\pi_1(Z)=\Z/2\Z$ unless $\varSigma^z_z$ is disconnected.

\begin{lem}\label{lem:14A.9.3} Let\/ $G$ be an arbitrary Lie group for which\/ $\pi_0(G)=\Z/2\Z$. The problem of extending a given homomorphism\/ $\vartheta:\idc G\to\SO(2)$ to a homomorphism\/ $\sigma:G\to\GO(2)$ whose associated map on connected components\/ $\pi_0(\sigma):\pi_0(G)\to\pi_0\bigl(\GO(2)\bigr)$ is bijective admits solutions if and only if the two equations below are satisfied for one/ every\/ $h\in G\smallsetminus\idc G$.%
\begin{subequations}
\label{eqn:14A.9.3}
\begin{flalign}
\makebox[.em][l]{($\forall x\in\idc G$)}
	&\label{eqn:14A.9.3a}	& \vartheta(hxh^{-1})&=\vartheta(x)^{-1} &
\\	&\label{eqn:14A.9.3b}	&      \vartheta(h^2)&=1                 &
\end{flalign}
\end{subequations}
The action by conjugation of\/ $\SO(2)$ on the set of all such extensions\/ $\sigma$ of\/ $\vartheta$ is transitive. \end{lem}

\begin{proof} The necessity of the conditions \eqref{eqn:14A.9.3} is clear. As to their sufficiency, fix $h\in G\smallsetminus\idc G$ so that these conditions are satisfied. Pick $H\in\GO(2)$ for which $\det H=-1$ and, for all $x\in\idc G$, set
\begin{equation*}
	\sigma(x)=\vartheta(x),\quad \sigma(xh)=\vartheta(x)H.
\end{equation*}
It follows from \eqref{eqn:14A.9.3} that the map thus defined $\sigma:G\to\GO(2)$ is a group homomorphism. For any other extension $\sigma'$ of $\vartheta$ whose associated map $\pi_0(G)\to\pi_0\bigl(\GO(2)\bigr)$ is bijective, the difference $\sigma'(h)H^{-1}$ will belong to $\SO(2)$. Then, any $P\in\SO(2)$ such that $P^2=\sigma'(h)H^{-1}$ will be an intertwiner for $\sigma$ and $\sigma'$ because
\begin{equation*}
	P\sigma(xh)P^{-1}=P\vartheta(x)HP^{-1}=\vartheta(x)PHP^{-1}=\vartheta(x)P^2H=\sigma'(x)\sigma'(h)=\sigma'(xh).	\qedhere%
\end{equation*} \end{proof}

\begin{prop}\label{prop:14A.9.4} Let\/ $\varSigma\tto Z$ be as in Proposition~\ref{prop:14A.9.2}. Let\/ $z\in Z$ be a base point, and let\/ $\vartheta:\idc{\varSigma^z_z}\to\SO(T_zZ)$ be a Lie-group homomorphism that makes the diagram \eqref{eqn:14A.9.1} commute. Suppose that the source fiber\/ $\varSigma^z$ has a finite fundamental group. Then\/ $\vartheta$ will satisfy the condition \eqref{eqn:14A.9.3a} in the statement of Lemma~\ref{lem:14A.9.3} for any given\/ $h\in\varSigma^z_z\smallsetminus\idc{\varSigma^z_z}$ if, and only if, the following identity is true, where\/ $c_h\in\Aut(\varSigma^z_z)$ denotes conjugation by\/ $h$,
\begin{equation}
\label{eqn:14A.9.4}
	\free\pi_1(c_h)=-\id\in\Aut\bigl(\free\pi_1(\varSigma^z_z)\bigr)
\end{equation}
in which case\/ $\vartheta$ will also automatically satisfy the other condition \eqref{eqn:14A.9.3b}. \end{prop}

\begin{proof} By the lemma~\ref{lem:14A.9.8} of the appendix, $\vartheta\circ c_h=\vartheta^{-1}$ iff $\free\pi_1(\vartheta)\circ(\free\pi_1(c_h)+\id)=0$. As observed toward the end of the proof of Proposition~\ref{prop:14A.9.2}, the hypothesis that $\pi_1(\varSigma^z)$ is finite implies that $\free\pi_1(\varSigma^z_z)$ has rank one. The commutativity condition \eqref{eqn:14A.9.1} satisfied by $\vartheta$ then implies that $\free\pi_1(\vartheta)$ is a monomorphism of free abelian groups, because as we know from the proof of Proposition~\ref{prop:14A.9.2} the boundary homomorphism $\partial_2:\pi_2(Z)\to\pi_1\bigl(\GO(T_zZ)\bigr)$ is injective. The equivalence of the condition \eqref{eqn:14A.9.3a} to the condition $\free\pi_1(c_h)+\id=0$ is then clear.

As to the claim that when $\pi_1(\varSigma^z)$ is finite \eqref{eqn:14A.9.3b} is a consequence of \eqref{eqn:14A.9.3a}, we first of all point out that it is not restrictive to suppose $\pi_1(\varSigma^z)=1$ viz.~$\varSigma\tto Z$ is source simply connected. Indeed since $\varSigma\tto Z$ is source connected, it possesses a source simply connected universal cover \(\varepsilon:\tilde\varSigma\onto\varSigma\). The finiteness of $\pi_1(\varSigma^z)$ implies that $\tilde\varSigma\tto Z$ itself is a proper (hence compact) Lie groupoid. The surjectivity of the Lie-group homomorphism $\varepsilon^z_z:\tilde\varSigma^z_z\to\varSigma^z_z$ guarantees the existence of an element $\tilde h\in\tilde\varSigma^z_z$ such that $\varepsilon(\tilde h)=h$. Necessarily, $\tilde h\in\tilde\varSigma^z_z\smallsetminus\within*(\tilde\varSigma^z_z)^{\scriptscriptstyle\circ}$. The homomorphism $\tilde\vartheta=\vartheta\circ\varepsilon^z_z:\within*(\tilde\varSigma^z_z)^{\scriptscriptstyle\circ}\to\SO(T_zZ)$ satisfies the commutativity condition \eqref{eqn:14A.9.1}, whence, assuming the validity of the claim in the source simply connected case, $\vartheta(h^2)=\vartheta(\varepsilon(\tilde h)^2)=\vartheta(\varepsilon(\tilde h^2))=\tilde\vartheta(\tilde h^2)=1$.

From the hypothesis $\pi_1(\varSigma^z)=1$, by looking at the long exact sequence of homotopy groups associated with the pointed fibration $(\varSigma^z_z,1_z)\hookto(\varSigma^z,1_z)\twoheadto(Z,z)$ we conclude that $\free\partial_2:\pi_2(Z)\to\free\pi_1(\varSigma^z_z)$ is an isomorphism of groups. Then
\begin{equation*}
	\im\bigl(\free\pi_1(\vartheta)\bigr)=\im(\partial_2)=\{\mkern 1mu 2e:e\in\pi_1\bigl(\GO(T_zZ)\bigr)\mkern 1mu\}
\end{equation*}
(the latter equality resulting from the proof of Proposition~\ref{prop:14A.9.2}). The truth, in the source simply connected case, of the implication `$\eqref{eqn:14A.9.3a}\seq\eqref{eqn:14A.9.3b}$' is now a consequence of the following lemma\emphpunct: \em Let\/ $G$ be a compact Lie group for which\/ $\pi_0(G)=\Z/2\Z$. Let\/ $h\in G\smallsetminus\idc G$ be a group element such that\/ $\free\pi_1(c_h)=-\id\in\Aut\bigl(\free\pi_1(G)\bigr)$. Then every homomorphism\/ $\vartheta:\idc G\to\mathbb T$~(= one-to\-rus\/ $\R/\Z$) with the property
\begin{equation*}
	\im{}\bigl(\free\pi_1(\vartheta):\free\pi_1(G)\to\Z\bigr)=2\Z
\end{equation*}
satisfies the identity\/ \eqref{eqn:14A.9.3b}. \em Proof\emphpunct: In the notations of the proof of Lemma~\ref{lem:14A.9.8}, we have the following commutative diagram of Lie group homomorphisms,
\begin{equation*}
\xymatrix{%
	{\idc G}
	\ar[d]^(.49){c_h}
	\ar@{->>}[r]^-\tau
	\ar`u/6pt[r]!/u19pt/`[rrr]^-\vartheta[rrr]
	&	\mathbb T^r/\pr(A)
		\ar[d]^(.45)\inv
		\ar[r]^-\cong
		&	\mathbb T^r
			\ar[d]^(.45)\inv
			\ar[r]^-{\vartheta'}
			&	\mathbb T
				\ar[d]^(.45)\inv
\\	{\idc G}
	\ar@{->>}[r]^-\tau
	&	\mathbb T^r/\pr(A)
		\ar[r]^-\cong
		&	\mathbb T^r
			\ar[r]^-{\vartheta'}
			&	\mathbb T}
\end{equation*}
where $\inv$ stands for inversion (which on any abelian Lie group is a homomorphism) and $\cong$ denotes a fixed but otherwise arbitrary Lie-group isomorphism; the commutativity of the left-hand square follows by Lemma~\ref{lem:14A.9.8} from the hypothesis $\free\pi_1(c_h)=-\id$ because
\begin{equation*}
	\free\pi_1(\tau\circ c_h)=\free\pi_1(\tau)\circ-\id=-\id\circ\free\pi_1(\tau)=\free\pi_1(\inv\circ\tau).
\end{equation*}
Since $\free\pi_1(\tau):\free\pi_1(G)\to\pi_1\bigl(\mathbb T^r/\pr(A)\bigr)$ is a group isomorphism, we have $\im\bigl(\free\pi_1(\vartheta')\bigr)=\im\bigl(\free\pi_1(\vartheta)\bigr)=2\Z$ so $\vartheta'$ must be of the form
\begin{flalign*}
&&	\vartheta'(\zeta_1,\dotsc,\zeta_r)=\zeta_1^{2l_1}\dotsm\zeta_r^{2l_r}
&&	\makebox[.em][r]{$\exists(l_1,\dotsc,l_r)\in\Z^r$.}
\end{flalign*}
Now $h^2\in\idc G$ is a fix-point for $c_h$. The element $(\bar\zeta_1,\dotsc,\bar\zeta_r)\in\mathbb T^r$ corresponding to $\tau(h^2)$ under $\cong$ must therefore be a fix-point for inversion on $\mathbb T^r$; equivalently, $\bar\zeta_i^2=1$ for every $i=1,\dotsc,r$. Then
\begin{equation*}
	\vartheta(h^2)=\vartheta'(\bar\zeta_1,\dotsc,\bar\zeta_r)=\bar\zeta_1^{2l_1}\dotsm\bar\zeta_r^{2l_r}=1^{l_1}\dotsm 1^{l_r}=1.
\makebox[.em][l]{\quad Q.E.D.}	\qedhere%
\end{equation*} \end{proof}

\begin{cor}\label{cor:14A.9.5} Let\/ $\varSigma\tto Z$ be any transitive Lie groupoid which is compact and source simply connected over a two-di\-men\-sion\-al base\/ $Z$. Then, $\varSigma\tto Z$ possesses tangent representations if, and only if, the two conditions below are satisfied at one/ every base point\/ $z\in Z$:%
\begin{enumerate}
\def\labelenumi{\rm(\roman{enumi})}%
 \item The second homotopy group of the source fiber\/ $\varSigma^z$ is trivial.
 \item The identity\/ $\pi_1(c_h)=-\id\in\Aut\bigl(\pi_1(\varSigma^z_z)\bigr)$ holds for one/ every isotropy group element\/ $h\in\varSigma^z_z\smallsetminus\idc{\varSigma^z_z}$, where\/ $c_h\in\Aut(\varSigma^z_z)$ denotes conjugation by\/ $h$.
\end{enumerate} \end{cor}

\begin{proof} For $\varSigma^z$ simply connected, the condition (b) in the statement of Proposition~\ref{prop:14A.9.2} is certainly satisfied, for then $\partial_2$ in the long exact sequence of homotopy groups associated with the pointed fibration $(\varSigma^z_z,1_z)\hookto(\varSigma^z,1_z)\twoheadto(Z,z)$ and consequently $\free\partial_2$ in \eqref{eqn:14A.9.1} must be onto. By the same token, $\partial_1:\pi_1(Z)\simto\pi_0(\varSigma^z_z)$. In particular, $\varSigma^z_z$ must be disconnected when $\pi_1(Z)=\Z/2\Z$. We conclude by Propositions~\ref{prop:14A.9.1}~and \ref{prop:14A.9.4}. \end{proof}

\subsection*{Homotopical properties of the ``homogeneous model''}

Since $Z$ is two-di\-men\-sion\-al, $\Hilb_+(TZ)\to Z$ is a bundle of compact, abelian, Lie groups and therefore admits exactly one multiplicative connection\textemdash necessarily flat (cf.~\cite[Example~2.2]{2016a}). The holonomy of this connection gives a homomorphism of Lie groupoids
\begin{equation}
\label{eqn:14A.10.1}
	\hol_Z:\Pi_1(Z)\longto\Aut\bigl(\Hilb_+(TZ)\bigr)
\end{equation}
from the fundamental groupoid of $Z$ (cf.~\cite[Example~5.1]{MM}) into the Lie groupoid whose arrows are the Lie group isomorphisms $\SO(T_zZ)\simto\SO(T_{z'}Z)$. The connection has trivial holonomy when this homomorphism is constant on each isotropy group $\Pi_1(Z)^z_z=\pi_1(Z,z)$ i.e.~factors through the pair groupoid $Z\times Z\tto Z$. There is a second canonical homomorphism of Lie groupoids over $Z$%
\begin{subequations}
\begin{equation}
\label{eqn:14A.10.2}
	\GO(TZ)\longto\Aut\bigl(\Hilb_+(TZ)\bigr)
\end{equation}
which to every tangent space isometry $H:T_zZ\simto T_{z'}Z$ attaches the isomorphism of Lie groups $\SO(T_zZ)\simto\SO(T_{z'}Z)\emphpunct, P\mapsto HPH^{-1}$. Since $Z$ is $2$\mdash dimensional, \eqref{eqn:14A.10.2} descends to an isomorphism
\begin{equation}
\label{eqn:14A.10.3}
	\GO(TZ)\big/\Hilb_+(TZ)\longsimto\Aut\bigl(\Hilb_+(TZ)\bigr),
\end{equation}
\end{subequations}
whose existence enables one to conclude that $\Aut\bigl(\Hilb_+(TZ)\bigr)\tto Z$ is source connected if, and only if, the manifold $Z$ is non-o\-ri\-ent\-able.

The following three properties are in fact equivalent: \em%
\begin{enumerate}
\itemsep=0pt%
\def\labelenumi{\rm(\roman{enumi})}%
 \item $\Hilb_+(TZ)\to Z$ is a trivializable Lie bundle.
 \item The unique multiplicative connection on\/ $\Hilb_+(TZ)\to Z$ has trivial holonomy.
 \item The manifold\/ $Z$ is orientable.
\end{enumerate} \em

Let $\Gamma(Z)\subset\mathscr G^Z_+=\Gamma^\infty\bigl(Z;\Hilb_+(TZ)\bigr)$ denote the subgroup of all \emph{holonomic} (i.e.~flat) smooth sections of the Lie bundle $\Hilb_+(TZ)\to Z$. These are precisely those $f\in\mathscr G^Z_+$ that are invariant under the holonomy homomorphism \eqref{eqn:14A.10.1} i.e.~satisfy $\hol_Z(p)\bigl(f(sp)\bigr)=f(tp)$ for every $p\in\Pi_1(Z)$. It follows immediately from the preceding observations that $\Gamma(Z)\cong\SO(2)$ (as $C^\infty$\mdash groups) when $Z$ is orientable and that $\Gamma(Z)=\{\pm\id\}$ otherwise. We intend to investigate the homotopical properties of the homogeneous $C^\infty$\mdash space
\begin{equation}
\label{eqn:14A.10.4}
	\mathscr G^Z_+\big/\Gamma(Z)=\Gamma^\infty\bigl(Z;\Hilb_+(TZ)\bigr)\big/\{\hol_Z\text{-invariant sections}\}.
\end{equation}

\begin{lem}\label{lem:14A.10.1} Let\/ $Z$ be a two-di\-men\-sion\-al connected Riemannian manifold. Then, within the homogeneous\/ $C^\infty$\mdash space \eqref{eqn:14A.10.4}, the\/ $C^\infty$\mdash path-con\-nect\-ed component of the identity class\/ $[\id]$ is strongly\/ $C^\infty$ contractible onto\/ $\{[\id]\}$. \end{lem}

\begin{proof} For $f\in\mathscr G^Z_+$, consider the following (free) lifting problem,
\begin{equation}
\label{eqn:14A.10.5}
 \begin{split}
\xymatrix@C=3em{%
	&	\Skew(TZ)
		\ar[d]^\exp
\\	Z \ar[r]^(.45)f
	\ar@{.>}[ur]^-{\tilde f}
	&	\Hilb_+(TZ)
}\end{split}
\end{equation}
where $\Skew(TZ)$ denotes the vector bundle over $Z$ whose fiber at $z\in Z$ is $\Skew(T_zZ)$, the vector space of skew-sym\-met\-ric linear operators on $T_zZ$, and where $\exp$ denotes the fi\-ber-bun\-dle map \(%
	\Skew(T_zZ)\to\SO(T_zZ)\sidetext (z\in Z)\emphpunct, A\mapsto\sum_{n=0}^\infty\frac1{n!}A^n
\). Since $\dim T_zZ=2$, this map makes $\Skew(TZ)$ into a covering space of $\Hilb_+(TZ)$. By the homotopy lifting property of covering maps, for $f$ lying in $\idc{\mathscr G}^Z_+$~(= the $C^\infty$\mdash path-con\-nect\-ed component of the identity section $\id$ within $\mathscr G^Z_+$), our lifting problem \eqref{eqn:14A.10.5} will certainly admit solutions $\tilde f$. Now, given $f\in\idc{\mathscr G}^Z_+$, pick a random lift $\tilde f$ of $f$, and for all $t\in\R$, set%
\begin{equation}
\label{eqn:2016b.2}
	\{f\}_t:=[{\exp}\circ t\tilde f]\makebox[.pt][l]{	$\pmod{\Gamma(Z)}$.}
\end{equation}
The equivalence class mod~$\Gamma(Z)$ of ${\exp}\circ t\tilde f$ is easily seen to be independent of the choice of $\tilde f$, hence this expression yields a well-de\-fined $C^\infty$ one-pa\-ram\-e\-ter family of maps $\{-\}_t:\idc{\mathscr G}^Z_+\to\mathscr G^Z_+/\Gamma(Z)\sidetext(t\in\R)$. We let the reader conclude. \end{proof}

\noindent{\it Comments to the proof.} It is clear that the existence of a holonomic section of $\Skew(TZ)$ other than the zero section implies the existence of a vec\-tor-bun\-dle trivialization $Z\times\R\simto\Skew(TZ)$ and hence of a Lie-bun\-dle trivialization $Z\times\SO(2)\simto\Hilb_+(TZ)$ and, thus, orientability of $Z$. (Accordingly, in the non-o\-ri\-ent\-able case, a lift of $f$, whenever it exists, must be unique.) It follows that the two holonomic sections $\pm\id$ lie within the same $C^\infty$\mdash path-com\-po\-nent of $\mathscr G^Z_+$ if and only if $Z$ is orientable.

\vskip\topsep Let $c:Z\approxto Z$ be an arbitrary metric preserving self-dif\-feo\-mor\-phism of $Z$. The Lie-bun\-dle automorphism of $\Hilb_+(TZ)$ given by \eqref{eqn:14A.8.11} transforms the multiplicative connection on $\Hilb_+(TZ)$ necessarily into itself. As a consequence, the $C^\infty$\mdash group automorphism $c^*:\mathscr G^Z_+\simto\mathscr G^Z_+$ \eqref{eqn:14A.8.12} carries the holonomic subgroup $\Gamma(Z)\subset\mathscr G^Z_+$ into itself. Of course $c^*$ also carries the identity $C^\infty$\mdash path-com\-po\-nent $\idc{\mathscr G}^Z_+\subset\mathscr G^Z_+$ into itself.

\begin{lem}\label{lem:14A.10.2} Let\/ $Z$ be as in the preceding lemma. Let\/ $c:Z\approxto Z$ be any metric preserving self-dif\-feo\-mor\-phism of\/ $Z$. Then, the\/ $C^\infty$ contraction provided by\/ \eqref{eqn:2016b.2} is equivariant with respect to the\/ $C^\infty$\mdash automorphism of\/ $\mathscr G^Z_+/\Gamma(Z)$ induced by\/ $c^*:\mathscr G^Z_+\simto\mathscr G^Z_+$. \end{lem}

\begin{proof} For any $\phi\in\Gamma^\infty\bigl(Z;\Skew(TZ)\bigr)$ let $c^*\phi$ denote the section of $\Skew(TZ)$ given by $z\mapsto(T_zc)^{-1}\phi(c(z))(T_zc)$. Directly from the definitions, we see that $c^*({\exp}\circ\phi)={\exp}\circ c^*\phi$. In particular, in the notations of the preceding proof, if $\tilde f$ is a lift of $f$ as in \eqref{eqn:14A.10.5} then $c^*\tilde f$ is a lift of $c^*f$. Now, if we let the $C^\infty$\mdash automorphism of $\mathscr G^Z_+/\Gamma(Z)$ induced by $c^*$ be also denoted by $c^*$, then as contended, we have
\begin{align*}
	c^*(\{[f]\}_t)=c^*([{\exp}\circ t\tilde f])
	&             =[c^*({\exp}\circ t\tilde f)]
\\	&             =[{\exp}\circ tc^*\tilde f]
	              =\{[c^*f]\}_t
	              =\{c^*([f])\}_t.	\qedhere%
\end{align*} \end{proof}

By construction $\mathscr G^Z_+/\Gamma(Z)$ is a right-coset space of the $C^\infty$\mdash group $\mathscr G^Z_+$ and thus carries a $C^\infty$ action of $\mathscr G^Z_+$ by left translations \(\tau^f:[g]\mapsto[fg]\sidetext (f\in\mathscr G^Z_+)\). Then, setting $\{-\}^f_t:=\tau^f\circ\{-\}^\id_t\circ(\tau^f)^{-1}$ provides a strong $C^\infty$ contraction of the $C^\infty$\mdash path-com\-po\-nent of $[f]$ onto $\{[f]\}$. Notice that since $c^*\circ\tau^f=\tau^{c^*f}\circ c^*$,
\begin{equation}
\label{eqn:14A.10.7}
	c^*\circ\{-\}^f_t=\{-\}^{c^*f}_t\circ c^*.
\end{equation}

We conclude the present subsection with a few words about the classification of the path-com\-po\-nents of the homogeneous $C^\infty$\mdash space $\mathscr G^Z_+/\Gamma(Z)$. For $f,g\in\mathscr G^Z_+$, we shall write $f\sim g$ to signify that $f$ and $g$ lie within the same path-com\-po\-nent.

To begin with, the projection $\mathscr G^Z_+\to\mathscr G^Z_+/\Gamma(Z)$ gives a surjective map of pointed sets
\begin{equation*}
	\pi_0(\mathscr G^Z_+)\longto\pi_0\bigl(\mathscr G^Z_+/\Gamma(Z)\bigr),
\end{equation*}
equivariant with respect to the left-trans\-la\-tion action of the group $\pi_0(\mathscr G^Z_+)$ on these sets. The fibers of this map are precisely the right cosets for the stabilizer $\{[\id]_\sim,[-\id]_\sim\}$ of the base point of $\pi_0\bigl(\mathscr G^Z_+/\Gamma(Z)\bigr)$. Hence $\pi_0\bigl(\mathscr G^Z_+/\Gamma(Z)\bigr)=\pi_0(\mathscr G^Z_+)\big/\{[\id]_\sim,[-\id]_\sim\}$. We are thus reduced to the study of the group $\pi_0(\mathscr G^Z_+)$. Notice that when $Z$ is orientable, $\id\sim-\id$ so the map is 1:1, whereas when $Z$ is not orientable, $\id\not\sim-\id$ so the map is 2:1.

Let us fix some point $z\in Z$. Clearly, every element of $\mathscr G^Z_+$ can be deformed near $z$ in a $C^\infty$ fashion into an element of the normal subgroup
\begin{equation*}
	\mathscr G^Z_{+,z}:=\{\mkern 1mu f\in\mathscr G^Z_+:f(z)=\id_z\in\SO(T_zZ)\mkern 1mu\}\subset\mathscr G^Z_+,
\end{equation*}
hence the inclusion of $\mathscr G^Z_{+,z}$ into $\mathscr G^Z_+$ induces a surjective homomorphism of $\pi_0(\mathscr G^Z_{+,z})$ onto $\pi_0(\mathscr G^Z_+)$, whose kernel can be given the following explicit description. Let $B\approx\R^2$ be any open ball in $Z$ centered at $z\approx 0$. Take $\phi\in\Gamma^\infty\bigl(Z;\Skew(TZ)\bigr)$ with $\supp_Z\phi\subset B$ for which $\phi(z)$ is a generator of $\ker(\exp_z)\subset\Skew(T_zZ)$, i.e., $\ker(\exp_z)=\Z\phi(z)$. Then,
\begin{equation*}
	\ker{}\bigl(\pi_0(\mathscr G^Z_{+,z})\longto\pi_0(\mathscr G^Z_+)\bigr)=\{\mkern 1mu[{\exp}\circ n\phi]_{\sim_z}:n\in\Z\mkern 1mu\},
\end{equation*}
where for $f,g\in\mathscr G^Z_{+,z}$ we write $f\sim_zg$ to signify the existence of some $C^\infty$ path connecting $f$ to $g$ within $\mathscr G^Z_{+,z}$.

Next, let us consider the following covering fibration of pointed manifolds \[%
\xymatrix@C=1.33em{%
 (F,0_z) \ar[r]^-\subset
 &	\bigl(\Skew(TZ),0_z\bigr) \ar[r]^-\exp
	&	\bigl(\Hilb_+(TZ),\id_z\bigr)
}\] whose fiber $F=\ker(\exp_z)$ is a discrete subgroup of the additive group of the vector space $\Skew(T_zZ)$. The associated long exact sequence of homotopy groups gives a short exact sequence of {\em pointed sets}\/ (not {\em groups!}) \[%
\xymatrix@C=1.33em{%
 1 \ar[r]
 &	\pi_1\bigl(\Skew(TZ)\bigr) \ar[r]
	&	\pi_1\bigl(\Hilb_+(TZ)\bigr) \ar[r]^-{\partial_1}
		&	\pi_0(F)=F \ar[r]
			&	\ast.}\]

\begin{lem*} For\/ $f,g\in\mathscr G^Z_{+,z}$ one has\/ $f\sim_zg$ if and only if\/ $\partial_1\circ\pi_1(f)=\partial_1\circ\pi_1(g):\pi_1(Z)\to F$. \end{lem*}

\begin{proof} The `only if' direction is clear. For the converse, observe that the identities%
\begin{equation*}
	\partial_1\circ\pi_1(fg)=\partial_1\circ\pi_1(f)+\partial_1\circ\pi_1(g),\quad \partial_1\circ\pi_1(g^{-1})=-\partial_1\circ\pi_1(g)
\end{equation*}
hold for all $f,g\in\mathscr G^Z_{+,z}$. It follows that $\im\pi_1(fg^{-1})\subset\ker\partial_1=\im\pi_1(\exp)$. The covering map lemma then entails that $fg^{-1}$ admits a unique lift to $\Skew(TZ)$ with initial condition $0_z$ and is therefore homotopic to the identity section within $\mathscr G^Z_{+,z}$. Right translating the homotopy by $g$ yields a homotopy between $f$ and $g$ within $\mathscr G^Z_{+,z}$. \end{proof}

\noindent{\it Caveat.} The maps $\partial_1\circ\pi_1(f):\pi_1(Z)\to F$ need not be group homomorphisms. A~priori, they are just maps of pointed sets i.e.~they only need to satisfy $1\mapsto 0_z$.

\vskip\topsep Now, given $f\in\mathscr G^Z_{+,z}$ it follows from the above that $f\sim\id$ iff $[f]_{\sim_z}\in\ker{}\bigl(\pi_0(\mathscr G^Z_{+,z})\to\pi_0(\mathscr G^Z_+)\bigr)$ iff an integer $n\in\Z$ exists for which $f\sim_z{\exp}\circ n\phi=({\exp}\circ\phi)^n$ iff $\partial_1\circ\pi_1(f)=n\partial_1\circ\pi_1({\exp}\circ\phi)$ for some $n\in\Z$. One easily sees that for every $p\in\pi_1(Z)$
\begin{equation*}
	[\partial_1\circ\pi_1({\exp}\circ\phi)](p)=%
\begin{cases}
	0_z      &\text{if $\hol_Z(p)=\id$,}
\\	2\phi(z) &\text{if $\hol_Z(p)\neq\id$;}
\end{cases}
\end{equation*}
the condition $f\sim\id$ is thus equivalent to $\partial_1\circ\pi_1(f):\pi_1(Z)\to F$ being a composite map of the form
\begin{equation}
\label{eqn:14A.10.9}
\xymatrix@C=2.2em{%
 \pi_1(Z) \ar[r]^-{\hol_Z}
 &	\Aut\bigl(\SO(T_zZ)\bigr) \ar[r]^-{\id\mapsto 0_z}
	&	2F \ar[r]^-\subset
		&	F.}
\end{equation}

\begin{lem}\label{lem:14A.10.3} Two elements\/ $f,g\in\mathscr G^Z_{+,z}$ will be\/ $C^\infty$ homotopic within\/ $\mathscr G^Z_+$\textemdash in other words, they will lie within the same\/ $C^\infty$\mdash path-con\-nect\-ed component of\/ $\mathscr G^Z_+$\textemdash if, and only if, the map\/ $\partial_1\circ\pi_1(f)-\partial_1\circ\pi_1(g):\pi_1(Z)\to F$ admits a factorization of the form\/ \eqref{eqn:14A.10.9}. In the case when\/ $Z$ is orientable, in particular, $f\sim g\aeq\partial_1\circ\pi_1(f)=\partial_1\circ\pi_1(g)$. \qed \end{lem}

\subsection*{Conclusions}

For an arbitrary Riemannian manifold $Z$ let $\SO(TZ)\tto Z$ stand for the (maximal open) source-con\-nect\-ed subgroupoid of the orthogonal linear groupoid $\GO(TZ)\tto Z$. When $Z$ is connected and orientable, for either choice of orientation, $\SO(TZ)\tto Z$ will consist precisely of the orientation preserving tangent space isometries $T_zZ\simto T_{z'}Z\sidetext(z,z'\in Z)$. On the other hand, when $Z$ is connected and non-o\-ri\-ent\-able, $\SO(TZ)\tto Z$ will coincide with $\GO(TZ)\tto Z$.

Given any orthogonal tangent representation $\alpha:\varSigma\to\SO(TZ)$ of a source-con\-nect\-ed transitive Lie groupoid $\varSigma\tto Z$, the transformation of long exact sequences of homotopy groups associated with the mapping of pointed fibrations \eqref{eqn:2016b.1} shows that $\pi_0(\alpha^z_z)$ must be onto. If $Z$ is orientable, $\SO(TZ)^z_z=\SO(T_zZ)$ so $\alpha^z_z$ has to take values in the orientation preserving isometries. If $Z$ is not orientable then $\SO(TZ)^z_z=\GO(T_zZ)$ so there have to be group elements $h\in\varSigma^z_z$ for which $\det\alpha(h)=-1$.

For $Z$ two-di\-men\-sion\-al now, given any realizable orthogonal tangent isotropy representation $\sigma:\varSigma^z_z\to\GO(T_zZ)$, either $Z$ is orientable and $\sigma$ is \emph{central} i.e.~$(\sigma)'_{\SO}=\SO(T_zZ)$ or $Z$ is non-o\-ri\-ent\-able and $\sigma$ is \emph{irreducible} i.e.~$(\sigma)'_{\SO}=\{\pm\id\}$. Notice that an arbitrary Lie-group homomorphism $\sigma:\varSigma^z_z\to\GO(T_zZ)$ will be either central or irreducible in this sense\textemdash depending on whether it takes values in $\SO(T_zZ)$ or not\textemdash and that its being so is really a property of the $\GO(T_zZ)$~conjugacy class of $\sigma$ rather than of $\sigma$ itself, because $(H\sigma H^{-1})'_{\SO}=H(\sigma)'_{\SO}H^{-1}$ for every $H\in\GO(T_zZ)$. In the central case, there is only one possibility for the special orthogonal commutant $(\alpha)'_{\SO}\subset\Hilb_+(TZ)$ of an arbitrary orthogonal tangent representation $\alpha:\varSigma\to\SO(TZ)$ such that $\alpha^z_z=\sigma$, namely, the whole $\Hilb_+(TZ)$. This is an immediate consequence of the explicit description of $(\alpha)'_{\SO}$ \eqref{eqn:14A.8.6a}. On the other hand, in the irreducible case, $(\alpha)'_{\SO}\subset\Hilb_+(TZ)$ consists of the images of the two scalar multiples of the constant identity section, $\pm\id$, and is therefore a trivial $\Z/2\Z$\mdash bundle over $Z$. Summarizing\emphpunct: \em In the notations of \S\ref{sec:4}, for every orthogonal tangent representation\/ $\alpha:\varSigma\to\SO(TZ)$ of a source-con\-nect\-ed, proper, transitive Lie groupoid\/ $\varSigma\tto Z$ with base a Riemannian $2$\mdash manifold\/ $Z=(Z,\gamma)$,
\begin{equation}
\label{eqn:14A.11.1}
	\Gamma\bigl((\alpha)'_{\SO},\pi^\alpha\bigr)=\Gamma(Z).
\end{equation}

\em Let $\mathscr C^\alpha$ denote the $C^\infty$\mdash path-com\-po\-nent of $\alpha$ within $\mathscr R^{\varSigma,\gamma}$. Clearly, $\mathscr C^\alpha$ has to be contained within the (closed) image of the (open) embedding \eqref{eqn:14A.8.9}, and, therefore, $C^\infty$\mdash isomorphic\textemdash via that embedding\textemdash to the identity $C^\infty$\mdash path-com\-po\-nent of $\mathscr G^Z_+/\Gamma(Z)$. On the account of Lemma~\ref{lem:14A.10.1}, by transport of structure along the same embedding, we obtain a canonical, strong, $C^\infty$ contraction $\{-\}^\alpha_t:\mathscr C^\alpha\to\mathscr C^\alpha$ of $\mathscr C^\alpha$ onto $\{\alpha\}$. The commutativity of \eqref{eqn:14A.8.15}, in combination with Lemma~\ref{lem:14A.10.2}, grants the validity of the equation below for any automorphism $\kappa\in\Aut(\varSigma,\gamma)$ of $\varSigma\tto Z$ that covers a Riemannian isometry on base level.
\begin{equation}
\label{eqn:14A.11.3}
	\kappa^*\circ\{-\}^\alpha_t=\{-\}^{\kappa^*\alpha}_t\circ \kappa^*
\end{equation}
From this it follows in particular that the contraction $\{-\}^\alpha_t$ must be $\Stab_{\Aut(\varSigma,\gamma)}(\alpha)=\{\kappa\in\Aut(\varSigma,\gamma):\kappa^*\alpha=\alpha\}$~equivariant. In conclusion:

\begin{stmt}\label{stmt:14A.11.1} For an arbitrary proper transitive Lie groupoid\/ \(\varSigma\tto Z\) which is source connected over a two-di\-men\-sion\-al Riemannian manifold\/ $Z=(Z,\gamma)$, the $C^\infty$\mdash path-con\-nect\-ed component\/ \(\mathscr C^\alpha\) of each element\/ \(\alpha\) of the\/ \(C^\infty\)\mdash space\/ \(\mathscr R^{\varSigma,\gamma}\) is\/ \(C^\infty\) contractible onto\/ $\{\alpha\}$ through a canonical, strong, $C^\infty$ contraction\/ $\{-\}^\alpha_t:\mathscr C^\alpha\to\mathscr C^\alpha$ that is equivariant with respect to the action of the stabilizer subgroup\/ $\Stab_{\Aut(\varSigma,\gamma)}(\alpha)\subset\Aut(\varSigma,\gamma)$. \qed \end{stmt}

\begin{stmt}\label{stmt:14A.11.2} Let\/ $\varSigma\tto Z=(Z,\gamma)$ be as in the preceding statement. Suppose, in addition, that the base\/ $Z$ is compact with finite fundamental group. Then, for any realizable orthogonal tangent isotropy representation\/ $\sigma:\varSigma^z_z\to\GO(T_zZ)$, the\/ $C^\infty$\mdash space\/ $\mathscr R^{\varSigma,\gamma}[\sigma]$ deformation retracts onto any given one of its points, say, \(\alpha\), through a canonical, \(\Stab_{\Aut(\varSigma,\gamma)}(\alpha)\)~equivariant, strong, \(C^\infty\) contraction\/ \(\{-\}^\alpha_t:\mathscr R^{\varSigma,\gamma}[\sigma]\to\mathscr R^{\varSigma,\gamma}[\sigma]\). \end{stmt}

\begin{proof} In virtue of \ref{stmt:14A.11.1}, it will suffice to show that $\mathscr R^{\varSigma,\gamma}[\sigma]$ is $C^\infty$\mdash path-con\-nected. Now, $\sigma$ satisfies the hypothesis of Proposition~\ref{prop:14A.8.9} that any other realizable representation in its $\GO(T_zZ)$~conjugacy class be intertwined to it by some element of $\SO(T_zZ)$; indeed, given $H\in\GO(T_zZ)$ with $\det H=-1$, it cannot happen that both $\pi_1(\sigma)$ and $\pi_1(c_H\circ\sigma)=\pi_1(c_H)\circ\pi_1(\sigma)=-\pi_1(\sigma)$ fit in the left-hand diagram \eqref{eqn:14A.9.1}, for as we know from the proof of Proposition~\ref{prop:14A.9.2} the bottom $\partial_2$ in that diagram is injective. Then, by Proposition~\ref{prop:14A.8.9} and the identity~\eqref{eqn:14A.11.1}, we have a $C^\infty$\mdash isomorphism \(\mathscr G^Z_+/\Gamma(Z)\approxto\mathscr R^{\varSigma,\gamma}[\sigma]\). This reduces our task to proving that when $Z$ is compact with $\pi_1(Z)$ finite the $C^\infty$\mdash space $\mathscr G^Z_+/\Gamma(Z)$ is connected.

Suppose that $\pi_1(Z)=1$. Then there is only one map $\pi_1(Z)\to\Z$ sending $1\mapsto 0$: by Lemma~\ref{lem:14A.10.3}, $\mathscr G^Z_{+,z}$ and, consequently, $\mathscr G^Z_+$ have to be connected.

Suppose, on the other hand, that $\pi_1(Z)=\Z/2\Z$. Then $Z$ is non-o\-ri\-ent\-able and thus $(\hol_Z)^z_z:\pi_1(Z)\to\Aut\bigl(\SO(T_zZ)\bigr)$, being non-triv\-i\-al, must be an isomorphism. It follows that a map $\pi_1(Z)\to\Z$ sending $1\mapsto 0$ is ``null'' in the sense of \eqref{eqn:14A.10.9} if, and only if, it takes values in $2\Z$. By Lemma~\ref{lem:14A.10.3}, the elements of $\pi_0(\mathscr G^Z_+)$ are then classified by the maps $\pi_1(Z)\to\Z/2\Z$ with $1\mapsto 0$. Clearly, there are only two such maps, corresponding respectively to the connected component of $\id$ and that of $-\id$. \end{proof}

\section{Examples}\label{sec:6}

We are now done with the more technical part of our paper, and ready to discuss some applications of our theory to the study of Cartan connections.

\begin{prop}\label{lem:14A.11.6} Let\/ \(\varGamma\tto M\) be an arbitrary regular Lie groupoid of rank two which is source proper, source connected, and whose orbits\/ \(O_x\emphpunct, x\in M\) are (closed, connected, $2$\babelhyphen{nobreak})manifolds with finite fundamental groups. Then, given\/ \(U,V\subset M\) and\/ \(\rho\in\Gamma^\infty(V;\mathscr R^\varGamma)^P\) as in the statement of Proposition~\ref{prop:14A.7.9}, it is possible to extend\/ \(\rho\mathbin|U\) to some global equivariant\/ \(C^\infty\) section of\/ \(\mathscr R^\varGamma\) if, and only if, some\/ \(r\in\Gamma^\infty(M;R^\varGamma)^P\) can be exhibited for which\/ \(\rho\mathbin|U\in\Gamma^\infty(U;\mathscr R^\varGamma[r])^P\) and such that for every base point\/ \(x\) the following three conditions, the first of which involves the boundary maps\/ \(\partial_2\) in the long exact sequences of homotopy groups associated with the two pointed fibrations\/ \((\varGamma^x_x,1_x)\xto\subset(\varGamma^x,1_x)\xto{t^x}(O_x,x)\) and\/ \(\bigl(\GL(T_xO_x),\id\bigr)\xto\subset\bigl(\GL(TO_x)^x,\id\bigr)\xto{t^x}(O_x,x)\), are satisfied for at least one representative, say, \(\sigma:\nef\varGamma^x_x\to\GL(T_xO_x)\) of the class\/ \(r(x)\in R^\varGamma_x\):
\begin{enumerate}
\def\labelenumi{\rm(\arabic{enumi})}%
 \item The composite map\/ \(\pi_2(O_x)\xto{\partial_2}\pi_1(\varGamma^x_x)\xto{\pi_1(\sigma)}\pi_1\bigl(\GL(T_xO_x)\bigr)\) equals either\/ \(\partial_2\) or\/ \(-\partial_2\).
 \item \(\forall h\in\varGamma^x_x\smallsetminus\idc{(\varGamma^x_x)}\emphpunct, \pi_1(\sigma)\circ\pi_1(c_h)=-\pi_1(\sigma):\pi_1(\varGamma^x_x)\to\pi_1\bigl(\GL(T_xO_x)\bigr)\), where\/ \(c_h\in\Aut(\varGamma^x_x)\) denotes conjugation by\/ \(h\).
 \item \(\forall h\in\varGamma^x_x\smallsetminus\idc{(\varGamma^x_x)}\emphpunct, \sigma(h^2)=1\).
\end{enumerate} \end{prop}

\noindent{\it Remarks.} Each one of the three properties (1)\textendash (3) is stable under conjugation by elements of $\GL(T_xO_x)$ and hence is actually a property of the class $r(x)$ rather than of the specific representative $\sigma\in r(x)$. Clearly, if any of the three properties holds at some base point $x$ then by the equivariance of $r$ it also holds at every other base point in the orbit of $x$.

\begin{proof} ({\it Necessity.})\hskip.5em By virtue of \ref{stmt:14A.6.4}, any global prolongation of $\rho\mathbin|U$ will be of the form $\rho^\alpha$ for a unique $\alpha\in\Rep(\varGamma;\varLambda)$. Composition with the equivariant map of $C^\infty$\mdash fibrations \eqref{eqn:14A.7.11} will produce a global equivariant $C^\infty$ section $r^\alpha$ of $R^\varGamma$ with the property that $\rho\mathbin|U\in\Gamma^\infty(U;\mathscr R^\varGamma[r^\alpha])^P$. By construction, a representative of the class $r^\alpha(x)$ will be provided by the restriction of the isotropy representation $\alpha^x_x:\varGamma^x_x\to\GL(\varLambda_x)=\GL(T_xO_x)$ to $\nef\varGamma^x_x$. Since $\alpha^x_x$ is a realizable tangent isotropy representation of a transitive, rank-two, source-con\-nect\-ed, compact Lie groupoid $\varGamma\mathbin|O_x\tto O_x$ with finite $\pi_1(O_x)$, it will follow at once from Proposition~\ref{prop:14A.9.1} and Lemma~\ref{lem:14A.9.3}\textemdash say, after choosing a random $\alpha$~invariant vector-bun\-dle metric on $\varLambda$\textemdash that the conditions (1)~to (3) must be satisfied for the restriction of $\alpha^x_x$ to the identity component $\idc{\varGamma^x_x}\subset\nef\varGamma^x_x$ of $\varGamma^x_x$.

({\it Sufficiency.})\hskip.5em Given \(r\) as in the statement of the proposition, one can obtain a global extension of $\rho\mathbin|U$ by arguing as in the proof of Proposition~\ref{prop:14A.7.9}. We only need to show that under the present circumstances each fiber \((\mathscr R^\varGamma[r])_x=\mathscr R^\varGamma_x[r(x)]\) is \(P^x_x\)~equivariantly \(C^\infty\) contractible. Let us endow \(\varLambda\) with an arbitrary vec\-tor-bun\-dle metric \(\phi\). By the orthogonalization trick of \S\ref{sec:1}, \(\mathscr R^\varGamma_x[r(x)]\) deformation retracts \(C^\infty\), strongly, and \(P^x_x\)~equivariantly onto \(\mathscr R^{\varGamma,\phi}_x[r(x)]\subset\mathscr R^\varGamma_x[r(x)]\). The statement~\ref{stmt:14A.11.2} applied to \(\varPi^x\tto P^x\) guarantees that \(\mathscr R^{\varGamma,\phi}_x[r(x)]=\mathscr R^{\varPi^x,\phi^x}[r(x)]\) will be \(P^x_x\)~equivariantly \(C^\infty\) contractible on the condition that it contains some \(P^x_x\)~invariant element. The existence of such an element is tantamount to the existence of some orthogonal tangent representation say $\alpha$ of $\varGamma\mathbin|O_x\tto O_x$ satisfying $\alpha^x_x\mathbin|\nef\varGamma^x_x\in r(x)$. From the assumptions (1)\textendash (3), thanks to Proposition~\ref{prop:14A.9.1} and Lemma~\ref{lem:14A.9.3}, we deduce that some of the representatives $\sigma:\nef\varGamma^x_x\to\GO(T_xO_x)$ of the class $r(x)$ must be the restrictions to \(\nef\varGamma^x_x\) of realizable orthogonal tangent isotropy representations of $\varGamma^x_x$. The existence of \(\alpha\) is now clear. \end{proof}

One may summarize the conclusions of the proposition informally by saying that, under the indicated assumptions, the primary obstruction is the only obstruction to solving the extension problem for longitudinal representations.

For an arbitrary Lie groupoid \(\varGamma\tto M\), and for any integer \(q\geqq 0\), let us write \(M_q\) for \(\{x\in M:\dim O_x=q\}\) the set of base points that lie on \(q\)\mdash dimensional orbits.

\begin{thm}\label{prop:14A.11.7} Let\/ \(\varGamma\tto M\) be a source-prop\-er Lie groupoid. Let\/ \(U,V\subset M\) be invariant open sets with\/ \(\overline U\subset V\), and let\/ \(\varPhi\in\Mcon(\varGamma\mathbin|V)\) be a multiplicative connection defined over\/ \(V\).

{\rm 0)}\hskip\labelsep One can always extend\/ \(\varPhi\mathbin|U\) over some invariant open neighborhood of\/ \(M_0\cup U\).

{\rm 1)}\hskip\labelsep Suppose, in addition, that\/ \(\varGamma\tto M\) is source connected. Then, one can extend\/ \(\varPhi\mathbin|U\) over some invariant open neighborhood of\/ \(M_0\cup M_1\) containing\/ \(U\).

{\rm 2)}\hskip\labelsep Still further, suppose that the source fibers of\/ \(\varGamma\tto M\) (besides being connected) have finite fundamental groups. Then, \(\varPhi\mathbin|U\) can be extended over some invariant open neighborhood of\/ \(M_0\cup M_1\cup M_2\cup U\) if, and only if, for every base point\/ \(x\in M_2\) the four conditions below, concerning the long exact sequence of homotopy groups associated with the pointed fibration\/ \((\varGamma^x_x,1_x)\xto\subset(\varGamma^x,1_x)\xto{t^x}(O_x,x)\), are simultaneously verified:
\begin{enumerate}
\def\labelenumi{\rm(\alph{enumi})}%
 \item \(\pi_2(\varGamma^x)=0\).
 \item \(F=\im{}\bigl(\free\partial_2:\pi_2(O_x)\to\free\pi_1(\varGamma^x_x)\bigr)\) sits inside\/ \(\free\pi_1(\varGamma^x_x)\) so that\/ \(\Z e\cap F\neq 0\seq 2e\in F\) for all\/ \(e\).
 \item \(\partial_1:\pi_1(O_x)\simto\pi_0(\varGamma^x_x)\).
 \item \(\free\pi_1(c_g)=-\id\in\Aut\bigl(\free\pi_1(\varGamma^x_x)\bigr)\) for\/ \(g\in\varGamma^x_x\smallsetminus\idc{(\varGamma^x_x)}\), where\/ \(c_g\in\Aut(\varGamma^x_x)\) denotes conjugation by\/ \(g\).
\end{enumerate} \end{thm}

\noindent{\it Comment.} When \(\varGamma\tto M\) is source simply connected, the conditions (b) and (c) are of course superfluous.

\begin{proof} We ask the reader to accept it on trust that the general validity of the theorem is a consequence of its validity in the special case that $\varGamma\tto M$ is regular (the reduction to the regular case being a straightforward application of \cite[Theorem~1.8]{2016a}). Alternatively, the reader may introduce an extra hypothesis of regularity on $\varGamma\tto M$ in the statement of the theorem. In any case, we shall throughout the proof assume that $\varGamma\tto M$ is regular.

Let $\rho^\varPhi\in\Gamma^\infty(V;\mathscr R^\varGamma)^P$ denote the equivariant local $C^\infty$ section of $\mathscr R^\varGamma$ with domain $V$ which by \ref{stmt:14A.6.4} corresponds to $\lambda^\varPhi_\varLambda\in\Rep(\varGamma\mathbin|V;\varLambda\mathbin|V)$. On the account of Section~\ref{sec:1}, in order that $\varPhi\mathbin|U$ may be prolonged globally it suffices that there exists some extension of $\rho^\varPhi\mathbin|U$ to a global equivariant $C^\infty$ section of $\mathscr R^\varGamma$. When the rank of $\varGamma\tto M$ is zero (resp.~one), the existence of such an extension is a direct consequence of Example~\ref{exmp:10/21*} (resp.~\ref{exmp:10/21**}). Suppose, on the other hand, that the rank is two.

(\ensuremath\Rightarrow.)\hskip.5em Suppose that \(\rho\in\Gamma^\infty(M;\mathscr R^\varGamma)^P\) exists extending \(\rho^\varPhi\mathbin|U\). By \ref{stmt:14A.6.4}, we must have \(\rho=\rho^\alpha\) for a unique \(\alpha\in\Rep(\varGamma;\varLambda)\). Let us endow \(\varLambda\) with some \(\alpha\)~invariant vec\-tor-bun\-dle metric. For every \(x\in M\) the transitive Lie groupoid \(\varGamma\mathbin|O_x\tto O_x\) satisfies the hypotheses of the propositions \ref{prop:14A.9.1}, \ref{prop:14A.9.2} and \ref{prop:14A.9.4} and admits a realizable orthogonal tangent isotropy representation \(\alpha^x_x:\varGamma^x_x\to\GO(T_xO_x)\). The necessity of the conditions (a)\textendash (d) is then a direct consequence of those propositions and of the remarks preceding Lemma~\ref{lem:14A.9.3}.

(\ensuremath\Leftarrow.)\hskip.5em Let \(\varSigma\tto Z\) be an arbitrary rank-two transitive Lie groupoid which is compact, source connected, and which admits tangent representations. Let us call a class \(r\in R^\varSigma\) \emph{realizable} if it is of the form \(r=[(z,\alpha^z_z)]\) for some \(z\in Z\) and for some \(\alpha:\varSigma\to\GL(TZ)\). For \(z\) fixed, and for any choice of Riemannian metric on \(Z\), every realizable class \(r\in R^\varSigma\) will admit representatives \((z,\alpha^z_z)\) that originate from {\em orthogonal}\/ tangent representations \(\alpha:\varSigma\to\GO(TZ)\). (This is a straightforward application of the orthogonalization trick, cf.~\ref{npar:14A.7.11}/b.) Under the additional assumption that the source fibers of \(\varSigma\tto Z\) have finite fundamental groups, there will be exactly one class in \(R^\varSigma\) which is realizable. We shall denote it \(r^\varSigma\). [By Proposition~\ref{prop:14A.9.2}, for \(\alpha:\varSigma\to\GO(TZ)\) arbitrary, the connected isotropy representation \(\idc{\alpha^z_z}:\idc{\varSigma^z_z}\to\SO(T_zZ)\) is uniquely determined. Its extension \(\alpha^z_z\) to all of \(\varSigma^z_z\) is unique up to within conjugation by the elements of \(\SO(T_zZ)\), by Lemma~\ref{lem:14A.9.3}.] Under the same assumption, for any Lie-group\-oid isomorphism \(\kappa:\varSigma'\simto\varSigma\) we will have \(r^{\varSigma'}=\kappa^*(r^\varSigma)\), the operation \(\kappa^*:R^\varSigma\simto R^{\varSigma'}\) being defined in the obvious way so as to correspond along the canonical maps \eqref{eqn:14A.7.11} to the inverse image operation \(\kappa^*:\mathscr R^\varSigma\simto\mathscr R^{\varSigma'}\).

Back to \(\varGamma\tto M\), observe that for each \(x\in M\) the Lie groupoid \(\varPi^x\tto P^x\) satisfies the hypotheses required of \(\varSigma\tto Z\) in the previous paragraph in order that the class \(r^\varSigma\in R^\varSigma\) may exist unique. Indeed, by \ref{prop:14A.9.2} and \ref{prop:14A.9.4} applied to \(\varGamma\mathbin|O_x\tto O_x\), the properties (a)\textendash (d) imply that \(\Rep(\varGamma\mathbin|O_x;TO_x)\) and, a~fortiori, \(\Rep(\varPi^x;TP_x)\) have to be non-emp\-ty. Let us set \(r^\varGamma(x):=r^{\varPi^x}\in R^{\varPi^x}=R^\varGamma_x\). Then, \(r^\varGamma\in\Gamma^\infty(M;R^\varGamma)^P\). [{\sc proof\emphpunct:} Source properness of $\varGamma\tto M$ implies local linearizability of $\varGamma\tto M$ and hence local triviality of $\varPi\tto P$ along $s:P\to M$ in the sense of \eqref{eqn:14A.7.4}. It will then be possible to cover $M$ with open sets $W$ for which the regular Lie groupoids $\varPi^W\tto P^W$ admit longitudinal tangent representations. Any such representation will correspond to an element of $\Gamma^\infty(W;\mathscr R^\varGamma)$, and the image of any such element under the canonical map \eqref{eqn:14A.7.11} will be $r^\varGamma\mathbin|W$. This shows that $r^\varGamma$ is $C^\infty$. As to its equivariance, viewing every arrow $k\in\varGamma(x,y)$ as a Lie-group\-oid isomorphism $k:\varPi^y\simto\varPi^x$, we have $r^\varGamma(y)=r^{\varPi^y}=k^*(r^{\varPi^x})=k\cdot r^\varGamma(x)$. {\sc q.e.d.}] Moreover, the definition of \(r^\varGamma\) entails that \(\rho^\varPhi\mathbin|U\in\Gamma^\infty(U;\mathscr R^\varGamma[r^\varGamma])^P\). We may then invoke Lemma~\ref{lem:14A.11.6} in order to get the existence of elements of \(\Gamma^\infty(M;\mathscr R^\varGamma)^P\) which extend $\rho^\varPhi\mathbin|U$, since for every $x$ the class $r^\varGamma(x)$ evidently satisfies the conditions (1)\textendash (3) in that lemma. \end{proof}

\begin{npar}\label{npar:2016b.1} Among the four conditions (a)\textendash (d), neither (a) nor (b) is redundant\textemdash i.e.~implied by the remaining three conditions\textemdash as the following examples illustrate.

The pair groupoid over the two-sphere, \(S^2\times S^2\tto S^2\), satisfies all the conditions but (a). That over the real projective plane, \(P^2\times P^2\tto P^2\), violates both (a) and (c).

For any transitive Lie groupoid \(\varSigma\tto Z\), given any closed normal subgroup \(G\) of one isotropy group \(\varSigma^z_z\), there will be a unique kernel \(K^G\subset\varSigma\tto Z\) for which \(K^G\cap\varSigma^z_z=G\). The quotient \(\varSigma/G:=\varSigma/K^G\tto Z\) will be a transitive Lie groupoid, with source fibration \(\varSigma^z_z/G\emto\varSigma^z/G\onto Z\). Now, for concreteness, take \(\varSigma\tto Z\) to be the action groupoid \(\SO(3)\ltimes S^2\tto S^2\), and let \(G\subset\varSigma^z_z\cong\SO(2)\) be an arbitrary non-triv\-i\-al finite subgroup (necessarily cyclic). Let \(\varepsilon:\varSigma\to\varSigma/G\) stand for the projection homomorphism. Then the transformation of long exact sequences of homotopy groups associated with the fi\-ber-bun\-dle map \(\varepsilon^z:\varSigma^z\to\varSigma^z/G\) shows that \(\partial_2^{\varSigma/G}:\pi_2(Z)\to\pi_1(\varSigma^z_z/G)\) equals the composition \(\pi_2(Z)\xto{\partial_2^\varSigma}\pi_1(\varSigma^z_z)\xto{\pi_1(\varepsilon^z_z)}\pi_1(\varSigma^z_z/G)\) and hence that under the identification \(\pi_1(\varSigma^z_z/G)\cong\Z\) the image of \(\partial_2^{\varSigma/G}\) corresponds to \(2|G|\Z\). The condition (b) is thus violated if we take \(\varGamma\tto M\) to be \(\varSigma/G\tto Z\). On the other hand, the condition (a) is satisfied, because \(\varSigma^z\xto{\varepsilon^z}\varSigma^z/G\) is a finite covering map and therefore \(0=\pi_2(\varSigma^z)\cong\pi_2(\varSigma^z/G)\).

Alternatively, we could have taken \(\varSigma\tto Z\) to be the action groupoid \(\SO(3)\ltimes P^2\tto P^2\). \end{npar}

\begin{cor}\label{cor:14A.11.8} Let\/ \(\varGamma\tto M\) be a Lie groupoid which is source proper, source connected, and whose source fibers have finite fundamental groups. Suppose that\/ \(\dim O_x\leqq 2\) for all\/ \(x\in M\), and that for every\/ \(x\in M_2\) the conditions (a)~to (d) of the preceding proposition are verified. Then, the\/ \(C^\infty\)\mdash space\/ \(\Mcon(\varGamma)\) is of the ``weak\/ \(C^\infty\) homotopy type'' of a point, in the sense that it is non-emp\-ty and the relative homotopy existence property below, which implies in particular that\/ \(\Mcon(\varGamma)\) has to be\/ \(C^\infty\)\mdash path-con\-nect\-ed, holds:
\begin{itemize}
 \item[(\textasteriskcentered)] Given any two\/ \(C^\infty\) parametric families of multiplicative connections\/ \(\varPhi_0,\varPhi_1:S\to\Mcon(\varGamma)\) indexed over a smooth manifold\/ \(S\) which coincide in a neighborhood of a given closed subset\/ \(C\subset S\), there exists some\/ \(C^\infty\) homotopy between\/ \(\varPhi_0\) and\/ \(\varPhi_1\) whose restriction to\/ \(C\) stays constant in time.
\end{itemize} \end{cor}

\begin{proof} It will evidently be enough to show that for every open subset \(W\subset\R\times S\) such that \(W\supset(\R\times C)\cup(\R\smallsetminus\mathopen]0,1\mathclose[\times S)\) the restriction of an arbitrary \(C^\infty\) parametric family $W\to\Mcon(\varGamma)$ to a sufficiently small neighborhood of \((\R\times C)\cup(\R\smallsetminus\mathopen]0,1\mathclose[\times S)\) can be prolonged to all of \(\R\times S\). Now, the latter extension property is a direct consequence of Theorem~\ref{prop:14A.11.7} applied to the groupoid \(\varGamma\times\R\times S\tto M\times\R\times S\). \end{proof}

A substantial part of the preceding result remains valid under significantly weaker assumptions.

\begin{thm}\label{prop:14A.11.9} Let\/ \(\varGamma\tto M\) be a Lie groupoid which is source proper, source connected, and whose orbit space\/ \(M/\varGamma\) is connected. Suppose that\/ \(\dim O_x\leqq 2\) for all\/ \(x\in M\), and moreover that there is some\/ \(x\) for which\/ \(\dim O_x\leqq 1\). Then, providing that it is non-emp\-ty, the\/ \(C^\infty\)\mdash space\/ \(\Mcon(\varGamma)\) is of the ``weak\/ \(C^\infty\) homotopy type'' of a point, in the sense that the relative homotopy existence property (\textasteriskcentered) from the preceding corollary holds. In particular, \(\Mcon(\varGamma)\) is\/ \(C^\infty\)\mdash path-con\-nect\-ed. \end{thm}

\begin{proof} By reasoning as in the proof of the corollary, and by the rank~$\leqq 1$ case of Theorem~\ref{prop:14A.11.7}, we are reduced to proving an assertion of the following form:

\begin{lem*} Let\/ \(\varGamma\tto M\) be as in the statement of the theorem. Let\/ \(U,V\subset M\) be invariant open neighborhoods of\/ \(M_{\leq 1}:=M_0\cup M_1\) with\/ \(\overline U\subset V\) such that every invariant component of\/ $V$ meets\/ $M_{\leq 1}$ i.e.~every connected component of\/ $V/\varGamma\subset M/\varGamma$ meets\/ $M_{\leq 1}/\varGamma$. Then the restriction over\/ $U$ of any local multiplicative connection\/ $\varPhi\in\Mcon(\varGamma\mathbin|V)$ admits some extension to a multiplicative connection defined on the whole\/ $\varGamma\tto M$. \end{lem*}

Before turning to the proof, we need to discuss the following general construction. Let $\mathscr R\to X$ be an arbitrary $\varOmega\tto X$ equivariant $C^\infty$\mdash fibration (cf.~\S\ref{sec:2}). Under the hypothesis that $\mathscr R\to X$ be locally trivial, there is a ``finest'' discrete obstruction \(\mathscr R\to C^{\mathscr R}\) for \(\mathscr R\to X\), which we call the \emph{component obstruction}, which is ``universal'' in the sense that any other discrete obstruction factors through it. Namely, let us set $C^{\mathscr R}:=\bigcup_{x\in X}\pi_0(\mathscr R_x)$, where $\pi_0(\mathscr R_x)$ denotes the set of ($C^\infty$\mdash path\babelhyphen{nobreak})con\-nect\-ed components of the fiber $\mathscr R_x$. Let us endow $C^{\mathscr R}$ with the quotient $C^\infty$\mdash structure induced along the natural projection $\mathscr R\to C^{\mathscr R}$. Let us make $\varOmega\tto X$ act along $C^{\mathscr R}\to X$ so that the natural projection becomes equivariant. Since $\mathscr R\to X$ admits local $C^\infty$ sections through every point, $C^{\mathscr R}\to X$ must be an \'etale differentiable mapping. The local triviality of $\mathscr R\to X$ implies that $C^{\mathscr R}\to X$ must be a covering mapping of Hausdorff differentiable manifolds\textemdash in particular, a discrete obstruction. When $\mathscr R=\mathscr R^\varGamma$, we may of course abbreviate $C^{\mathscr R^\varGamma}$ into $C^\varGamma$.

We are now ready to prove the lemma. Since $M_{\leq 1}$ is closed within $M$, its complement $M_2=M\smallsetminus M_{\leq 1}$ is an invariant open subset of $M$. The restriction $\varGamma\mathbin|M_2\tto M_2$ is a regular Lie groupoid of rank two which is source proper and source connected. Let us consider the \(\varGamma\mathbin|M_2\tto M_2\)~equivariant \(C^\infty\)\mdash fibration \(\mathscr R:=\mathscr R^{\varGamma\mathbin|M_2}\to M_2\), along with the associated component obstruction \(C:=C^{\varGamma\mathbin|M_2}\to M_2\). Let us put \(U_2:=U\cap M_2\emphpunct, V_2:=V\cap M_2\).

Since $\Mcon(\varGamma)\neq\emptyset$, there exists some multiplicative connection on $\varGamma\tto M$, say, $\varPsi$. Let $\rho^{\varPsi\mathbin|M_2}\in\Gamma^\infty(M_2;\mathscr R)^{\varGamma\mathbin|M_2}\emphpunct, c^{\varPsi\mathbin|M_2}\in\Gamma^\infty(M_2;C)^{\varGamma\mathbin|M_2}$ denote the corresponding equivariant global $C^\infty$ sections of $\mathscr R,~C$. By construction, the invariant subbundle $\mathscr R[c^{\varPsi\mathbin|M_2}]\subset\mathscr R$ has ($C^\infty$\mdash path\babelhyphen{nobreak})con\-nect\-ed fibers $(\mathscr R[c^{\varPsi\mathbin|M_2}])_x=c^{\varPsi\mathbin|M_2}(x)$. Since $\mathscr R$ is locally transversely trivializable, so must be $\mathscr R[c^{\varPsi\mathbin|M_2}]\to M_2$, which then by \ref{stmt:14A.11.1} satisfies the hypotheses of Proposition~\ref{prop:14A.7.5}. If we can show that $\rho^{\varPhi\mathbin|V_2}$ takes values in $\mathscr R[c^{\varPsi\mathbin|M_2}]$, Proposition~\ref{prop:14A.7.5} and \S\ref{sec:1} are going to guarantee the existence of a multiplicative connection on $\varGamma\mathbin|M_2\tto M_2$ that extends $\varPhi\mathbin|U_2$, which we can then paste to $\varPhi\mathbin|U$ so as to get the required global extension of $\varPhi\mathbin|U$. Now, by Proposition~\ref{prop:14A.11.7}, over a suitable invariant open neighborhood say $W$ of $M_{\leq 1}$ with $\overline W\subset V$ there has to be some $C^\infty$ homotopy between $\varPhi\mathbin|W$ and $\varPsi\mathbin|W$. Thus $\rho^\varPhi\mathbin|W_2\in\Gamma^\infty(W_2;\mathscr R[c^{\varPsi\mathbin|M_2}])^{\varGamma\mathbin|M_2}$, where $W_2:=W\cap M_2$. It follows that $c^{\varPhi\mathbin|M_2}$ and $c^{\varPsi\mathbin|M_2}$ agree on all of $U_2$, because they agree on $W_2$, and $W_2$ meets each invariant component of $U_2$. This proves the lemma and, hence, the theorem. \end{proof}

\appendix

\section{Technical complements}\label{sec:A}

\paragraph*{I.} We begin by reviewing some general notions relating to \(C^\infty\)\mdash spaces (cf.~\S\ref{npar:14A.6.3}) and then prove an auxiliary result for use in Section~\ref{sec:4}, more precisely in the proof of \ref{stmt:14A.8.6}.

Recall that the \emph{subspace\/ \(C^\infty\)\mdash structure} on an arbitrary subset \(\mathscr Z\) of a given \(C^\infty\)\mdash space \(\mathscr X\) consists by definition of those \(C^\infty\) maps \(S\to\mathscr X\) whose images are contained in \(\mathscr Z\). By a \emph{\(C^\infty\)\mdash embedding} of another \(C^\infty\)\mdash space \(\mathscr Y\) into \(\mathscr X\), we mean an injective map \(e:\mathscr Y\emto\mathscr X\) with the property that the \(C^\infty\) maps \(S\to\mathscr Y\) are those which upon composition with \(e\) produce \(C^\infty\) maps \(S\to\mathscr X\), equivalently, one which induces a \(C^\infty\)\mdash isomorphism between \(\mathscr Y\) and the \(C^\infty\)\mdash subspace $e(\mathscr Y)\subset\mathscr X$.

Every \(C^\infty\)\mdash space \(\mathscr X\) has a natural topology, the \emph{\(C^\infty\)\mdash space topology}; by definition, this is the finest topology on \(\mathscr X\) making \(C^\infty\) maps from smooth manifolds into \(\mathscr X\) continuous. Clearly, \(C^\infty\) maps between \(C^\infty\)\mdash spaces are continuous with respect to this topology. The \(C^\infty\)\mdash path-com\-po\-nents of a \(C^\infty\)\mdash space \(\mathscr X\) are evidently open (thus, closed) relative to the \(C^\infty\)\mdash space topology on \(\mathscr X\), hence they coincide with the (ordinary topological) connected components of \(\mathscr X\). Our notation $\pi_0(\mathscr X)$ for the set of \(C^\infty\)\mdash path-com\-po\-nents of \(\mathscr X\), which we used freely in Sections \ref{sec:5}~and \ref{sec:6}, is therefore legitimate.

For any equivalence relation \(\mathscr E\subset\mathscr X\times\mathscr X\) on a given \(C^\infty\)\mdash space \(\mathscr X\), the \emph{induced\/ \(C^\infty\)\mdash structure} on the quotient \(\mathscr X/\mathscr E\) consists by definition of those maps \(S\to\mathscr X/\mathscr E\) which can be lifted locally to \(C^\infty\) maps of \(S\) into \(\mathscr X\). A map out of \(\mathscr X/\mathscr E\) will be \(C^\infty\) iff so is its composition with the quotient projection \(\mathscr X\to\mathscr X/\mathscr E\).

The \emph{mapping space} \(\mathscr M(\mathscr X,\mathscr Y)\) for a pair of \(C^\infty\)\mdash spaces \(\mathscr X\),~\(\mathscr Y\) is the \(C^\infty\)\mdash space obtained by endowing the set of \(C^\infty\) maps from \(\mathscr X\) into \(\mathscr Y\) with the \(C^\infty\)\mdash structure relative to which the \(C^\infty\) maps \(S\to\mathscr M(\mathscr X,\mathscr Y)\) are those that correspond to \(C^\infty\) maps \(S\times\mathscr X\to\mathscr Y\).

\begin{lem}\label{lem:A.1} Let\/ $G$ and\/ $H$ be compact Lie groups. The orbits of the standard conjugation action of\/ $\idc H$ on\/ $\Hom(G,H)$ are precisely the connected components of\/ $\Hom(G,H)$. The\/ \(C^\infty\)\mdash structure on\/ $\Hom(G,H)$ coincides with the unique differentiable manifold structure which makes every map\/ $\idc H\to\Hom(G,H)\emphpunct, h\mapsto hfh^{-1}$ submersive. \end{lem}

\begin{proof} Our proof will be based on the following remark\emphpunct: \em Let\/ \(g_1,\dotsc,g_n\) be representatives for the connected components of\/ \(G\smallsetminus\idc G\). The map below is an\/ \(H\)~equivariant\/ \(C^\infty\)\mdash embedding,
\begin{equation}
\label{eqn:2016b.3}
	\Hom(G,H)\longto L(\mathfrak g,\mathfrak h)\times H\times\dotsb\times H,
\quad	f\mapsto\bigl(\mathrm Lf,f(g_1),\dotsc,f(g_n)\bigr)
\end{equation}
the action\/ \(H\circlearrowright L(\mathfrak g,\mathfrak h)\times H^{\times n}\) being given by\/ $\bigl(h;\lambda,(h_i)_{i=1}^n\bigr)\mapsto\bigl(\Ad_H(h)\circ\lambda\mathord,\mskip\medmuskip (hh_ih^{-1})_{i=1}^n\bigr)$. \em Proof\emphpunct: The map in question is equivariant because of the naturality of $\exp$,
\begin{flalign*}
&&	\begin{split}
\xymatrix@C=3em{%
	\mathfrak g
	\ar[d]^{\exp_G}
	\ar[r]^{\mathrm Lf}
	&	\mathfrak h
		\ar[d]^{\exp_H}
		\ar[r]^{\Ad_H(h)}
		&	\mathfrak h
			\ar[d]^{\exp_H}
\\	G
	\ar[r]^f
	&	H
		\ar[r]^{c_h}
		&	H
}\end{split}
&&	\makebox[.em][r]{[$c_h(x)=hxh^{-1}$].}
\end{flalign*}
The map is injective, because a Lie-group homomorphism $f$ of $G$ into $H$ is determined by its values $f(g_1),\dotsc,f(g_n)$ and its restriction to the identity component of \(G\), which in turn is determined by $\mathrm Lf$ (by basic Lie theory). The map is also clearly $C^\infty$. The only thing left for us to check is that any smooth map with values in the image of \eqref{eqn:2016b.3}, say,
\begin{equation*}
	S\xto{(\lambda,(h_i)_{i=1}^n)}L(\mathfrak g,\mathfrak h)\times H^{\times n},
\quad	s\mapsto\left(\lambda(s),\bigl(h_i(s)\bigr)_{i=1}^n\right)=\left(\mathrm Lf(s),\bigl(f(s)(g_i)\bigr)_{i=1}^n\right)
\end{equation*}
gives rise to a smooth map $S\times G\to H\emphpunct, (s,g)\mapsto f(s)(g)$. Evidently, it will be enough to show that $S\xto\lambda L(\mathfrak g,\mathfrak h)$ gives rise to a smooth map $S\times\idc G\to H$. Without loss of generality, we may assume that $G$ is connected. Let $B$ be any open ball centered at $0\in\mathfrak g$ so small that $\exp_G$ restricts to a diffeomorphism of $B$ onto some open neighborhood $U\subset G$ of $1_G$. Since the exponential map of a connected compact Lie group is onto, and because of the naturality of the exponential map, we see that for every $g=\exp_G(X)\in G$ the restriction
\begin{align*}
	S\times gU\to H,
\quad	(s,gu)\mapsto f(s)(gu)	&=f(s)(g)f(s)(u)
\\				&=f(s)\bigl(\exp_G(X)\bigr)f(s)\bigl(\exp_G(\exp_G^{-1}(u))\bigr)
\\				&=\exp_H\bigl(\mathrm Lf(s)X\bigr)\exp_H\bigl(\mathrm Lf(s)\exp_G^{-1}(u)\bigr)
\\				&=\exp_H\bigl(\lambda(s)X\bigr)\exp_H\bigl(\lambda(s)\exp_G^{-1}(u)\bigr)
\end{align*}
is smooth. Now the open sets $gU$ cover the whole $G$, so the desired conclusion follows. Our remark is thus proven.

It is well known that the connected components of \(\Hom(G,H)\) for the com\-pact-open topology coincide with the \(\idc H\)~orbits \(\mathscr O_f:=\{c_h\circ f:h\in\idc H\}\) for the standard conjugation action of \(\idc H\) on \(\Hom(G,H)\). (For instance, compare \cite[Lemma~38.1]{CF}, or \cite{LW}.) Since the \(C^\infty\)\mdash space topology on \(\Hom(G,H)\) is finer than the com\-pact-open topology, these must be closed open subsets of \(\Hom(G,H)\) for the \(C^\infty\)\mdash space topology. The equivariant \(C^\infty\)\mdash embedding \eqref{eqn:2016b.3} restricts to a \(C^\infty\)\mdash isomorphism \(\mathscr O_f\approxto O_f\) between each \(\idc H\)~orbit \(\mathscr O_f\subset\Hom(G,H)\) and the corresponding orbit in \(L(\mathfrak g,\mathfrak h)\times H^{\times n}\), \[%
	O_f:=\left\{\mkern 1mu\left(\Ad_H(h)\circ\mathrm Lf\mathord,\mskip\medmuskip\bigl(hf(g_i)h^{-1}\bigr)_{i=1}^n\right):h\in\idc H\mkern 1mu\right\}.
\] Now for \(\idc H\) compact the orbits of the action \(\idc H\circlearrowright L(\mathfrak g,\mathfrak h)\times H^{\times n}\) are smooth submanifolds of \(L(\mathfrak g,\mathfrak h)\times H^{\times n}\), so \(\Hom(G,H)\) admits an open cover by ``charts'' and must therefore be a differentiable manifold. \end{proof}

\paragraph*{II.} We proceed to establish another auxiliary result, which we need in \S\ref{sec:5} in the proof of Propositions \ref{prop:14A.9.2}~and \ref{prop:14A.9.4}. Let us begin with the following remark:

\begin{lem*} Let\/ $G$ be a connected Lie group. Let\/ $\vartheta:G\to\mathbb T$~(= one-to\-rus\/ $\R/\Z$) be any Lie-group homomorphism whose associated homomorphism of fundamental groups\/ $\pi_1(\vartheta):\pi_1(G)\to\pi_1(\mathbb T)$ is zero. Then, there exists a\/ $C^\infty$ path\/ $\R\ni t\mapsto\vartheta_t\in\Hom(G,\mathbb T)$ joining\/ $\vartheta_1=\vartheta$ and the constant homomorphism $\vartheta_0=1_G$. \end{lem*}

\begin{proof} Consider the following lifting problem: \[%
\xymatrix@C=3em{%
 \ast \ar[r]^0
 \ar[d]_{1_G}
 &	\R
	\ar[d]^\pr
\\ G \ar[r]^\vartheta
 \ar@{.>}[ur]^{\tilde\vartheta}
 &	\mathbb T
}\] Since by assumption $\im\pi_1(\vartheta)\subset\{0\}=\im\pi_1(\pr)$, the covering map lemma implies that there has to be a unique solution $\tilde\vartheta$ to this problem, which is also necessarily a Lie-group homomorphism of $G$ into $\R=(\R,+)$. Then, setting \[%
	\vartheta_t(g):=\pr\bigl(t\tilde\vartheta(g)\bigr)
\] defines the required homotopy from $\vartheta$ to $1_G$ through homomorphisms $G\to\mathbb T$. \end{proof}

By the general structure theory of compact Lie groups (see e.g.~\cite[\S V.8]{BtomD}) every compact, connected Lie group $G$ fits in a short exact sequence
\begin{equation}
\label{eqn:2016b.4}
	1\to A\to\mathbb T^r\times K\to G\to 1
\end{equation}
where $A\subset Z(\mathbb T^r\times K)$ is a finite central subgroup of the product of a torus $\mathbb T^r=(\R/\Z)^r$ and a compact, (connected and) simply connected Lie group $K$. Note that $\pi_1(G)$ must be finitely generated and precisely of rank $r$ since it fits in a short exact sequence $1\to\Z^r=\pi_1(\mathbb T^r\times K)\to\pi_1(G)\to\pi_0(A)=A\to 1$ with $A$ finite. Letting $T\subset A$ denote the kernel of the composite homomorphism $A\xto\subset\mathbb T^r\times K\xto\pr\mathbb T^r$, which we may also regard as a central subgroup $T$ of $K$, we get an alternative presentation of the form \eqref{eqn:2016b.4} where, now, $A$ is a finite subgroup of $Z(\mathbb T^r\times K)$ for which the composite homomorphism $A\xto\subset\mathbb T^r\times K\xto\pr\mathbb T^r$ is injective, and $K$ is a compact, connected Lie group with finite fundamental group $T$.

\begin{lem}\label{lem:14A.9.8} For any connected Lie group\/ $G$ and any compact, connected, abelian Lie group\/ $T$, the homomorphism
\begin{equation}
\label{eqn:2016b.5}
	\free\pi_1:\Hom(G,T)\longto\Hom{}\bigl(\free\pi_1(G)\mathord,\mskip\thickmuskip\free\pi_1(T)\bigr),
\quad	\vartheta\mapsto\free\pi_1(\vartheta)
\end{equation}
has path-con\-nect\-ed fibers\textemdash in other words, any two Lie group homomorphisms\/ $\vartheta_0,\vartheta_1\in\Hom(G,T)$ will be\/ $C^\infty$\mdash homotopic through Lie group homomorphisms if, and only if, $\free\pi_1(\vartheta_0)=\free\pi_1(\vartheta_1)$. If in addition\/ $G$ is compact, then the same homomorphism is bijective. \end{lem}

\begin{proof} After composition with any set of projections $T\cong\mathbb T^n\xto{\pr_i}\mathbb T\sidetext(i=1,\dotsc,n)$, we are reduced to the special case $T=\mathbb T$. Since $\mathbb T$ is abelian, the ratio $\vartheta_1{\vartheta_0}^{-1}$ of two Lie group homomorphisms $G\to\mathbb T$ is itself one such homomorphism, and we have $\pi_1(\vartheta_1{\vartheta_0}^{-1})=\pi_1(\vartheta_1)-\pi_1(\vartheta_0)$. The first claim is then an immediate consequence of the last lemma.

As to the second claim, the injectivity of \eqref{eqn:2016b.5} for compact $G$ is clear since any homotopy through Lie group homomorphisms from a compact Lie group into an abelian Lie group is constant. Regarding the surjectivity, let $\varphi:\free\pi_1(G)\to\free\pi_1(\mathbb T)=\Z$ be a random homomorphism of groups. When $G=\mathbb T^r$ for some $r\geqq 0$, the existence of $\vartheta:G\to\mathbb T$ for which $\free\pi_1(\vartheta)=\varphi$ is obvious: \[%
	\vartheta(\zeta_1,\dotsc,\zeta_r):=\zeta_1^{\varphi(e_1)}\dotsm\zeta_r^{\varphi(e_r)}
\] where $e_i$ denotes the $i$\mdash th standard basis vector in $\Z^r=\pi_1(\mathbb T^r)=\free\pi_1(\mathbb T^r)$. Same story in the case when $G=\mathbb T^r/\Lambda$ for some finite subgroup $\Lambda\subset\mathbb T^r$, for then there have to be Lie group isomorphisms $G\cong\mathbb T^r$ (as $\mathbb T^r/\Lambda$ is an $r$\mdash dimensional, compact, connected, abelian Lie group, hence an $r$\mdash torus). For general $G$, in the notations of \eqref{eqn:2016b.4} fwd., we have the following homomorphisms of short exact sequences of Lie groups,
\begin{equation*}
 \newdir{ >}{{}*!/-5pt/@{>}}%
\xymatrix@C=1em{%
	1 \ar[r]
	&	A \ar@{ >->}[r]
		\ar@{->>}[d]
		&	\mathbb T^r\times K \ar@{->>}[r]
			\ar@{->>}[d]
			&	G \ar[r]
				\ar@{=}[d]
				&	1
\\	1 \ar[r]
	&	A/T \ar@{ >->}[r]
		\ar[d]^\simeq
		&	\mathbb T^r\times K/T \ar@{->>}[r]
			\ar[d]^\pr
			&	G \ar[r]
				\ar@{-->}[d]^\tau
				&	1
\\	1 \ar[r]
	&	\pr(A) \ar@{ >->}[r]
		&	\mathbb T^r \ar@{->>}[r]
			&	\mathbb T^r/\pr(A) \ar[r]
				&	1}
\end{equation*}
which induce transformations between the associated long exact sequences of homotopy groups,
\begin{equation*}
 \newdir{ >}{{}*!/-5pt/@{>}}%
\xymatrix@C=1em{%
	1 \ar[r]
	&	\Z^r \ar@{ >->}[r]
		\ar[d]
		&	\pi_1(G) \ar@{->>}[r]
			\ar@{=}[d]
			&	A \ar[r]
				\ar@{->>}[d]
				&	1
\\	1 \ar[r]
	&	\Z^r\oplus T \ar@{ >->}[r]^-{j'}
		\ar[d]
		&	\pi_1(G) \ar@{->>}[r]
			\ar[d]^{\pi_1(\tau)}
			&	A/T \ar[r]
				\ar[d]^\simeq
				&	1
\\	1 \ar[r]
	&	\Z^r \ar@{ >->}[r]
		&	\pi_1\bigl(\mathbb T^r/\pr(A)\bigr) \ar@{->>}[r]
			&	\pr(A) \ar[r]
				&	1\makebox[.em][l]{.}}
\end{equation*}
Notice the inclusions $\ker\bigl(\pi_1(\tau)\bigr)\subset j'(0\oplus T)\subset\tor\pi_1(G)$, the last of which is valid because $T$ is finite. In fact, since $\pi_1\bigl(\mathbb T^r/\pr(A)\bigr)$ is free,
\begin{equation*}
	\ker\bigl(\pi_1(\tau)\bigr)=\tor\pi_1(G)
\end{equation*}
so that $\free\pi_1(\tau)$ is an isomorphism of $\free\pi_1(G)$ onto all of $\pi_1\bigl(\mathbb T^r/\pr(A)\bigr)$. Now, let \[\bar\varphi:\pi_1\bigl(\mathbb T^r/\pr(A)\bigr)\to\Z\] denote the unique homomorphism such that $\bar\varphi\circ\free\pi_1(\tau)=\varphi$. By the above, we know that there has to be some Lie-group homomorphism $\bar\vartheta:\mathbb T^r/\pr(A)\to\mathbb T$ for which $\pi_1(\bar\vartheta)=\bar\varphi$. Then, setting $\vartheta:=\bar\vartheta\circ\tau$, as desired, we have
\begin{equation*}
	\free\pi_1(\vartheta)=\pi_1(\bar\vartheta)\circ\free\pi_1(\tau)=\bar\varphi\circ\free\pi_1(\tau)=\varphi.	\qedhere%
\end{equation*} \end{proof}

{\footnotesize
\bibliographystyle{abbrv}
\bibliography{bib/gtrentin,bib/2016a}
}%
\end{document}